\documentclass[12pt,reqno]{amsart}

\usepackage{amsthm,amsmath,amssymb,amsfonts}
\usepackage{amsaddr}
\usepackage{natbib}
\usepackage{txfonts}
\usepackage{fullpage}
\newtheorem{theorem}{Theorem}

\newtheorem{proposition}{Proposition}

\newtheorem{lemma}{Lemma}

\newtheorem{example}{Example}%
\newtheorem{remark}{Remark}%

\newtheorem{assumption}{Assumption}

\usepackage[utf8]{inputenc} 
\usepackage[T1]{fontenc}
\usepackage[setpagesize=false]{hyperref}
\usepackage{url}
\usepackage{booktabs}
\usepackage{amsfonts}
\usepackage{nicefrac}
\usepackage{xcolor}
\usepackage{adjustbox}
\usepackage{natbib}

\usepackage{amsmath,amssymb,amsthm,bm}
\usepackage{graphicx}
\graphicspath{ {./img/} }
\usepackage[subrefformat=parens]{subcaption}
\usepackage{mathtools}
\mathtoolsset{showonlyrefs=true}

\newcommand{\indep}{\perp \!\!\! \perp}

\newcommand{\abs}[1]{\left\lvert#1\right\rvert}
\newcommand{\norm}[1]{\left\lVert#1\right\rVert}

\newcommand{\rbr}[1]{\left(#1\right)}
\newcommand{\sbr}[1]{\left[#1\right]}
\newcommand{\cbr}[1]{\left\{#1\right\}}

\newcommand{\N}{\mathbb{N}}

\newcommand{\R}{\mathbb{R}}

\newcommand{\E}{\mathbb{E}}
\newcommand{\pr}{\mathbb{P}}

\newcommand{\iid}{\overset{\mathrm{iid}}{\sim}}
\newcommand{\cvas}{\overset{\mathrm{a.s.}}{\longrightarrow}}

\newcommand{\prox}{\mathrm{prox}}

\DeclareMathOperator*{\argmax}{arg max}
\DeclareMathOperator*{\argmin}{arg min}

\usepackage{xcolor}
\newcommand{\rc}{\color{black}}

\newcommand{\bc}{\color{black}}
\newcommand{\cc}{\color{black}}
\newcommand{\ccc}{\color{black}}

\renewcommand{\hat}{\widehat}
\renewcommand{\d}{\mathrm{d}}

\def\cL{\mathcal{L}}
\def\bzero{\bm{0}}
\def\bb{\bm{b}}
\def\bD{\bm{D}}
\def\bx{\bm{x}}
\def\bI{\bm{I}}
\def\bL{\bm{L}}
\def\bP{\bm{P}}
\def\bU{\bm{U}}
\def\bX{\bm{X}}
\def\bY{\bm{Y}}
\def\bZ{\bm{Z}}
\def\btheta{\bm{\theta}}
\def\bbeta{\bm{\beta}}
\def\bxi{\bm{\xi}}
\def\bvarepsilon{\bm{\varepsilon}}
\def\bSigma{\bm{\Sigma}}
\newcommand{\as}{\overset{\mathrm{a.s.}}{\longrightarrow}}
\newcommand{\pconv}{\overset{{p}}{\to}}
\newcommand{\dconv}{\overset{{d}}{\to}}
\def\d{\mathrm{d}}

\begin{document}

\title{{\Large Moment-Based Adjustments of Statistical Inference in High-Dimensional Generalized Linear Models}}

\date{}

\allowdisplaybreaks
\sloppy
\raggedbottom

\let\origmaketitle\maketitle
\def\maketitle{
  \begingroup
  \def\uppercasenonmath##1{} % this disables uppercasing title
  \let\MakeUppercase\relax % this disables uppercasing authors
  \origmaketitle
  \endgroup
}

\author{Kazuma Sawaya$^1$, Yoshimasa Uematsu$^2$, Masaaki Imaizumi$^{1,3,4}$
}
\address{
$^1$The University of Tokyo, Tokyo, Japan \\
$^2$Hitotsubashi University, Tokyo, Japan\\
$^3$RIKEN Advanced Intelligene Project, Tokyo, Japan\\
$^4$Kyoto University, Kyoto, Japan
}

\begin{abstract}
We develop a statistical inference method for generalized linear models (GLMs) in high-dimensional settings, where the number of unknown coefficients $p$ is of the same order as the sample size $n$. In this regime, constructing confidence intervals requires estimating unknown hyperparameters, such as the signal strength. However, existing estimators for the hyperparameters are not stably applicable to GLMs when $p/n$ is close to or greater than $1$, both theoretically and empirically. In this study, we develop an estimator for the hyperparameter that addresses the issue and establish an inferential framework, provided that the link function of the GLM exhibits an asymmetry property. The proposed estimator utilizes the moments of the output variable of GLMs and a convex surrogate loss. Our framework is theoretically valid even when the limit of $p/n$ exceeds $1$, ensuring the strong consistency of the hyperparameter estimator and asymptotically attaining the exact coverage probability of the confidence intervals. Our numerical experiments support these theoretical results. 
\end{abstract}

\maketitle

\section{Introduction}
We consider a pair $(\bX,Y)$ with $p$-dimensional random feature $\bX$ and random response $Y$ following the generalized linear model (GLM):
\begin{align}
\label{eq:model}
    \E[Y\mid\bX=\bx]=g\rbr{\bx^\top\bbeta}, ~~~ \forall \bx \in \mathbb{R}^p,
\end{align}
where $g:\R \to \R$ is an inverse link function that {\cc is monotonically increasing}, and $\bbeta = (\beta_1,\dots,\beta_p)^\top \in\R^p$ is an unknown deterministic coefficient vector. 
This is one of the most popular classes of statistical models, including linear regression model, logistic regression model, Poisson regression model, and others. Given a sample of size $n$ and assuming $n$-dependent dimension $p\equiv p(n)$, a GLM and its specific models in the \textit{proportionally high-dimensional regime}, where there exists $\kappa>0$ such that
\begin{align}
    n,p(n) \to \infty \mbox{~~~and~~~} p(n)/n \to \kappa\in(0,\infty), \label{eq:high_dim}
\end{align}
have recently been studied \citep{bayati2011dynamics,rangan2011generalized,el2013robust,thrampoulidis2018precise,sur2019modern,barbier2019optimal,bellec2025observable}.

{\rc 
Our interest lies in \textit{statistical inference} (or simply \textit{inference}), such as hypothesis testing and constructing confidence intervals, for $\bbeta$.
In high-dimensional settings \eqref{eq:high_dim}, a maximum likelihood estimator (MLE) is influenced by dimensionality $\kappa$ and unknown parameters in a complex manner \cite{sur2019likelihood,sur2019modern,salehi2019impact,barbier2019optimal,feng2022unifying}, which is contrast to low-dimensional settings \cite{mccullagh1980regression,cordeiro1983improved} or other (sparse) high-dimensional settings \cite{van2014asymptotically,tian2017asymptotics,fei2021estimation,cai2023statistical,ponnet2024robust}.
Indeed, the asymptotic bias and variance can be characterized as fixed points of a nonlinear system called the \textit{state evolution (SE) system}. 
The several unknown hyperparameters in the system must be estimated for inference, which is a major challenge for practical use.
}

Despite the development of such inferential frameworks, some more challenges remain. 
First, estimating the SE hyperparameters becomes hard as the dimensionality increases, particularly when $\kappa \geq 1$. Existing estimators, including \cite{sur2019modern,yadlowsky2021sloe,chen2024method}, have been theoretically justified only for cases where $\kappa < 1$. Moreover, some methods exhibit numerical instability even when $\kappa$ is less than $1$ but close to $1$, as we verify in Section \ref{sec:numerical-exprmnt}.
Second, even if the SE hyperparameters were known, the inference method with a high-dimensional MLE proposed by \cite{sur2019modern} would not be applicable to certain GLM settings. This is because, in some GLMs, the negative log-likelihood function is non-convex, rendering existing methods relying on convexity inapplicable.

In this study, we develop a \textit{moment-based adjustment} for statistical inference in high-dimensional GLMs and prove its asymptotic validity.
This methodology is based on the following techniques.
First and most importantly, we construct a simple yet effective estimator of the SE hyperparameter, which can be computed from the data by leveraging a moment equation on the output variable $Y$.
{\rc This method enables the estimation of the SE hyperparameters regardless of the value of $\kappa > 0$, while requiring the link function $g(\cdot)$ to be \textit{asymmetric}, i.e., $g_0(x):=(g(x)+g(-x))/2$ being strictly monotone.}
The moment-based estimator remains applicable even when $\kappa \geq 1$, as it neither depends on the inverse of sample covariance matrices nor on the estimator of $\bbeta$.
The second technique is the use of {\cc a} \textit{surrogate} loss function, which is a convex and consistent loss function for GLMs.
This loss function mitigates the non-convexity of the negative log-likelihood function in a certain class of GLMs.
We leverage this loss and derive an SE system for GLMs by applying the approximate message passing (AMP) algorithm.

Our contributions are summarized as follows:
\setlength{\leftmargini}{17pt} 
\begin{itemize}
  \setlength{\parskip}{0cm}
  \setlength{\itemsep}{0cm}
    \item We perform statistical inference for GLMs in a high-dimensional setting.
    We also prove its asymptotic validity; specifically, the proposed confidence interval asymptotically attains the exact coverage probability. To this aim, we leverage the convex surrogate loss function and derive the SE system for a general class of GLMs from the AMP algorithm. 
    \item We develop a moment-based estimator for hyperparameters in the SE system using an additive form of GLMs.
    We prove the consistency of the estimator for any $\kappa\in(0,\infty)$ under the restriction that the link function is asymmetric.
    Additionally, we verify that the estimator remains numerically stable even when $\kappa$ is large.
\end{itemize}

\subsection{Related Works and Comparison}

Statistical estimation and inference in the high-dimensional regime \eqref{eq:high_dim} have been actively studied.
For example, linear models have been considered by \cite{bayati2011dynamics,el2013robust,thrampoulidis2015lasso,thrampoulidis2018precise,takahashi2018statistical,el2018impact,miolane2021distribution,hastie2022surprises,mei2022generalization,bellec2025error}; {\cc while logistic regression models were studied by  }\cite{sur2019likelihood,sur2019modern}, and others \citep{tan2023multinomial,bellec2025observable,loureiro2021learning,sawaya2024high}. In particular, 
{\rc
\cite{barbier2019optimal} deduced the Bayes-optimal estimation and generalization errors for GLMs.
Additionally, \cite{loureiro2021learning} provided a rigorous formula for the asymptotic training loss and generalization error in more general teacher-student settings.
The relationship between our study and these two papers can be summarized as follows.
First, our primary focus is on estimating unknown hyperparameters essential for statistical inference, rather than identifying the risk, which differentiates our motivation.
Second, our proposed estimator with the convex surrogate loss is applicable to a wider class of MLEs, while 
its theoretical foundations 
have already been established in the prior research.
}

Regarding the estimation of SE hyperparameters, the first heuristic approach called \textit{ProbeFrontier} was proposed by \cite{sur2019modern} though it is limited to logistic regression. This method estimates the signal strength (SS) appearing in the SE systems by utilizing the phase transition properties of the MLEs in logistic regression \cite{candes2020phase}. \cite{yadlowsky2021sloe} developed the Signal Strength Leave-One-Out Estimator (SLOE) for the SS by adopting a different representation of the SE system for logistic regression.  
{\rc
It is based on the non-regularized MLE, and therefore the estimation errors can diverge as $\kappa$ increases and approaches to one.
A subsequent paper \cite{chen2024method}, appearing after the preprint release of our study, developed an estimator named Method-of-Moment Inference (MoMI), by extending our moment-based {\cc inference by} incorporating information from $\bX$. 
This study also studied asymptotic variance of the coefficient estimator via the bootstrap method. 
A key difference from our results is that their method {\ccc is difficult to} handle the $ \kappa\geq 1 $ setting when the covariance matrix of $\bX$ is unknown.
{\ccc For the regime $\kappa\geq 1$, they construct an estimator based on Chebyshev-polynomial approximations and higher-order $U$-statistics.
This estimator is consistent in principle, but computationally challenging in practice, and the convergence is expected to be substantially slower than the standard $n^{-1/2}$ rate.
Thus, they are studying the limits of what can be learned from data in principle for the $\kappa \geq 1$ regime, the so-called information-theoretic perspective, without considering computational cost or statistical efficiency.
In contrast, our estimator offers a computationally simple estimator for the signal strength parameter in the regime, and also we study its statistical property with the cost of an additional condition on the link function.
}
 We also note that \citep{chen2024method} include a discussion of universality with respect to the distribution of $\bX$ and extend their framework to doubly-robust functionals, thereby adopting a broader scope.

Table \ref{tab:method_comparison} summarizes the difference between our study and related methods. Accordingly, our method differs from existing approaches by achieving two key objectives: ensuring applicability even when $\kappa$ is greater than 1 and enabling the estimation of hyperparameters beyond the SS needed for inference, e.g., variance of Gaussian noise.
These are justified both theoretically and numerically.
The intuitive reason our method can accommodate the $\kappa\geq 1$ setting, albeit at the cost of imposing an asymmetric assumption, is that it relies solely on the moment information of $Y$. As a result, since our method does not use the information of $\bX${\cc; in particular,} the inverse of the sample covariance matrix $(n^{-1}\sum_{i=1}^n\bX_i\bX_i^\top)^{-1}$ does not appear. This is a key advantage, as it ensures robustness against the instability of coefficient estimators and the inverse matrices as $\kappa$ increases.
On the other hand, a drawback of using only the moment information of $Y$ is that SS cannot be identified when the link function is symmetric.
}
\begin{table}[t]
    \centering
        \caption{Comparison of methods based on $\kappa$ values.
        \label{tab:method_comparison}}
    \begin{tabular}{l|cccc}
    \hline 
        {\textbf{Method}}  & {Model} & { $\kappa < 1$} & { $ \kappa \geq 1$} & {Condition on the link $g(\cdot)$} \\
        \hline\hline 
        ProbeFrontier \cite{sur2019likelihood} & Logistic & $\checkmark$ &   & sigmoid $1/(1+\exp(-x))$\\ 
        SLOE \cite{yadlowsky2021sloe} & GLM & $\checkmark$ &   & \\ 
        MoMI \cite{chen2024method} & GLM & $\checkmark$ & {\ccc Limited}  & \\ 
        \textbf{Ours} & GLM & $\checkmark$ & $\checkmark$  & restricted to asymmetric\\
        \hline
    \end{tabular}
\end{table}

Without SE systems, \cite{bellec2025observable} introduced an inference framework that remains valid even when the link function is unknown under the identification condition $\mathrm{Var}(\bX^\top\bbeta)=1$. 
{\rc
However, when the link function is specified (i.e., the GLM case), this approach is not directly applicable since we cannot constrain SS to be 1.}
As for classification problems, \cite{mai2019high,mai2019large,seddik2021unexpected} established feasible methods for inference though they require prior knowledge of the covariance matrix of the covariates.

Various theoretical tools have been proposed for analyzing statistical models in the regime \eqref{eq:high_dim}, including (i) AMP algorithms \citep{donoho2009message,bolthausen2014iterative,bayati2011dynamics}, (ii) convex Gaussian minimax theorem \citep{thrampoulidis2015lasso,thrampoulidis2018precise}, (iii) leave-one-out techniques \citep{el2013robust}, (iv) second-order Poincar\'e inequalities \citep{chatterjee2009fluctuations}, and (v) second-order Stein's formulae \citep{bellec2021second,bellec2022derivatives}. In addition to (i), \cite{rangan2011generalized} proposed 
generalized AMP applicable to GLMs. We use this to characterize the asymptotic behavior of an estimator of GLMs.

In a high-dimensional regime different from \eqref{eq:high_dim},
statistical inferences for GLMs have also been studied \citep{guo2016tests,jankova2016confidence,salehi2019impact,cai2021statistical,li2023estimation,sawaya2025provable}, typically under sparsity conditions.
{\rc
We note that the regime \eqref{eq:high_dim} admits the number of nonzero coefficients $s$ to be proportional to $p$, while classical high-dimensional statistics typically focus on sparse models with $s\log p=o(n)$.
}

\subsection{State Evolution (SE)-based Confidence Interval: Logistic Regression Case}\label{sec:surcandes}

We review statistical inference frameworks that are valid in the high-dimensional regime \eqref{eq:high_dim}, especially those aimed at confidence intervals.
This framework for for a logistic regression is developed by \cite{sur2019modern,zhao2022asymptotic}. 
Specifically, the studies showed that a properly corrected maximum likelihood estimator (MLE) is asymptotically normally distributed {\cc under the Gaussian design $\bX_i\iid\mathcal{N}(\bzero,\bSigma(n))$ with $\bSigma(n)$ a covariance matrix}.
{\cc Here, the vector $\bbeta=\bbeta(n)$ in the model is assumed to satisfy 
\[\lim_{n\to\infty}\mathrm{Var}(\bX_i(n)^\top\bbeta(n))=\lim_{n\to\infty}\bbeta(n)^\top\bSigma(n)\bbeta(n)=\gamma^2<\infty,\]
where we call $\gamma^2$ the \textit{signal strength} (SS).}
Let $\hat{\bbeta}^{\mathrm{MLE}} = (\hat{\beta}_1,\dots,\hat{\beta}_p)^\top$ be the MLE.
Then, as $n,p(n)\to\infty$ with $p(n)/n\to\kappa\in(0,1)$, we have the following for any $j=1,\dots,p$ satisfying $\sqrt{n}\tau_j\beta_j=O(1)$ if the MLE exists:
\begin{align}
    \frac{\sqrt{n}(\hat{\beta}_j^\mathrm{MLE}-\mu\beta_j)}{\sigma { /\tau_j}}\overset{d}{\to}\mathcal{N}(0,1), \label{def:asym_norm}
\end{align}
where $\tau_j^2:=\mathrm{Var}(X_{ij}\mid\bX_{i\setminus j} )$ is the conditional variance of the $j$-th element of $\bX_i$ given $\bX_{i\setminus j} := (X_{i,1},\dots,X_{i,j-1},X_{i,j+1},\dots,X_{i,p})$.

In this logistic regression case, the inverse link and the integrated inverse link functions are defined as $g(t) = 1/(1 + \exp(-t))$ and $G(t):= \log(1+\exp(t))$, respectively. Here, $(\mu, \sigma, \eta)$ are the \textit{state evolution} (SE) parameters 
that satisfy the following SE system: %{\cc in terms of logistic sigmoid function $g(t)=1/(1+\exp(-t))$}, 
\begin{align}
\label{eq:SE-system-SC}
    \begin{cases}
        \kappa^2 \sigma^2 &= \mathbb{E}_{(Q_1,Q_2)} \sbr{2 g(Z) \rbr{\eta g(D)}^2},\\
        0 &= \mathbb{E}_{(Q_1,Q_2)} \sbr{g(Z) Z (\eta g(D))},\\
        1-\kappa &= \mathbb{E}_{(Q_1,Q_2)} \sbr{ 2 g(Z)( 1 + \eta g'(D))^{-1} },
    \end{cases}
\end{align}
with $D=\mathrm{prox}_{\eta G}\rbr{\mu Z +\sqrt{\kappa}\sigma Q_2}$, $Z=\gamma Q_1$, and $(Q_1,Q_2)\sim \mathcal{N}_2(\bzero,\bI_2)$. 
{
In this system, the covariance structure of $\bX_i$ appears only through $\gamma^2$.}
Remarkably, SE-based inference reduces the inference problem on high-dimensional coefficients to an SE system with a small number of scalar SE parameters. The SE system is derived as a fixed-point equation of an iterative algorithm for an optimization problem on estimators, 
which is based on the approximate message passing, the details of which are provided in \cite{feng2022unifying}.

On the {\cc basis} of the asymptotic normality of the MLE in \eqref{def:asym_norm}, we may construct a confidence interval for each $\beta_j$ for $j=1,\ldots,p$, provided that the SE parameters are available.
Specifically, given $(\mu,\sigma)$, the MLE-centric confidence interval with confidence level $1 - \alpha \in (0,1)$ is formed as
\begin{align}
    \mathrm{CI}_{(1-\alpha)} := \left[ \frac{\hat{\beta}^{\mathrm{MLE}}_j}{\mu} - z_{(1-\alpha/2)} \frac{\sigma { /\tau_j}}{\sqrt{n}\mu},  \frac{\hat{\beta}^{\mathrm{MLE}}_j}{\mu} + z_{(1-\alpha/2)} \frac{\sigma { /\tau_j}}{\sqrt{n}\mu}\right], \label{def:CI}
\end{align}
where $z_{\alpha/2}$ is the $(\alpha/2)$-quantile of a standard Gaussian distribution. 
Then we obtain $\mathrm{Pr}(\beta_j \in \mathrm{CI}_{(1-\alpha)}) \to 1 - \alpha$ as $n,p(n) \to \infty$ in the sense of \eqref{eq:high_dim}.
Importantly, we need to know the SS parameter $\gamma^2$ in advance, because the SE parameters $(\mu, \sigma)$ depend on $\gamma^2$.
{\cc
\begin{remark}
    The results presented in this section are intended solely to clarify the motivation. We emphasize that some of our subsequent results do not cover logistic regression.
\end{remark}}

\subsection{Notation}
Define $\R_+=(0,\infty)$. 
For a vector $\bb = (b_1,\dots,b_p)\in \R^p$, $b_j$ denotes {\cc the} $j$-th element of $\bb$ for $j=1,\dots,p$. 
For function $f: \R \to \R$, $f'(\cdot)$ denotes the derivative of $f(\cdot)$.
We say $f(\cdot)$ is $C$-\textit{Lipschitz} with $C > 0$ if $f(\cdot)$ is Lipschitz continuous and its Lipschitz constant is $C$.
For any convex function $f:\R\to\R$ and constant $\eta>0$, we define the proximal operator $\prox_{\eta f}:\R\to\R$ as
$\prox_{\eta f}(x)=\argmin_{z\in\R}\cbr{\eta f(z)+(x-z)^2/2}$. 
{\cc Let $\Phi:\R\to\R$ be the cumulative distribution function of the standard Gaussian distribution.}

\section{Preliminary}

\subsection{Generalized Linear Model}

Suppose that we observe {\cc $n$ i.i.d.}\  pairs $\{(\bX_i,Y_i)\}_{i=1}^n$ of a feature vector $\bX_i\equiv\bX_i(n)\in\R^{p}$ and a target variable $Y_i \equiv Y_i(n)\in \mathcal{Y}$ that follow the GLM \eqref{eq:model}, where $\mathcal{Y}$ is a response space, such as $\R,\R_+,\{0,1\},\{0,1,2,\ldots\}$, and so on.  
Hereafter, we drop the dependence on $n$ whenever it is clear from context.
We assume that the feature vector is generated independently from $\bX_i \sim \mathcal{N}_p(\bzero,\bSigma)$ using the covariance matrix $\bSigma\equiv\bSigma(n) \in \R^{p\times p}$.
The GLM can represent several models by specifying the inverse link function $g(\cdot)$ and distribution of $Y_i$ for a given $\bX_i$: for example, $Y_i\mid \bX_i\sim\mathrm{Ber}(g(\bX_i^\top\bbeta))$ with $g(t) = 1/(1 + \exp(-t))$ for the logistic regression model, and $Y_i\mid \bX_i\sim\mathrm{Pois}(g(\bX_i^\top\bbeta))$ with $g(t) = \exp(t)$ for the Poisson regression model.

If a random variable $Y_i\mid \bX_i$ has a density function (continuous case) or a mass function (discrete case) $f(\cdot \mid \bX_i)$, then the maximum likelihood estimator (MLE) of $\bbeta$ is defined as $\hat{\bbeta}^{\mathrm{MLE}}(n)=\argmax_{b\in\R^p}\sum_{i=1}^n \log f(Y_i\mid \bX_i)$.
In a low-dimensional setting where $n\to\infty$ with fixed $p < \infty$, the MLE $\hat{\bbeta}^{\mathrm{MLE}}$ asymptotically obeys the normal distribution as $\sqrt{n}(\hat{\bbeta}^{\mathrm{MLE}}-\bbeta)\overset{d}{\to}\mathcal{N}_p(\bzero,\mathcal{I}_{\bbeta}^{-1})$, where $\mathcal{I}_{\bbeta}$ is the Fisher information matrix at the true coefficient $\bbeta$. Once the matrix is estimated consistently, we can consider the inference of $\bbeta$ conventionally. 
By contrast, this classical approach is no longer valid in a proportionally high-dimensional setting \eqref{eq:high_dim}, as presented in \cite{sur2019modern}.

\subsection{Our Goal}

Our goal is to construct a valid confidence interval for {\cc each coordinate of} the unknown coefficient vector $\bbeta$ of the GLM \eqref{eq:model} that is valid in the high-dimensional regime \eqref{eq:high_dim}. 
Specifically, for $\alpha \in (0,1)$ we construct a random set $\mathrm{CI}_{(1-\alpha)} \subset \R$ such that we obtain $\mathrm{Pr}(\beta_j \in \mathrm{CI}_{(1-\alpha)}) \to 1 - \alpha$ {\cc in} the limit \eqref{eq:high_dim}.

\section{Method}

We {\cc construct} a confidence interval for GLMs in the high-dimensional regime {\cc in} three steps, as overviewed in Section \ref{sec:surcandes}:
\begin{enumerate}
    \item We construct a surrogate estimator for the coefficient {\cc vector} $\bbeta$ of GLM, and derive its SE systems.
    \item We develop moment-based estimators for parameters in the derived SE system. This step constitutes the core of our novelty. 
    \item We construct a confidence interval based on the surrogate estimator and the moment-based estimator.
\end{enumerate}
Each step is detailed in the following subsections.

\subsection{Surrogate Estimator for Coefficient Vector \texorpdfstring{$\bbeta$}{beta}}
\label{sec:surrogate}

\subsubsection{Definition}

Given an inverse link function $g(\cdot)$, we develop an estimator for the true coefficient vector $\bbeta$ by using a \textit{surrogate loss} function  \citep{auer1995exponentially,agarwal2014least}, which is characterized by $\ell: \mathbb{R}^p \times \mathbb{R}^p \times \mathcal{Y}\to\R$ with
\begin{align}
\ell(\bb;\bx,y):=G\rbr{\bx^\top \bb}-y \bx^\top \bb, \label{def:surrogate_loss}
\end{align}
where $G(\cdot)$ is a function satisfying $G'=g$.
Using the loss, we define the surrogate estimator 
\begin{align}
\label{eq:estimator}
    \hat{\bbeta}(n):=\argmin_{\bb\in\R^p}\sum_{i=1}^n\ell(\bb;\bX_i,Y_i).
\end{align}

The use of a surrogate estimator is justified in two ways.
First, if $g(\cdot)$ monotonically increases, the surrogate loss is convex in $\bb$.
This property provides computational and statistical advantages for a broader class of GLMs and plays a crucial role in the derivation of SE systems.
Second, the minimizer of surrogate risk coincides with the true coefficient $\bbeta$, which is expressed as follows:
\begin{proposition}[Lemma 1 in \cite{agarwal2014least}]
\label{prop:agarwal}
Consider an $\mathbb{R}^p \times \mathcal{Y}$-valued random element $(\bX,Y)$ that follows GLM \eqref{eq:model}.
If $g(\cdot)$ is increasing, then the coefficient vector $\bbeta$ in \eqref{eq:model}  satisfies $\bbeta=\argmin_{\bb\in\R^p}\E\sbr{\ell(\bb;\bX,Y)\mid\bX}.$
\end{proposition}
This claims that the surrogate loss estimator is a reasonable extension of MLE for GLMs. 
Modeling the distribution of $Y_i \mid \bX_i$ to follow an exponential dispersion family \citep{jorgensen1987exponential} and choosing $g(\cdot)$ as the inverse of the canonical link function make the surrogate estimator \eqref{eq:estimator} equivalent to the MLE, covering logistic and Poisson regression. 

{\ccc
We also consider the following ridge-penalized surrogate estimator, which is applicable to the setting $n\le p$:
\begin{align}
\label{eq:estimator-ridge}
\hat{\bbeta}^\lambda(n)\in\argmin_{\bb\in\R^p}\cbr{\sum_{i=1}^n\ell(\bb;\bX_i,Y_i)+{\ccc\frac{n\lambda}{2}\|\bb\|^2 }},
\end{align}
for any $\lambda\ge0$.
If $\lambda=0$, it reduces to \eqref{eq:estimator}.
}

\subsubsection{SE System for Surrogate Estimator of GLMs}

We derive an SE system associated with the {\ccc general} surrogate estimator \eqref{eq:estimator-ridge}.
As this estimator is given by the minimization problem of the convex cost function, we can derive the SE system, as described in Section \ref{sec:surcandes}.
{\cc
Assume that $g(\cdot)$ is differentiable almost everywhere.}
Recall that $G(\cdot)$ is a function satisfying $G'=g$, and $\gamma^2=\lim_{n\to\infty}\mathrm{Var}(\bX_i(n)^\top\bbeta(n))$ is the signal strength (SS) parameter.
{\ccc Denote $\mathsf{Sym}_+(p)$ be the class of $p\times p$ symmetric positive definite matrices.}
\begin{proposition}
\label{prop:SE-GLM}
{\ccc 
Suppose that $\|\hat\bbeta^\lambda\|=O_p(1)$, and the empirical distributions $p^{-1}\sum_{j=1}^p\delta_{\sqrt{n}\beta_j}$ and $n^{-1}\sum_{i=1}^n\delta_{\varepsilon_i}$ converge in the $2$-Wasserstein sense to the distributions of $\Bar{\beta}$ and $\Bar{\varepsilon}$ with finite second order moment, respectively.
Assume one of the following two regimes:
\begin{enumerate}
    \item[(i)] $\lambda>0, \quad \bSigma=\bI_p, \quad \text{and} \quad \kappa\in(0,\infty),$
    \item[(ii)] $\lambda=0, \quad \text{arbitrary}\ \   \bSigma\in\mathsf{Sym}_+(p), \quad \text{and} \quad \kappa\in(0,1).$ 
\end{enumerate}
Then, the SE system of $\hat\bbeta^\lambda$ is given by,}
    \begin{align}
\label{eq:SE-system}
\begin{cases}
\kappa^2\sigma^2&=\eta^2\E_{(Q_1,Q_2,\Bar{Y})}\sbr{\rbr{\Bar{Y}-g\rbr{L}}^2}, \\
\gamma^2\lambda\mu&=\E_{(Q_1,Q_2,\Bar{Y})}\sbr{Z\rbr{\Bar{Y}-g\rbr{L}}},\\
    1-\kappa+\lambda\eta&=\E_{(Q_1,Q_2,\Bar{Y})}\sbr{({\displaystyle 1+\eta g'\rbr{L}})^{-1}},
\end{cases}
\end{align}
where $L=\prox_{\eta G}(\mu Z+\sqrt{\kappa}\sigma Q_2+ \eta \Bar{Y})$, $Z=\gamma Q_1$, $(Q_1,Q_2)\sim \mathcal{N}_2(\bzero,\bI_2)$ and $\Bar{Y}=h(Z,\bar{\varepsilon})$.
\end{proposition}
{\ccc
We define the solution to \eqref{eq:SE-system} as $(\mu_\lambda,\sigma_\lambda,\eta_\lambda)$, and we write $(\mu,\sigma,\eta)=(\mu_0,\sigma_0,\eta_0)$.}

In view of \eqref{eq:SE-system}, the SE parameters $(\mu_{\ccc\lambda},\sigma_{\lambda}^2,\eta_{\ccc\lambda})$ depend on three components: (i) $\kappa = \lim_{n\to\infty}p(n)/n$, (ii) the SS parameter  $\gamma^2=\lim_{n\to\infty}\mathrm{Var}(\bX_i(n)^\top\bbeta(n))$, and (iii) the distribution of $\bar{Y}$ determined by the model.
Some numerical plots of the equations \eqref{eq:SE-system} will be presented in Section \ref{sec:numerical_se}.
{\cc Note that, by taking $g(t)={1}/({1+\exp(-t)})$ and $\bar{Y}=\mathbf{1}\{g(Z)\ge \varepsilon\}$ with $\varepsilon\sim\mathrm{Unif}[0,1]$, one can verify that \eqref{eq:SE-system-SC} is a special case of \eqref{eq:SE-system}.
This can be seen in Section H.3.1 of \citep{sur2019modern}.}

{\cc Constructing data-driven approximations to} the parameters of this equation is a non-trivial task.
Although we can set $\kappa=p(n)/n$ and regard it as known, the estimation of $\gamma^2$ and other parameters (if present) of the distribution of $\bar{Y}$ for each GLM is non-trivial. 
In the next subsection, we propose a method for estimating them.

\begin{remark}\label{remark:example}
    If the distribution of ${Y}\mid\bX$ has a single parameter (e.g., Bernoulli, Poisson, and exponential distribution), then 
    {\bc
    one can specify $\varepsilon\sim\mathrm{Unif}[0,1]$ and}
    the model assumption $\E[Y \mid \bX]=g(\bX^\top\bbeta)$ implies that $\gamma^2$ fully characterizes the SE parameters $(\mu_{\ccc\lambda},\sigma_{\lambda}^2,\eta_{\ccc\lambda})$. Thus, in this case, the SS parameter $\gamma^2$ is the only unknown parameter to be estimated. Meanwhile, estimating the distribution of ${Y}\mid\bX$ is difficult if it has multiple parameters. Even in such cases, we propose a method for their estimation, focusing on the case of ${Y}\mid\bX$ following a normal distribution; see Section \ref{sec:moment_regression}.
\end{remark}

\subsection{Moment-Based Estimation for SE Parameter} \label{sec:hyper_parameter_estimation}

We propose a novel estimator for the SS parameter $\gamma^2$ and other parameters of the distribution of $\bar{Y}$ in \eqref{eq:SE-system} to access the SE parameters $(\mu_{\ccc\lambda},\sigma_{\lambda}^2)$ for inference.
Once we estimate $\gamma^2$ and the parameters of $\Bar{Y}$, we can solve system \eqref{eq:SE-system} by substituting the estimators to obtain the SE parameters.

\subsubsection{Additive Form of GLM}

We convert the GLM \eqref{eq:model} into the \textit{additive model} with observations $(Y_i,\bX_i)$:
\begin{align}
    Y_i=g(\bX_i^\top\bbeta)+e_i\quad{\rm with}\quad e_i = Y_i - \mathbb{E}[Y_i \mid \bX_i],\quad i=1,\dots,n,\label{eq:glm_model2}
\end{align}
where $e_i$ is a noise variable that satisfies $\E[e_i\mid \bX_i]=0$.
Note that the parameters of the distribution of $e_i\mid\bX_i$ may depend on $\bX_i,~i=1,\dots,n$. 
While the original GLM \eqref{eq:model} highlights the conditional mean, the additive model \eqref{eq:glm_model2} provides a more analyzable form of the distribution by introducing the variable $e_i$. 
This structure is essential for our following estimation method.

\subsubsection{Moment-Based Estimation when $Y\mid\mathbf{X}$ is Gaussian} \label{sec:moment_regression}

We first apply our moment-based method to the Gaussian output:   ${Y}\mid \bX \sim \mathcal{N}(g(\bX^\top \bbeta),\sigma_{{e}}^2)$ with an unknown variance parameter $\sigma_{{e}}^2>0$. 
This is a typical case of GLMs including a nonlinear regression setup, but it has additional unknown parameter $\sigma_{{e}}^2$ to be estimated, as well as the SS parameter $\gamma$.
We first consider this rather complicated case where we estimate both $\gamma$ and $\sigma_{{e}}^2$ as it serves as an appropriate introduction to our approach.

At the beginning, we characterize the parameters $\gamma^2$ and $\sigma_{{e}}^2$ by the moments of $Y$.
In the limit, the additive model \eqref{eq:glm_model2} has the form
\begin{align}
\label{eq:limit-additive}
    \bar{Y} = g
    (\gamma Q) + \bar{e} \quad \text{with}\quad Q \sim \mathcal{N}(0, 1) \quad \text{and}\quad \bar{e} \sim \mathcal{N}(0, \sigma^2_e),
\end{align}
where $Q$ and $\bar{e}$ are independent.
Then, we obtain the condition with second and fourth moments as
\begin{align} \label{eq:regression_moment}
    \Psi(\gamma,\sigma_e):=
    \begin{bmatrix}
    \E[\bar{Y}^2] - \E[g(\gamma Q)^2] - \sigma_e^2 \\ 
    \E[\bar{Y}^4] - \E[g(\gamma Q)^4] - 6 \E[g(\gamma Q)^2] \sigma_e^2 - 3\sigma_e^4
    \end{bmatrix}
    =\bzero\in\R^2.
\end{align}
This simultaneous equation characterizes the parameter $(\gamma, \sigma_e^2)$ as a solution.
{\cc
The sufficient conditions for the uniqueness of the solution to this system are discussed in Appendix \ref{sec:identification}. 
}
We consider estimation of $(\gamma, \sigma_e^2)$ using an empirical analog of the equation \eqref{eq:regression_moment}.
Since $g(\cdot)$ is known, we can simulate $h_2 (\varsigma) := \E[g({\varsigma Q})^2]$ and $h_4 (\varsigma) := \E[g({\varsigma Q})^4]$ by generating ${\varsigma Q}\sim \mathcal{N}(0,\varsigma^2)$ for each $\varsigma>0$. 
Using the observations $\{Y_i\}_{i=1}^n$, we then define the estimators as the solution to
\begin{align} \label{eq:empirical_equation}
&(\hat{\gamma}, \hat{\sigma}_e) := \left\{ (\varsigma, \varsigma_e) \in \mathbb{R}_+^2 : \Psi_n(\varsigma,\varsigma_e)=\bzero
\right\},\\
\mbox{where~}&\Psi_n(\varsigma,\varsigma_e)=
\begin{bmatrix}
        n^{-1} \sum_{i=1}^n Y_i^2 -  h_2(\varsigma) - \varsigma_e^2 \\
        n^{-1} \sum_{i=1}^n Y_i^4 - h_4(\varsigma) - 6 h_2(\varsigma) \varsigma_e^2 - 3 \varsigma_e^4
\end{bmatrix}
\in\R^2.
\end{align}
The equation in \eqref{eq:empirical_equation} can be viewed as an empirical analogue of \eqref{eq:regression_moment}, provided that the convergence $n^{-1} \sum_{i=1}^n Y_i^a \overset{\rm a.s.}{\longrightarrow} \E[\bar{Y}^a]$ for $a \in \{2,4\}$ is true under some assumptions. 
The solution is obtained by root-finding algorithms, such as the Gauss-Newton method.

An advantage of the moment-based estimator is the stability for any $\kappa\in(0,\infty)$.
This is because the moment-based estimator is independent of both the estimator $\hat\bbeta$ and an inverse of high-dimensional matrix {\cc$(\bm{X}^\top \bm{X})^{-1}$ with $\bm{X} = (\bX_1,\dots,\bX_n)^\top$, which becomes unstable} as $\kappa$ increases {\cc and is ill-defined when $\kappa > 1$}.
Additionally, since this method relies only on observations through $Y$, it is computationally efficient.
{\cc
\begin{remark}
    Above, we construct the SS estimator using second- and fourth-order moments so as to ensure identifiability of the SS when $g$ is an odd function.
    In contrast, if $g$ is not an odd function, we can simplify the estimation equation by only using the first and second moments of $\bar{Y}$ as
\begin{align}
    \begin{bmatrix}
    \E[\bar Y]-\E[g(\gamma Q)]\\
    \E[\bar{Y}^2] - \E[g(\gamma Q)^2] - \sigma_e^2
    \end{bmatrix}
    =\bzero\in\R^2.
\end{align}
\end{remark}
}

\subsubsection{Moment-Based Estimation for {\cc Single-Parameter} Cases} \label{sec:moment_other}

We next study moment-based estimation in other GLMs with a single parameter.
The problem becomes simpler since we only need to estimate $\gamma$, unlike the Gaussian case. 
Examples {\cc and the meaning of {single-parameter}} are provided in Remark \ref{remark:example}.

In this case, we only consider the first moment using the limit form \eqref{eq:limit-additive}:
\begin{align}
    \tilde\Psi(\gamma):=\E[\Bar{Y}]-\E[g(\gamma Q)]=0, \label{eq:limit-additive2}
\end{align}
where we have used the fact that $\E[e_1]=\E[\E[e_1\mid\bX_1]]=0$ by the assumption. The estimator of $\gamma$ is obtained as a solution to the {\cc following} linear equation with $h_1(\varsigma) := \E[g(\varsigma Q)]$:
\begin{align}
\label{eq:gamma-hat}
    \textstyle \hat{\gamma}:=\cbr{\varsigma \in \R_+ :\tilde\Psi_n(\varsigma)=0},
    \mbox{~~where~}\tilde\Psi_n(\varsigma):=n^{-1}\sum_{i=1}^nY_i-h_1(\varsigma).
\end{align}
Similarly to Section \ref{sec:moment_regression}, this estimator is stable for any $\kappa\in(0,\infty)$ and computationally efficient {\cc as it does not rely on quantities (e.g., $\hat\bbeta$ and $(\bX^\top\bX)^{-1})$ that suffer from instability in high dimensions}.{\cc
\begin{example}
    In a Poisson regression model, the estimator of $\gamma$ is defined as the solution to $n^{-1}\sum_{i=1}^n Y_i = \mathbb{E}\left[\exp({Z_{\varsigma})}\right]$.
Since it follows that
$\mathbb{E}\left[\exp({Z_{\varsigma})}\right] = \exp({\varsigma^2/2})$ by Gaussianity, the estimator admits the closed-form expression
\[
\hat{\gamma}^2 
= 2 \log\left( \frac{1}{n}\sum_{i=1}^n Y_i \right).
\]
\end{example}
}
If necessary, we can construct an estimator considering higher-order moments as in Section \ref{sec:moment_regression}.
{\cc We provide such an example where $Y\mid X$ follows gamma distributions.
\begin{example}
    We consider another example of GLMs with the following distribution
\[
Y \mid \bm{X} \sim \mathrm{Gamma}(k,\theta),
\]
where $k>0$ and $\theta>0$ denote the shape and scale parameters, respectively. In this case, the conditional mean and variance satisfy
\[
\mathbb{E}[Y \mid \bm{X}] = g(\bm{X}^\top \bm{\beta}) = k\theta,
\qquad
\mathrm{Var}(Y \mid \bm{X}) = k\theta^{2}.
\]
Since we can write the moments of $Y$ as $\E[ Y]=\E[\mathbb{E}[Y \mid \bm{X}]]=\E[g(\bm{X}^\top \bm{\beta})]$ and $\mathrm{Var}( Y)=\E[\mathrm{Var}(Y \mid \bm{X})]+\mathrm{Var}(\E[Y \mid \bm{X}]) =k^{-1}\E[g(\bm{X}^\top \bm{\beta})^2]+\E[g(\bm{X}^\top \bm{\beta})^2]-\E[g(\bm{X}^\top \bm{\beta})]^2$, 
% Thus, we have
% \begin{align}
%     \E[ Y]&=\E[\mathbb{E}[Y \mid \bm{X}]]=\E[g(\bm{X}^\top \bm{\beta})],\\
%     \mathrm{Var}( Y)&=\E[\mathrm{Var}(Y \mid \bm{X})]+\mathrm{Var}(\E[Y \mid \bm{X}]) \\
%     &=k^{-1}\E[g(\bm{X}^\top \bm{\beta})^2]+\E[g(\bm{X}^\top \bm{\beta})^2]-\E[g(\bm{X}^\top \bm{\beta})]^2.
% \end{align}
we obtain the following moment equations regarding $(\gamma,k)$
\begin{align*}
    \begin{cases}
    \ \E[\bar Y] = \E[g(\gamma Q)],\\
    \ \E[\bar{Y}^2] = \left(1+k^{-1}\right)\E[g(\gamma Q)^2],
    \end{cases}
\end{align*}
which provide a valid estimator as a solution to their sample analogue.
\end{example}
}

\subsection{Confidence Interval} \label{sec:theory_CI}

We construct a confidence interval for $\bbeta$ by using the estimators in Sections \ref{sec:hyper_parameter_estimation} and $\hat{\bbeta}^{\ccc\lambda}$ in \eqref{eq:estimator}.
Let $(\hat{\mu}_{\ccc\lambda},\hat{\sigma}_{\lambda}^2,\hat\eta_{\ccc\lambda})$ be the solutions to the SE system \eqref{eq:SE-system}, with $\gamma^2$ replaced by $\hat{\gamma}^2$. 
We also introduce an estimator $\hat{\tau}_j$ of the conditional variance $\tau_j^2$, the details of which are provided in Section \ref{sec:tauj}.
Then, for any $j\in\{1,\dots,p\}$, our confidence interval with a confidence level $1-\alpha \in (0,1)$ is defined as:
\begin{align}
\label{eq:correct-CI}
    \mathcal{CI}_{{1-\alpha},j}^\lambda := \sbr{\frac{\hat{\beta}_j^\lambda}{\hat{\mu}_\lambda}-z_{(1-\alpha/2)}\frac{\hat{\sigma}_\lambda}{\sqrt{n}\hat{\mu}_\lambda\hat{\tau}_j}, \frac{\hat{\beta}_j^\lambda}{\hat{\mu}_\lambda}+z_{(1-\alpha/2)}\frac{\hat{\sigma}_\lambda}{\sqrt{n}\hat{\mu}_\lambda\hat{\tau}_j}}.
\end{align}
%{\mc [Need some explanation on the covariance estimator.]}

\section{Theory}
\label{sec:theory}
We establish the asymptotic properties of the estimators derived above. Specifically, we prove the asymptotic normality of the surrogate estimator and the consistency of the moment-based estimator, which justifies the validity of the confidence intervals.
%The validations of the surrogate/moment-based estimators are necessary to prove the validity of the confidence interval.

\subsection{Asymptotic Normality of Surrogate Estimator}

We derive the theoretical results of the surrogate estimator $\hat{\bbeta}^{\ccc\lambda}$. 
First, we give the following assumptions.
\begin{assumption}
We consider the following conditions: \label{asmp:basic} 
\setlength{\leftmargini}{27pt} 
    \begin{itemize}
      \setlength{\parskip}{0cm}
  \setlength{\itemsep}{0cm}
    \item[(A1)] $\bX \sim\mathcal{N}_p(\bzero,\bSigma)$, and $\mathrm{Var}(\bX^\top\bbeta)$ converges to $\gamma^2\in(0,\infty)$ as $n\to\infty$.
    {\cc
    Additionally, the empirical measure of the coordinates of $\bm{\beta}$, $\frac{1}{p}\sum_{j=1}^p \delta_{\sqrt{n}\beta_j},$
converges in the $2$-Wasserstein distance to a limiting distribution with finite second moment.}
    \item[(A2)] {\cc The} inverse link function $g:\R\to\R$ is monotonically increasing and $L$-smooth (i.e., the derivative of $g(\cdot)$ is $L$-Lipschitz continuous). 
    \item[(A3)] We can write $Y=h(\bX^\top\bbeta,\varepsilon)$ with a non-random function $h: \R^2 \to \R$ and an $\R$-valued random variable $\varepsilon$ independent of $\bX$, such that $\varepsilon$ has a finite second moment and $h(\cdot,\cdot)$ is Lipschitz continuous with respect to the first argument.
    \item[(A4)] There exists a positive solution to the SE system \eqref{eq:SE-system}.
\end{itemize}
\end{assumption}

{\cc (A1) is standard in developing the asymptotic theory based on approximate message passing (see, e.g., \citep{bayati2011dynamics,sur2019modern,feng2022unifying}).}
(A3) is a technical requirement for the Lipschitzness of SE parameters with respect to $\gamma^2$. Examples of $h(\cdot,\cdot)$ are provided in Section \ref{sec:h-example}. 
{\cc
As noted there, the function $h$ in the logistic regression model is not continuous, and therefore an additional approximation argument is required for a formal justification (though it is out of the scope of this paper).
(A4) is standard and often imposed either implicitly or explicitly in the related literature (e.g., \citep{loureiro2021alearning,li2026spectrum}).
That said, there are several papers that do study existence/uniqueness of SE solutions in specific settings.
For instance, \citep{bellec2023existence} establish uniqueness of the SE solution for certain $M$-estimators under a linear model.
However, to the best of our knowledge, analogous general results for the SE system arising from \emph{generalized} linear models are not currently available.}
In the sequel, we focus only on the situation in which the estimator \eqref{eq:estimator} exists asymptotically; that is, {\cc$\|\hat\bbeta\|=O_p(1)$}.

We obtain the asymptotic normality of the adjusted test statistics for each coordinate of $\hat\bbeta$ with the oracle SE parameters $\mu$ and $\sigma^2$.
\begin{proposition} \label{prop:asym_norm_betahat}
    Suppose that we know a solution to the system of nonlinear equations \eqref{eq:SE-system}. 
    {\ccc
    Assume one of the following two regimes:
\begin{enumerate}
    \item[(i)] $\lambda>0, \quad \bSigma=\bI_p, \quad \text{and} \quad \kappa\in(0,\infty),$
    \item[(ii)] $\lambda=0, \quad \text{arbitrary}\ \   \bSigma\in\mathsf{Sym}_+(p), \quad \text{and} \quad \kappa\in(0,1).$ 
\end{enumerate}}
    Under (A1), (A2), and (A4) in Assumption \ref{asmp:basic}, and $\|\hat\bbeta^{\ccc\lambda}\|=O_p(1)$,  we have the following for $j = 1,\dots,p(n)$ satisfying $\sqrt{n}\tau_j|\beta_j|=O(1)$:
    \begin{align}
        \frac{\sqrt{n}(\hat{\beta}_j^{\ccc\lambda}-\mu_{\ccc\lambda}\beta_j)}{\sigma_{\ccc\lambda}/\tau_j}\dconv\mathcal{N}(0,1).
    \end{align}
{\cc Furthermore, if {\ccc $\max_{j\in[p]}\sqrt{n}\tau_j|\beta_j|=O(1)$} holds, we have
\begin{align}
    \max_{j\in[p]} \sup_{t\in\R} \abs{\pr\rbr{\frac{\sqrt{n}(\hat{\beta}_j^{\ccc\lambda}-\mu_{\ccc\lambda}\beta_j)}{\sigma_{\ccc\lambda}/\tau_j}<t}-\Phi(t)}\to0.
\end{align}
}
\end{proposition}
This is the first result to demonstrate the marginal asymptotic normality of estimators for {\cc \textit{general}} GLM settings {\cc beyond logistic regression MLE} in a high-dimensional setting.
This differs from the classical result $\sqrt{n}(\hat\bbeta^{\rm MLE}-\bbeta)\dconv\mathcal{N}(\bzero,\mathcal{I}_{\bbeta}^{-1})$ in a setup with $p < \infty$. 
We regard the marginal convergence as {\cc the consequence of a straightforward application} of \cite{zhao2022asymptotic}{\cc's proof technique}. 
Although the values $\mu$ and $\sigma$ are unknown in practice, a feasible construction of the {\cc statistic has been discussed in Section~\ref{sec:hyper_parameter_estimation} and its} asymptotic normality is provided in Proposition \ref{prop:t-stat}.

{\ccc The reason we allow only isotropic designs in the penalized case is that the
equivariance of the estimator under orthogonal transformations of the design
matrix, which is a key ingredient of our proof, is lost under a general
covariance structure.}

\subsection{Consistency of Moment-Based Estimator}

We show the consistency of the estimators for the SS parameters in the cases of Sections \ref{sec:moment_regression} and \ref{sec:moment_other}.
We need the following assumption.
\begin{assumption} \label{asmp:g_estimation}
    We consider the following conditions:
\setlength{\leftmargini}{27pt} 
    \begin{itemize}
      \setlength{\parskip}{0cm}
  \setlength{\itemsep}{0cm}
        \item [(A5)] For any $\epsilon>0$, we have either: 
        \begin{align}
            \inf_{(\varsigma,\varsigma_e):|\varsigma-\gamma|+|\varsigma_e-\sigma_e|>\epsilon}\|\Psi(\varsigma,\varsigma_e)\|>0,
        \end{align}
        {\cc if $Y\mid\bX$ follows the Gaussian distribution with mean $g(\bX^\top\bbeta)$ as discussed in Section \ref{sec:moment_regression}}, or 
        \begin{align}
            \inf_{\varsigma:|\varsigma-\gamma|>\epsilon}|\tilde\Psi(\varsigma)|>0,
        \end{align}
        {\cc if $Y\mid\bX$ follows the distribution having single parameter $g(\bX^\top\bbeta)$ as discussed in Section \ref{sec:moment_other}}. 
        \item [(A6)] There exist constants $C,c,\epsilon>0$ such that {\cc$|g(z)|\le\exp(C|z|^{2-\epsilon}+c)$} for any $z \in \R$. 
    \end{itemize}
\end{assumption}
Condition (A5) identifies the SS parameters $(\gamma,\sigma_e^2)$ by the moment equations. This condition is not satisfied if the moments of $Y$ do not have sufficient information, such as the case of logistic regression. 
{\rc
However, a later study \cite{chen2024method} proposes a method that utilizes information from $\bX$ to handle logistic regression using our moment-based technique especially when $n>p$. }
We also discuss a sufficient condition for (A5) in Section \ref{sec:identification}.
Condition (A6) is required for the existence of $\E[g(\gamma Q)^q]$ for $q \geq 1$ used in \eqref{eq:limit-additive} and \eqref{eq:limit-additive2}.
{\cc We note that Assumption~(A6) is implied by Assumption~(A2), which imposes $L$-smoothness on $g$. However, we state it separately here in order to make explicit the minimal conditions required for the consistency of the signal strength estimator.}

We then obtain the following results for the consistency of the estimators for the SS parameters:

\begin{theorem}
\label{thm:consistent_moment}
    Under Assumption \ref{asmp:basic} (A1) and Assumption \ref{asmp:g_estimation}, $\hat{\gamma}^2$ defined in \eqref{eq:gamma-hat} satisfies $\hat\gamma^2\overset{\mathrm{a.s.}}{\longrightarrow}\gamma^2$
    and $(\hat\gamma^2,\hat\sigma_e^2)$ defined in \eqref{eq:empirical_equation} satisfies $(\hat\gamma^2,\hat\sigma_e^2)\overset{\mathrm{a.s.}}{\longrightarrow}(\gamma^2,\sigma_e^2)$, as {\cc $n\to\infty$}.
{\cc
    In addition, if $\E[Y_i^2]<\infty$ and $h_{\ccc1}'(\gamma)\neq0$,
    we have, for $\hat{\gamma}$ in \eqref{eq:gamma-hat}, as $n\to\infty$,
    \begin{align}
    \label{eq:gamma-asyN}
        \sqrt{n}(\hat\gamma-\gamma)=\frac{1}{h_1'(\gamma)}\cdot\frac{1}{\sqrt{n}}\sum_{i=1}^n(Y_i-h_1(\gamma))+o_p(1).
    \end{align}
    Also, if $\E[Y_i^8]<\infty$ holds and
    \begin{align}
        A:=	\begin{pmatrix}
   -h_2'(\gamma) & -2\sigma_e \\
   -(h_4'(\gamma)+6h_2'(\gamma)\sigma_e^2) & -12\sigma_e(h_2(\gamma)+\sigma_e^2)
\end{pmatrix}
    \end{align}
    is invertible,
    we have, for $(\hat\gamma,\hat\sigma_e)$ in \eqref{eq:empirical_equation}, as $n\to\infty$,
    \begin{align}
    \label{eq:sigma-asyN}
        \sqrt{n}
        \begin{pmatrix}
            \hat{\gamma}-\gamma\\\hat\sigma_e-\sigma_e
        \end{pmatrix}
        =-A^{-1}\cdot\frac{1}{\sqrt{n}}\sum_{i=1}^n
        \begin{pmatrix}
            Y_i^2-\E[\bar Y^2]\\
            Y_i^4-\E[\bar Y^4]
        \end{pmatrix}
        +o_p(1).
    \end{align}
}
\end{theorem}
{\cc
The result above holds regardless of $p$ or $\kappa$.
The asymptotic linearity of \eqref{eq:gamma-asyN} and \eqref{eq:sigma-asyN} implies the asymptotic normality of the SS estimators and the  $\sqrt{n}$-consistency (i.e., $|\hat{\gamma} - \gamma | = O_p( n^{-1/2})$).
The same rate of convergence is also established in \citep{chen2024method}'s estimator of signal strength.
}

{\cc
\begin{remark}
Assumption~(A5), which is required to ensure identifiability for SS estimation, excludes logistic regression. Consequently, the consistency of the SS estimator (Theorem~\ref{thm:consistent_moment}), the normal approximation results based on the estimated SS (Proposition~\ref{prop:t-stat} below), and the resulting confidence interval (Theorem~\ref{thm:conf-intrvl} below) also exclude logistic regression.

On the other hand, the remaining theoretical results (Propositions~\ref{prop:agarwal}--\ref{prop:asym_norm_betahat}) do not depend on this assumption and therefore continue to hold for general GLMs, including logistic regression.
\end{remark}
}
\subsection{Asymptotic Validity of Proposed Confidence Interval}

At last, we demonstrate the asymptotic validity of our bias correction for statistical inference. 
To observe this, we use the asymptotic normality with the estimated SE parameters {\ccc under the following additional assumption}.

{\ccc
\begin{assumption}
\label{asmp:stable-sol}
Let $\bm\vartheta_0=(\mu_\lambda,\sigma_\lambda^2,\eta_\lambda)$ denote the selected positive solution to the state-evolution system at the
true signal strength $\gamma$. 
There exist a compact interval $\mathcal I_\gamma\subset(0,\infty)$ satisfying $\gamma\in\operatorname{int}(\mathcal I_\gamma)$, 
and a compact set $\Theta\subset\mathbb R\times(0,\infty)^2$
containing $\bm\vartheta_0$ such that the following hold:
\begin{enumerate}
    \item for every $\varsigma\in\mathcal I_\gamma$, the SE
    system obtained by replacing $\gamma$ with $\varsigma$ admits at least
    one solution in $\Theta$;
    \item at $\varsigma=\gamma$, $\bm\vartheta_0$ is the unique solution
    in $\Theta$; and
    \item whenever $\hat\gamma\in\mathcal I_\gamma$, the plug-in
    SE parameters $\hat{\bm\vartheta}=(\hat\mu_\lambda,\hat\sigma_\lambda^2,\hat\eta_\lambda)$
    are selected as a measurable solution in $\Theta$.
\end{enumerate}
\end{assumption}
This is a local stability condition for the solution of the SE system.
}
\begin{proposition}
\label{prop:t-stat}
Suppose that all the conditions in Assumptions \ref{asmp:basic}, \ref{asmp:g_estimation}, {\ccc and \ref{asmp:stable-sol}} hold, {\ccc $\E[h(0,\varepsilon)^2]<\infty$} holds, and $\hat{\tau}_j^2$ is a uniformly consistent estimator of the conditional variance $\tau_j^2$. 
{\ccc
    Assume one of the following two regimes:
\begin{enumerate}
    \item[(i)] $\lambda>0, \quad \bSigma=\bI_p, \quad \text{and} \quad \kappa\in(0,\infty),$
    \item[(ii)] $\lambda=0, \quad \text{arbitrary}\ \   \bSigma\in\mathsf{Sym}_+(p), \quad \text{and} \quad \kappa\in(0,1).$ 
\end{enumerate}}
Then, for every $j = 1,\dots,p(n)$ satisfying $\sqrt{n}\tau_j|\beta_j|=O(1)$, we obtain the following {\cc under $\|\hat\bbeta^{\ccc\lambda}\|=O_p(1)$}:
\begin{align*}
    \frac{\sqrt{n}(\hat{\beta}_j^{\ccc\lambda}-\hat{\mu}_{\ccc\lambda}\beta_j)}{\hat{\sigma}_{\ccc\lambda}/\hat{\tau}_j}\overset{d}{\to}\mathcal{N}(0,1).
\end{align*}
{\cc Furthermore, if {\ccc $\max_{j\in[p]}\sqrt{n}\tau_j|\beta_j|=O(1)$} holds, we have 
\begin{align}
    \max_{j\in[p]} \sup_{t\in\R} \abs{\pr\rbr{\frac{\sqrt{n}(\hat{\beta}_j^{\ccc\lambda}-\hat\mu_{\ccc\lambda}\beta_j)}{\hat\sigma_{\ccc\lambda}/\hat\tau_j}<t}-\Phi(t)}\to0.
\end{align}
}
\end{proposition}

{\rc
To the best of our knowledge, this is the first result establishing the asymptotic normality with the estimated SE parameters in the GLM literature.
{\cc Accordingly, its proof technique is new in the literature, and}
the proof relies on the Lipschitz continuity of SE parameters with respect to SS.
}
Consequently, the proposed confidence interval asymptotically achieves a confidence level $(1-\alpha)$.

\begin{theorem}
\label{thm:conf-intrvl}
    {\cc Consider} the settings of Proposition \ref{prop:t-stat}.  
    Then, for any confidence level $(1-\alpha)\in(0,1)$ and for every $j = 1,\dots,p(n)$ satisfying $\sqrt{n}\tau_j\beta_j=O(1)$, we obtain $\mathrm{Pr}\rbr{\beta_j \in \mathcal{CI}_{{1-\alpha},j}^{\ccc\lambda}}\to1-\alpha$, as $n,p(n)\to\infty$ where $p(n)/n\to\kappa\in(0,1)$.
\end{theorem}

\section{Numerical Experiments}
\label{sec:numerical-exprmnt}
\subsection{Empirical Performance of \texorpdfstring{$\hat{\gamma}$}{gammahat} and \texorpdfstring{$\hat{\sigma}_{e}^2$}{sigma2hat}} \label{sec:exp_empirical}
We fix $n=4,000$ and vary $p=\kappa n$ over $\kappa \in \{0.1,0.2,0.3,0.4,0.5\}$. 
The SS parameter is also fixed at $\gamma^2=1$. We independently generate $n$ realizations of the feature vector $\bX\in\R^p$ from the $p$-variate normal distribution $\mathcal{N}_p(\bzero, \bSigma)$ with $\bSigma_{ij}=0.5^{|i-j|}$. The true regression coefficients $\beta_j$ for $j=1,\ldots,p$ are independently drawn from a normal distribution with the variance determined by $\gamma^2$. 
We calculate $g(\bX^\top\bbeta)$ conditional on $\bX$ and $\bbeta$ and finally draw the response $Y$ from a given distribution.

\begin{figure}[htbp]
    \centering
            \includegraphics[width=0.7\linewidth]{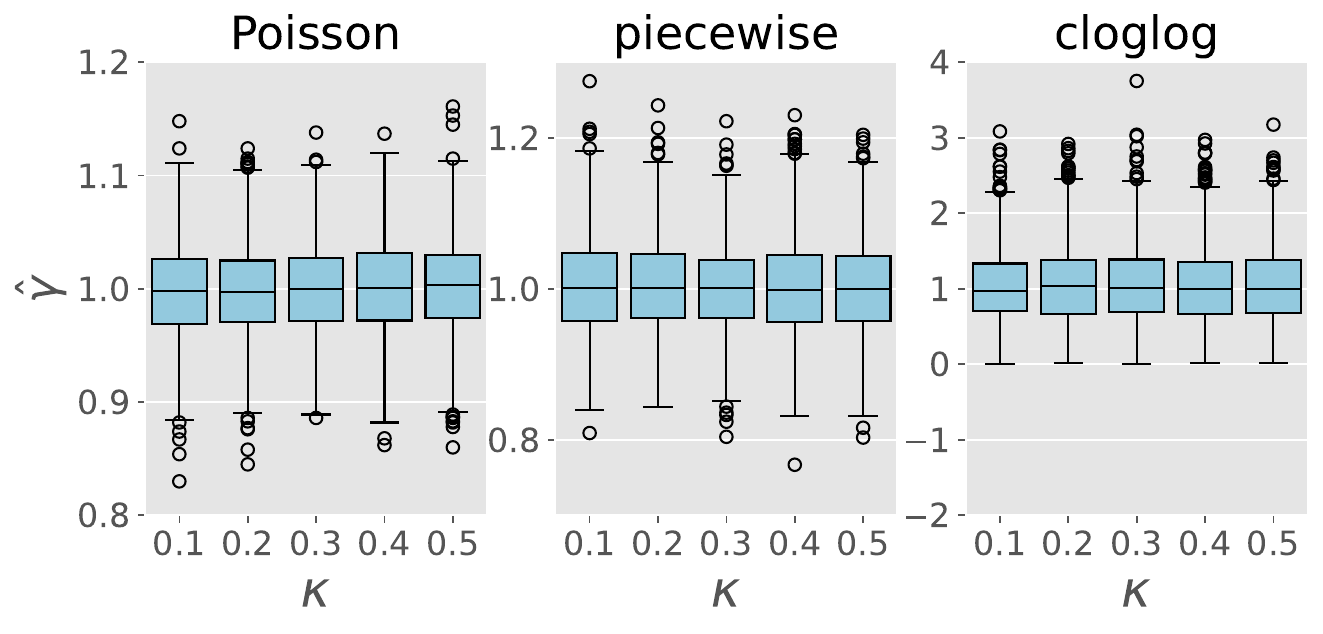}
            \includegraphics[width=0.7\linewidth]{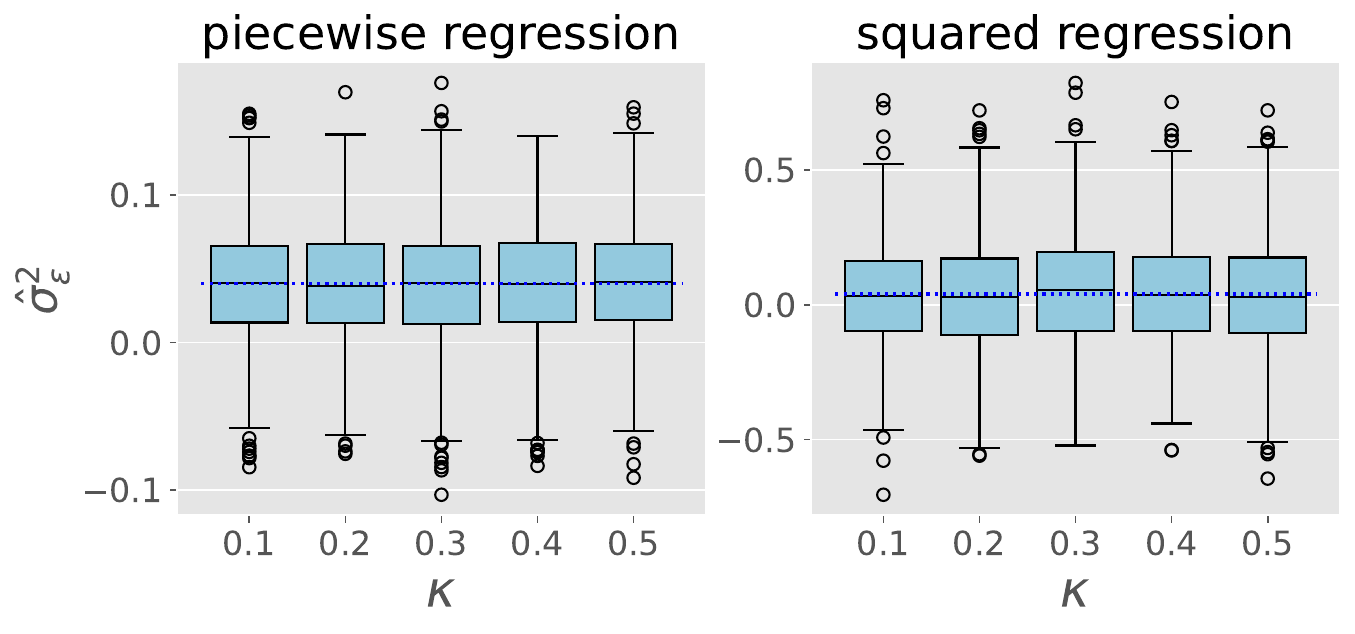}
        \caption{Numerical performance of $\hat{\gamma}$ and $\hat{\sigma}_{e}^2$ over 1,000 simulations. The true values are $\gamma=1$ and $\sigma^2_\varepsilon =  0.04$.}\label{fig:gamma_sig}
\end{figure}

We consider three cases to check the performance of $\hat{\gamma}$ in \eqref{eq:gamma-hat}:
{\ccc
\begin{enumerate}
    \item[(i)] Poisson regression $g(t)=\mathrm{exp}(t)$, $Y\mid \bX\sim\mathrm{Pois}(g(\bX^\top\bbeta))$;
    \item[(ii)] piecewise regression  $g(t)=\min(5t,0.1t)$, $Y\mid \bX\sim \mathcal{N}(g(\bX^\top\bbeta),0.04)$;
    \item[(iii)] complementary log-log (cloglog) regression $g(t)=1-\exp(-\exp(t))$,  $Y\mid\bX\sim\mathrm{Bern}(g(\bX^\top\bbeta))$.
\end{enumerate}
}
We further investigate two situations to evaluate the behavior of $\hat{\sigma}_{e}^2$ in \eqref{eq:empirical_equation}: 
{\ccc
\begin{enumerate}
    \item[(i)] piecewise regression  $g(t)=\min(5t,0.1t)$, $Y\mid \bX\sim \mathcal{N}(g(\bX^\top\bbeta),0.04)$;
    \item[(ii)] squared regression $g(t)=t^2$, $Y\mid \bX\sim \mathcal{N}(g(\bX^\top\bbeta),0.04)$. 
\end{enumerate}
} 
These results are summarized in Figure \ref{fig:gamma_sig}; the upper and lower panels demonstrate stability and consistency of $\hat{\gamma}$ and $\hat{\sigma}_{e}^2$, respectively, regardless of the values of $\kappa$. 

\subsection{Coverage Proportion of Proposed Confidence Interval}

\begin{figure}[htbp]
  \centering
 \includegraphics[width=0.8\textwidth]{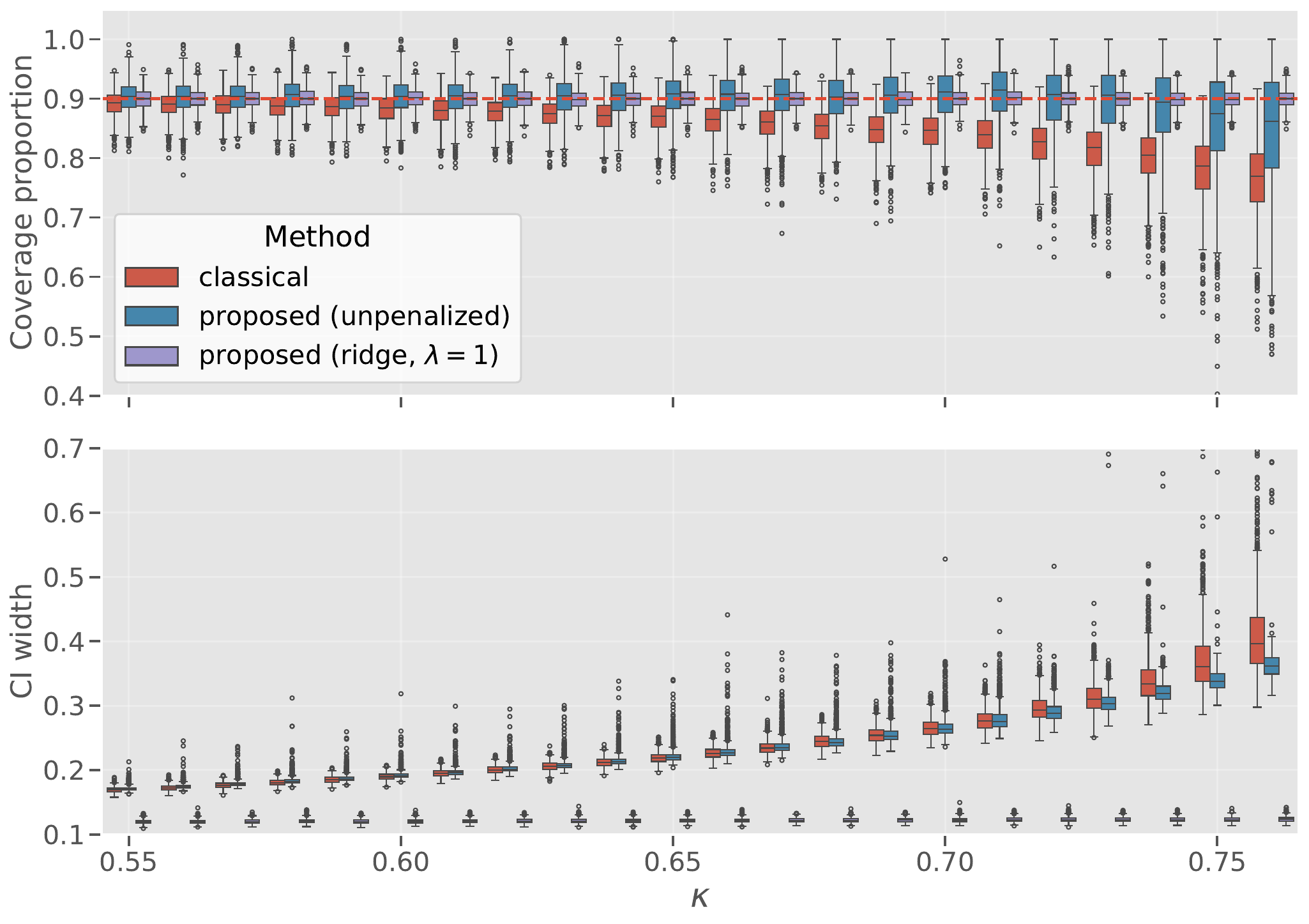}
 \caption{\ccc Empirical coverage probabilities (top) and interval widths (bottom) of nominal 90\% confidence intervals for the regression coefficients, constructed using the classical method (red), the proposed method using the unpenalized estimator (blue), and the proposed method using the penalized estimator (purple), over 1,000 Monte Carlo repetitions with $n=1,000$.}\label{fig:CI}
\end{figure}

We compare the empirical coverage proportion of the proposed confidence interval \eqref{eq:correct-CI} with that of the classical confidence interval with $n=1,000$, {\ccc$\kappa\in\{0.55, 0.56,\ldots,0.76\}$}, and $\gamma=1$.
The proposed confidence interval is calculated by solving a system of nonlinear equations using the estimated SS.
The classical confidence interval for $\beta_j$, $j=1,\ldots,p$, is constructed as
$
    [\hat{\beta}_j \pm z_{(1-\alpha/2)}({\hat{\mathcal{I}}_{jj}^{-1}(\hat{\bbeta})/n})^{1/2}]
$, where $\hat{\mathcal{I}}(\hat{\bbeta})=n^{-1}\sum_{i=1}^ng'(\bX_i^\top\hat{\bbeta})\bX_i\bX_i^\top$ is the empirical Fisher information matrix evaluated at $\hat{\bbeta}$. We analyze the Poisson regression model.
Figure \ref{fig:CI} illustrates that our proposed interval achieves theoretical coverage in all cases; whereas the coverage of the conventional method decreases as $\kappa$ increases.

{\cc
Similarly, Figure~\ref{fig:CI-ridge} presents the results for confidence intervals constructed using the ridge-penalized estimator with $\lambda=1$. In this experiment, we consider an isotropic Gaussian design, fix $n=1,000$, and vary {\ccc$\kappa\in\{0.1, 0.2,\ldots,2.0\}$}, reporting the empirical coverage probabilities and interval widths.
}

{\cc

\begin{figure}[htbp]
  \centering
 \includegraphics[width=0.7\textwidth]{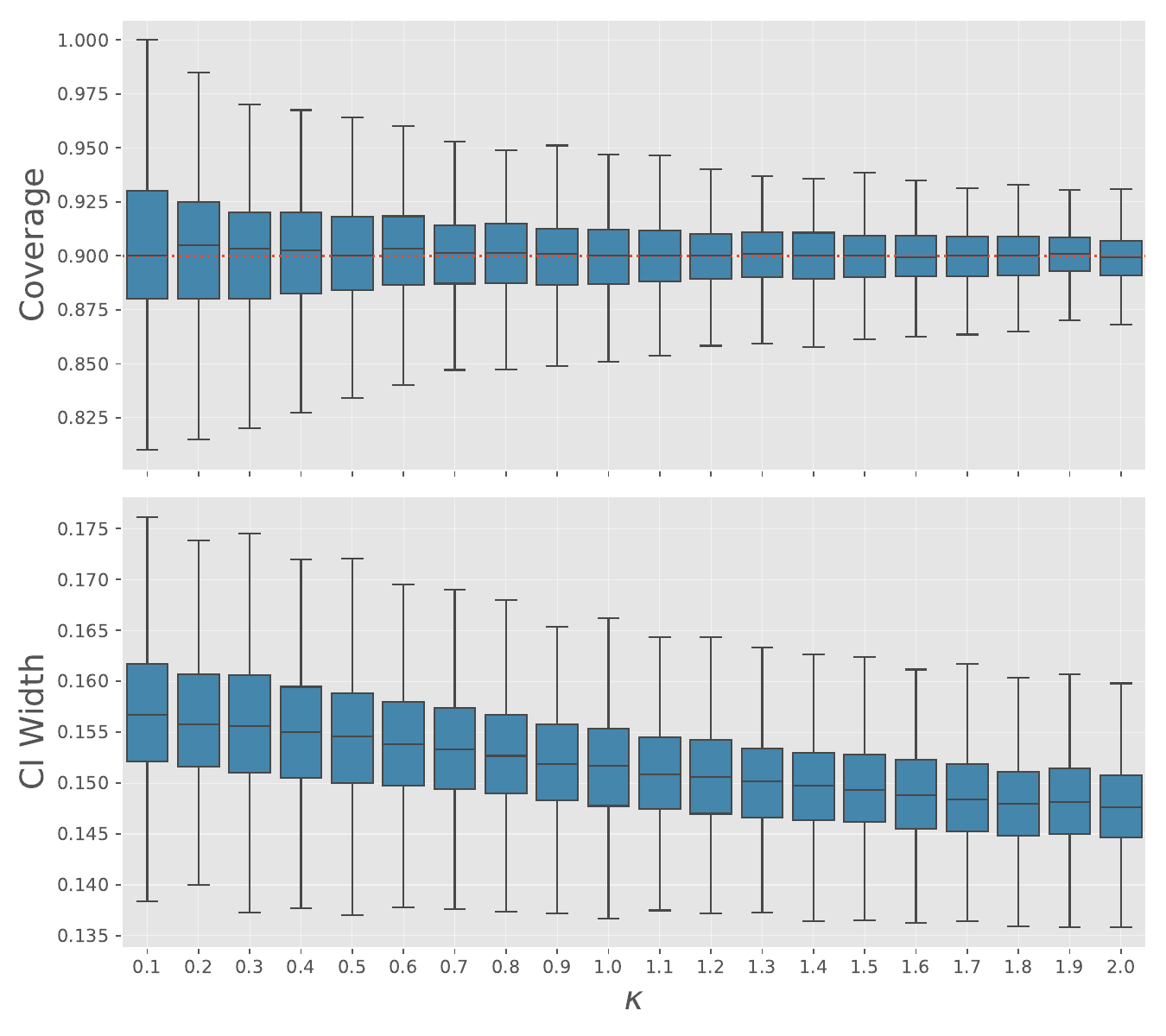}
 \caption{\ccc Empirical coverage probabilities (top) and interval widths (bottom) of nominal 90\% confidence intervals constructed from the ridge-penalized estimator with $\lambda=1$, over 1,000 Monte Carlo repetitions with $n=1,000$.}\label{fig:CI-ridge}
\end{figure}
}

\subsection{Comparison with Other Methods.}

\begin{figure}[htbp]
  \centering
 \includegraphics[width=0.75\hsize]{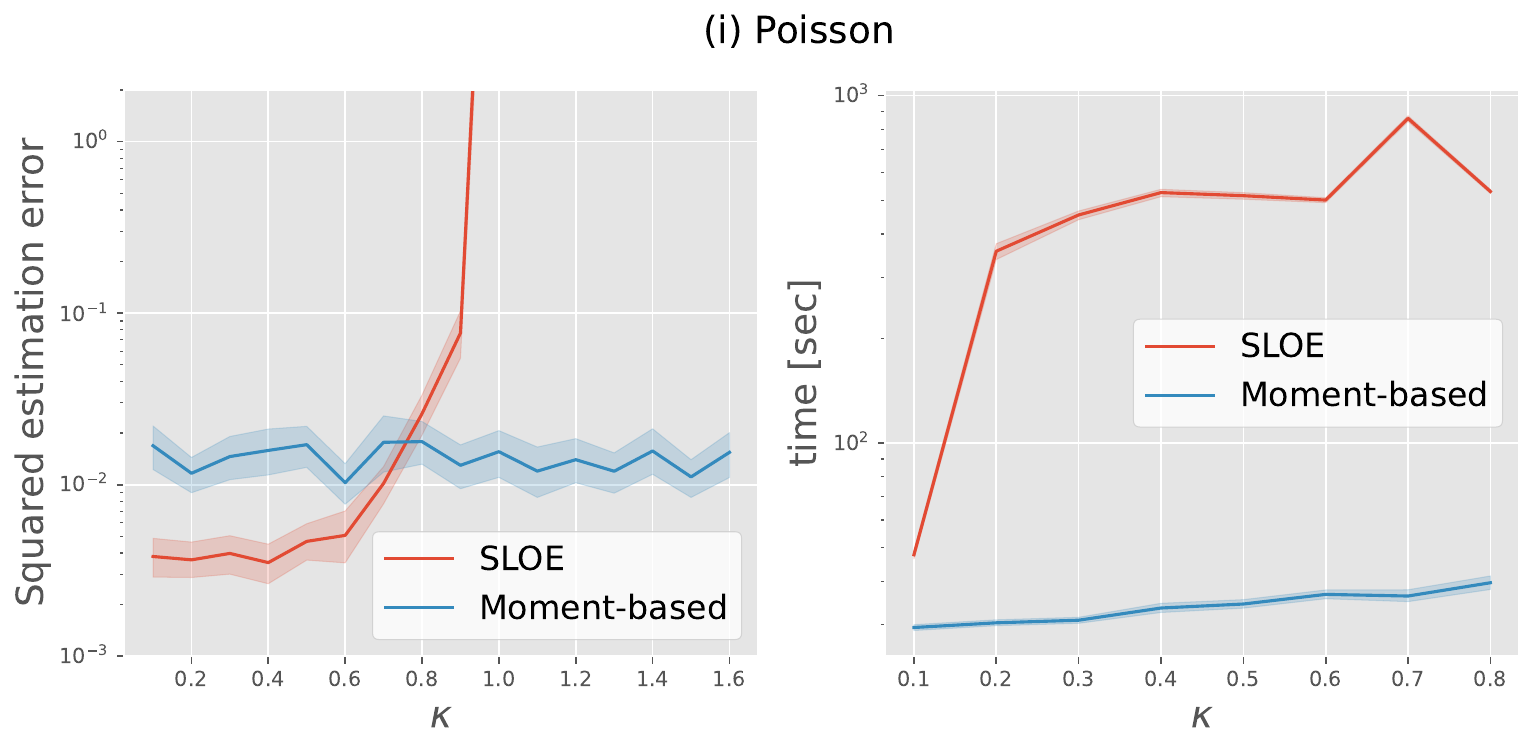} 
% \caption{Scaled squared estimation error of SLOE and our proposed method, in the Poisson regression. The LOO-based method is our extension of the existing method; see Section \ref{sec:SLOE-GLM}.}\label{fig:SS_pois}
%\end{figure}
%\begin{figure}
%  \centering
 \includegraphics[width=0.75\hsize]{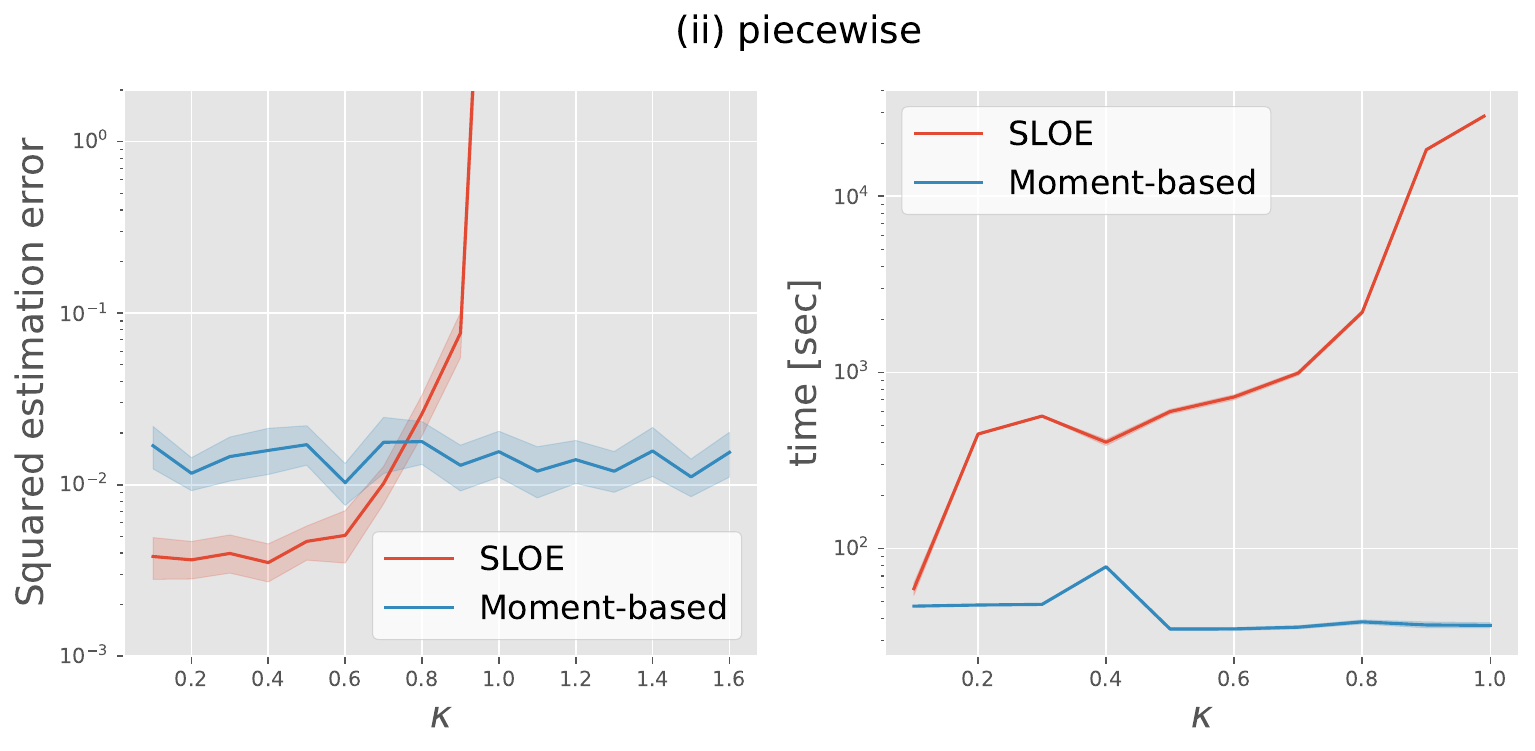}
 \caption{Scaled squared estimation error of SLOE and our proposed method, in the piecewise regression. The LOO-based method is our extension of the existing method; see Section \ref{sec:SLOE-GLM}.}\label{fig:SS_poisson_piece}
\end{figure}

{\rc
We compare the moment-based SS estimator $\hat\gamma^2$ with the previous leave-one-out-based SS estimator, SLOE $\tilde\gamma_\mathrm{SLOE}^2$ \cite{yadlowsky2021sloe}. 
Since the original SLOE is designed for logistic regression, we extend it to GLMs; see Section \ref{sec:SLOE-GLM} for details.
Here, note that the estimand of SLOE is not exactly SS but rather the corrupted signal strength, $\gamma_\mathrm{c}^2=\hat\bbeta^\top\bSigma\hat\bbeta$. 
In the setting with $n=1,000$ and $\gamma=2$, we plot the scaled squared estimation errors $(\hat\gamma/\gamma-1)^2$ and $(\tilde\gamma_\mathrm{SLOE}/\gamma_\mathrm{c}-1)^2$, respectively, over $100$ replications.
We consider the Poisson regression and piecewise regression as in Section \ref{sec:exp_empirical}.

Figure \ref{fig:SS_poisson_piece} indicates that our moment-based estimator performs well even when $\kappa$ exceeds $1$ in both cases while SLOE diverges as $\kappa$ increases and approaches $1$. 
Additionally, regarding the computational load, SLOE requires a longer computation time as $\kappa$ increases; whereas our proposed method remains nearly invariant.
}

\begin{figure}[htbp]
  \centering
 \includegraphics[width=0.75\textwidth]{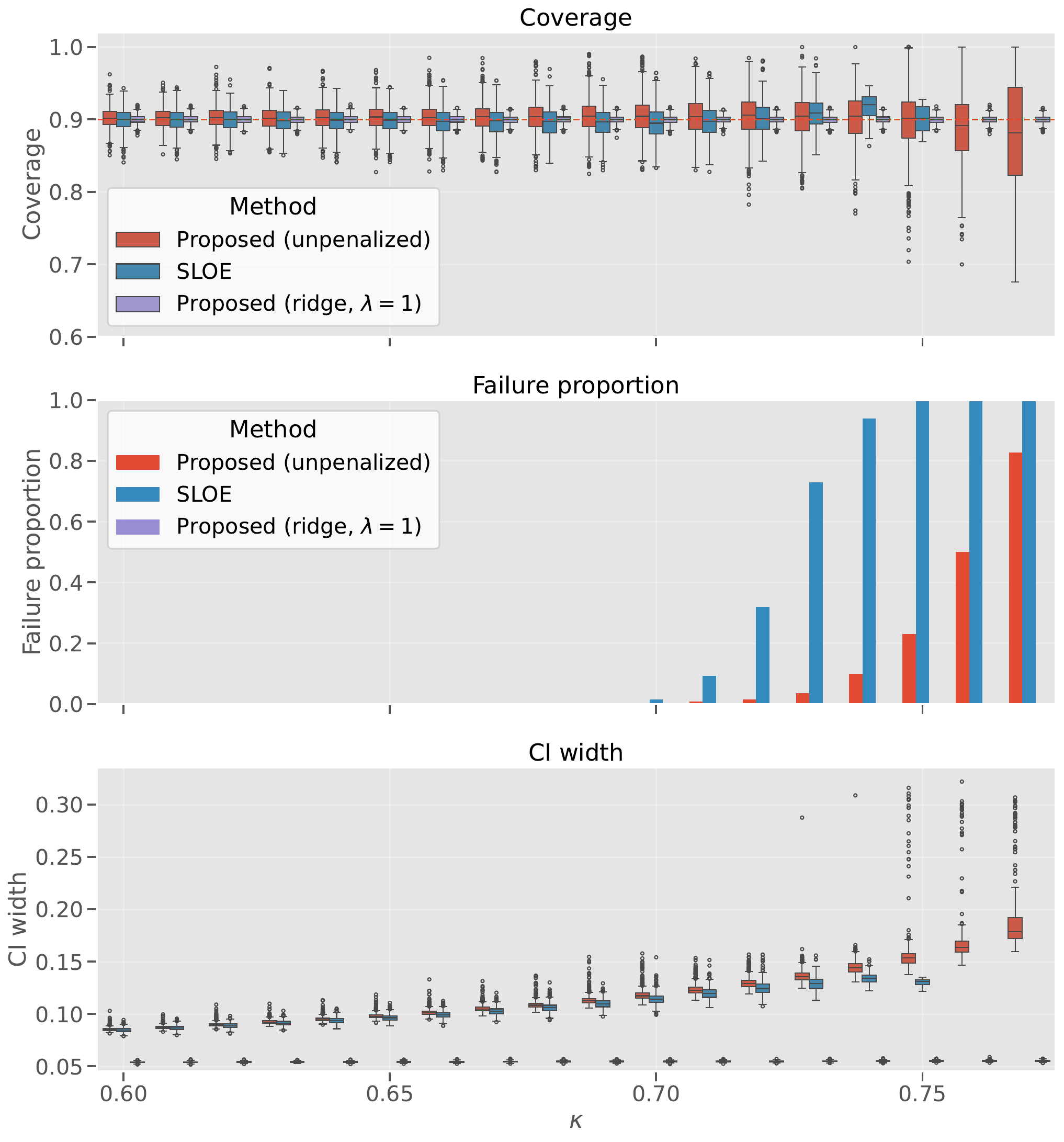}
 \caption{\ccc Comparison of confidence intervals based on the proposed method using the unpenalized estimator (red), SLOE (blue), and the proposed method using the ridge-penalized estimator (purple) for the Poisson regression MLE based on 1,000 Monte Carlo repetitions.
}\label{fig:vsSLOE}
\end{figure}

{\cc
Finally, we report the empirical coverage results of the resulting confidence intervals constructed using the proposed method and SLOE, respectively. 
In Figure~\ref{fig:vsSLOE}, we present the empirical coverage probabilities, average interval widths, and the proportion of failures in solving the fixed-point SE equations. The experiments are conducted for the (penalized) Poisson regression MLE with $n=5,000$, $p=\lfloor \kappa n \rfloor$, and $\gamma=1$, based on 1,000 Monte Carlo repetitions.

Empirically, we observe the following phenomenon. When the estimation error of $\gamma$ is upward biased, the solution of the plug-in fixed-point equation based on the estimated $\gamma$ may diverge, even though the fixed-point equation with the true $\gamma$ remains well-defined. This instability appears more prominently for the SLOE-based method. In Figure~\ref{fig:vsSLOE}, the panel labeled “failure proportion” reports the fraction (out of 1,000 repetitions) for which such divergence occurs.

Overall, as the regime approaches the boundary where $\hat{\bm{\beta}}$ itself diverges, the standard errors associated with the SLOE-based approach tend to inflate substantially. 
In contrast, the proposed method provides uniformly stable estimation of $\gamma$ across different values of $\kappa$, resulting in more {\ccc stable} confidence interval construction. }

{\ccc
Figure~\ref{fig:vsSLOE} also reports the results obtained from the ridge-penalized estimator. 
Ridge regularization substantially reduces the failure proportion over the range of $\kappa$ considered, indicating that regularization effectively mitigates the instability of the unpenalized procedure.
This provides a practical motivation for using the penalized estimator when the unpenalized estimator is close to its instability boundary.
We note, however, that our current theoretical coverage guarantee for the ridge-penalized procedure requires an isotropic design, whereas the unpenalized theory permits a general covariance structure when $\kappa<1$.}

\begin{figure}[htbp]
  \centering
 \includegraphics[width=\textwidth]{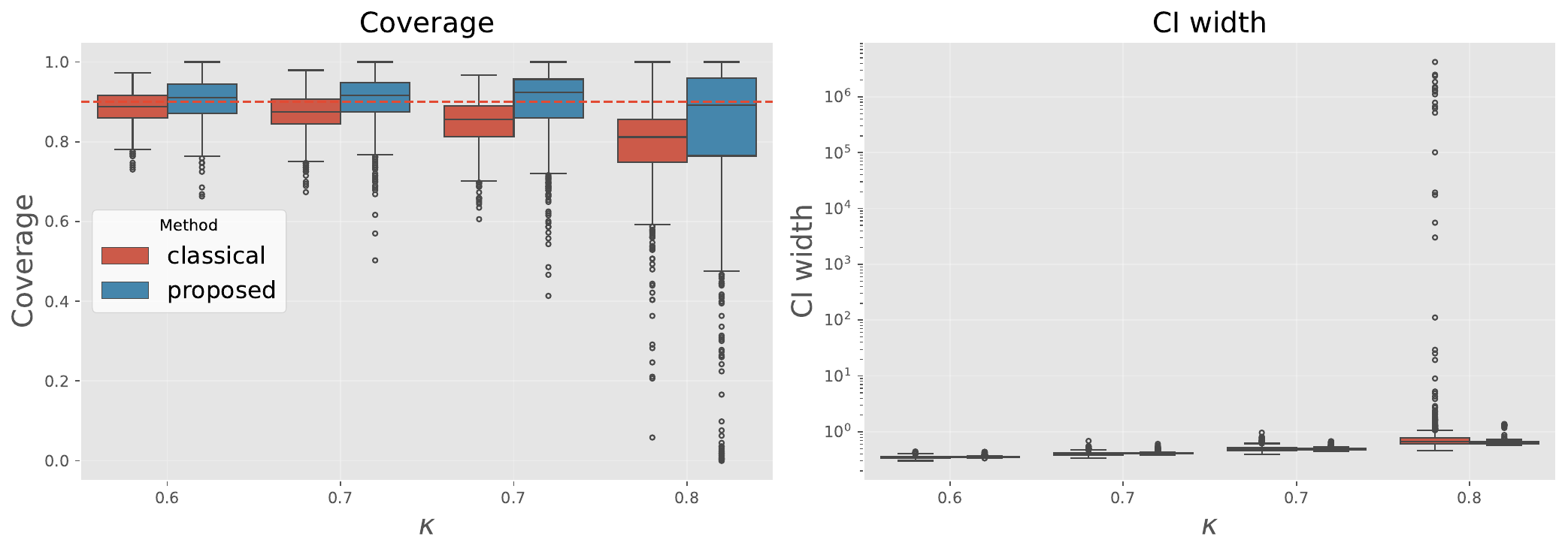}
 \caption{\ccc Empirical coverage probabilities (left) and interval widths (right) of nominal 90\% confidence intervals for the regression coefficients, constructed using the classical method (red) and the proposed method (blue), over 1,000 Monte Carlo repetitions.
}\label{fig:real-CI}
\end{figure}

{\cc
\subsection{Pseudo real-data analysis}

We consider an application based on the Cleveland Clinic Heart Disease dataset \citep{detrano1989international} from the UCI Machine Learning Repository \citep{Dua:2019}. 
The dataset contains $303$ observations with $13$ covariates. The original response variable represents the presence of heart disease, taking integer values from $0$ (no presence) to $4$.

To emulate a proportionally high-dimensional regime, for each fixed $\kappa$ and $i=1,\ldots,303$, we generate an independent random vector from 
$\mathcal{N}\left(\bm{0}, \bm{I}_{\lfloor 303\kappa \rfloor - 13}\right)$
and concatenate it with the original $13$ covariates to form the augmented covariate vector $\bm{X}_i$.
We then standardize $\bX$ so that each column has sample mean zero and unit norm.
For the coefficient vector $\bm{\beta}$, we assign centered normal random variables to the first half of the coordinates including the original 13 covariates, set the second half (corresponding to the artificially added covariates) to zero, and rescale $\bm{\beta}$ so that $\gamma^2 = 1$.

Since the original response variable exhibits excessive zeros, which may violate the assumed conditional structure of $Y \mid \bm{X}$, we instead generate artificial responses according to
$Y \mid \bm{X} \sim \mathrm{Pois}\left(\exp(\bm{X}^\top \bm{\beta})\right),$
so as to focus on universality with respect to the design distribution while avoiding misspecification of $Y\mid \bX$.

Under this setup, for a target level $\alpha$, we construct confidence intervals for each coordinate of $\bm{\beta}$ based on the Poisson regression MLE, using both the classical method and the proposed method. For each value of $\kappa$, we repeat the Monte Carlo simulation (with respect to the artificial generation of $\bX$ and $Y$) 1,000 times. The results are shown in Figure~\ref{fig:real-CI}.
These results demonstrate that the proposed method achieves coverage probabilities close to the target level while maintaining reasonable interval widths.}

\section{Conclusion}

We have developed a statistical inference method for GLMs in high dimensions. 
Our approach extends the SE-based inference method designed for logistic regression model to GLMs.
Specifically, we have proposed a surrogate estimator for GLMs, an associated SE system, and a method to estimate the necessary parameters in the system. 
Our methodology works well in terms of both theory and experimentation. 
One limitation of our method is that it requires Gaussianity of the covariate; however, this limitation {\cc might be relaxed by extending} studies on the universality of methodologies in terms of data, for example, \cite{montanari2022universality,venkataramanan2022estimation,lahiry2024universality,tsuda2026universality}.

{\cc
\section*{Acknowledgments}

We are sincerely grateful to the editor, associate editor, and anonymous reviewers for their careful reading of the manuscript and for their insightful and constructive comments. Their suggestions have significantly improved the clarity, rigor, and overall presentation of the paper. }

\section*{Funding}

The first author was supported by the Japan Society for the Promotion of Science through Research Fellow DC1 (24KJ0841) and the Japan Science and Technology Agency through ACT-X (JPMJAX24CC).
The third author was supported by Japan Society for the Promotion of Science through KAKENHI (24K02904), and the Japan Science and Technology Agency through CREST (JPMJCR21D2),  FOREST (JPMJFR216I), and BOOST (JPMJBY24A9).

\appendix

\section{Supporting Information} 

\subsection{Numericals of SE System} \label{sec:numerical_se}

\begin{figure}[h!]
    \centering
    \begin{adjustbox}{max width=\textwidth}
    \begin{minipage}[b]{0.5\linewidth}
        \centering
        \begin{subfigure}[b]{0.45\linewidth}
            \centering
            \includegraphics[width=\linewidth]{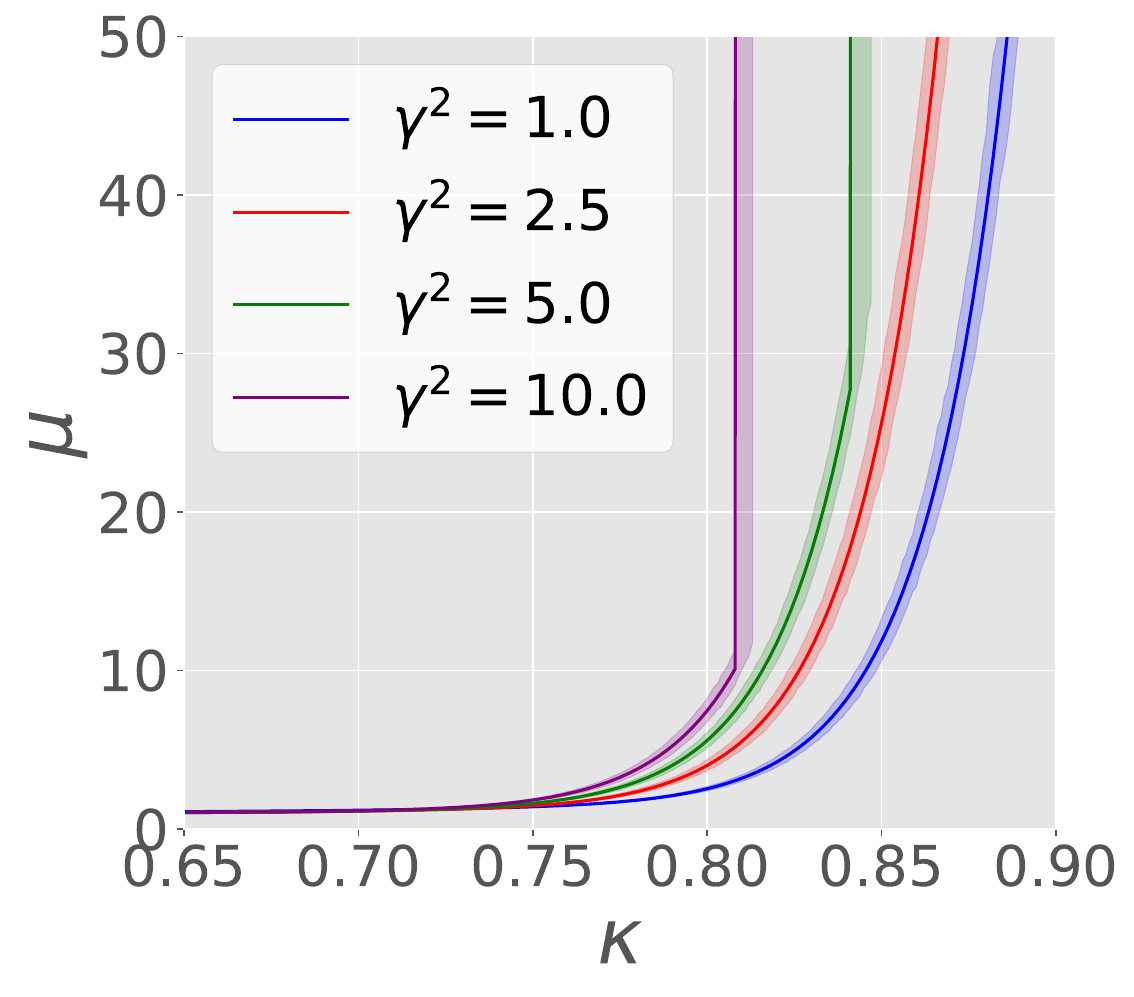}
        \end{subfigure}
        \hfill
        \begin{subfigure}[b]{0.45\linewidth}
            \centering
            \includegraphics[width=\linewidth]{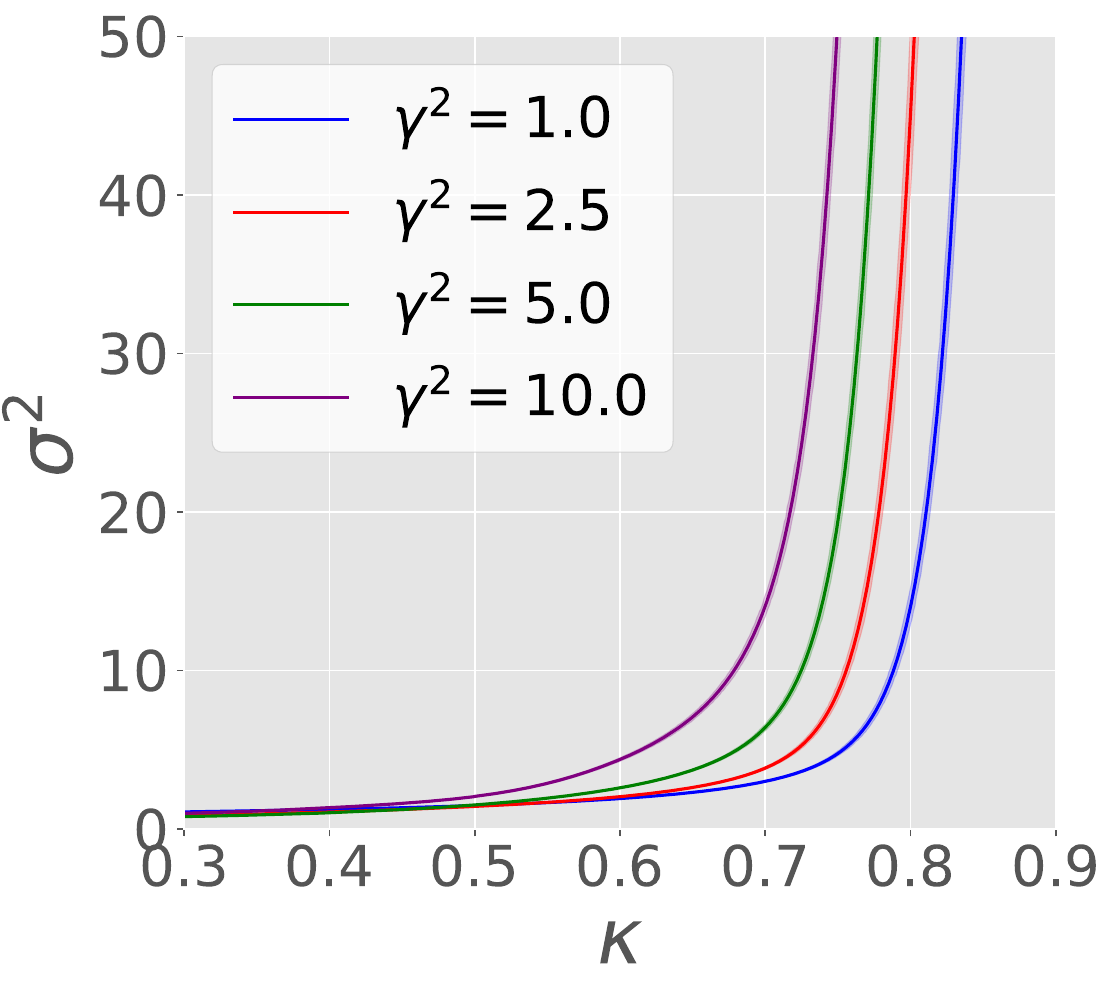}
        \end{subfigure}
        \caption{Computed bias $\mu$ (left) and variance $\sigma^2$ (right) for Poisson regression model $Y\mid\bX\sim\mathrm{Pois}(g(\bX^\top\bbeta))$ with $g(t) = \exp(t)$.  \label{fig:pois_SE}}
    \end{minipage}
    \hfill
    \begin{minipage}[b]{0.5\linewidth}
        \centering
        \begin{subfigure}[b]{0.45\linewidth}
            \centering
            \includegraphics[width=\linewidth]{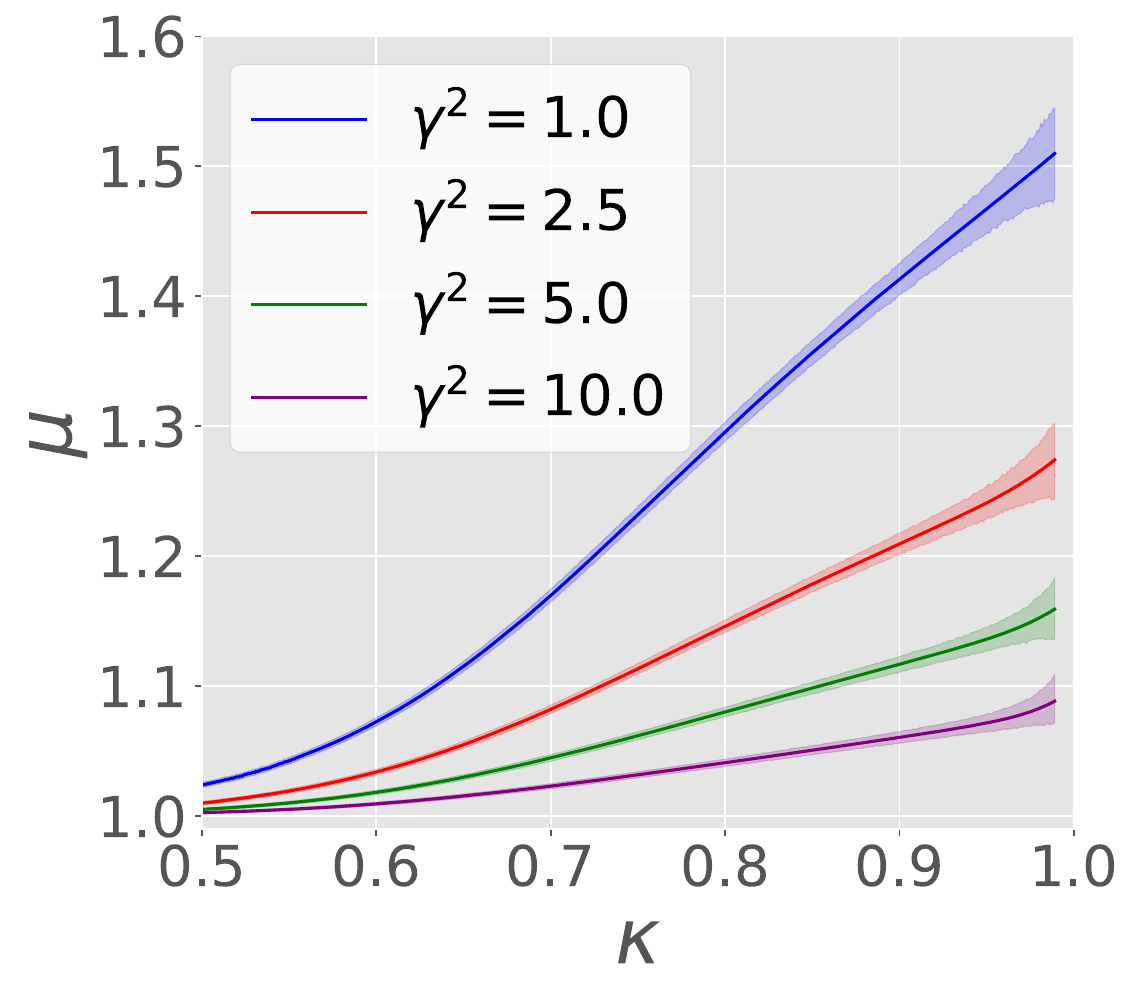}
        \end{subfigure}
        \hfill
        \begin{subfigure}[b]{0.45\linewidth}
            \centering
            \includegraphics[width=\linewidth]{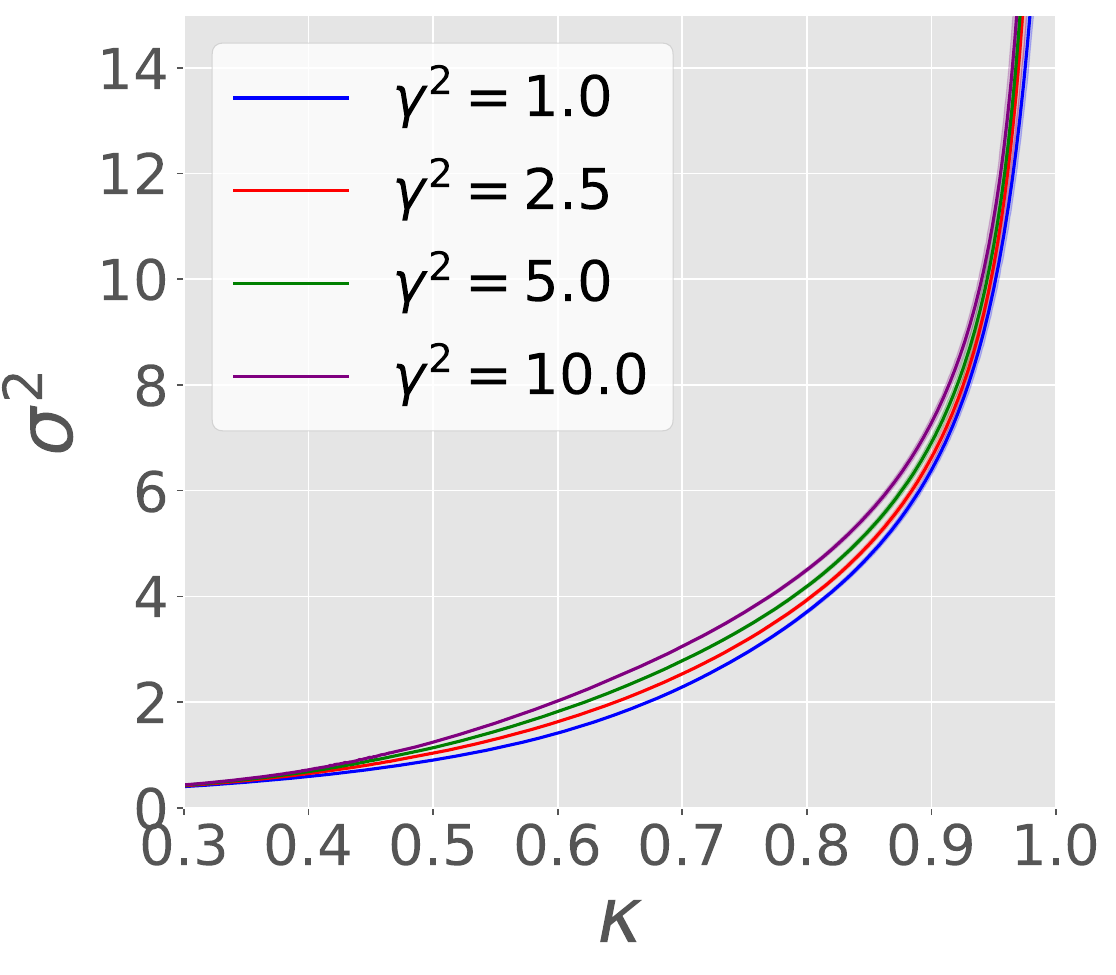}
        \end{subfigure}
        \caption{Computed bias $\mu$ (left) and variance $\sigma^2$ (right) for regression model $Y\mid\bX\sim\mathcal{N}(g(\bX^\top\bbeta),0.2)$ with $g(t)=\min(5t,0.1t)$.\label{fig:piece_SE}} 
    \end{minipage}
    \end{adjustbox}
\end{figure}

We demonstrate the solutions to the system for Poisson regression in Figure \ref{fig:pois_SE} and piecewise regression in Figure \ref{fig:piece_SE}.
These plots are numerical solutions to the SE system over $100$ simulations with their mean and $95\%$ bootstrap confidence intervals.

\subsection{Construction of \texorpdfstring{$\hat{\tau}_j^2$}{hattau2}}
\label{sec:tauj}
To construct the corrected confidence interval \eqref{eq:correct-CI} with correlated features, we have to estimate the conditional variance parameter,  $\tau_j^2=\mathrm{Var}(\bX_{ij}|\bX_{i\setminus j})$, for each $j=1,\ldots,p$. In this article, we follow \cite{zhao2022asymptotic}. They consider the residual sum of squares $\mathrm{RSS}_j$ obtained by regression of $\bX_j=(\bX_{1j},\ldots,\bX_{nj})^\top\in\R^n$ onto a sub-vector of the input $\bX_{\setminus j}=(\bX_1,\ldots,\bX_{j-1},\bX_{j+1},\ldots,\bX_{p})\in\R^{n\times(p-1)}$. 
Owing to the Gaussianity of $\bX_i$, it satisfies
\begin{align*}
    \mathrm{RSS}_j=\bX_j^\top \bP_{\bX_{\setminus j}}^\perp \bX_j\sim\tau_j^2\chi_{n-p+1}^2,
\end{align*}
where $\bP_{\bX_{\setminus j}}^\perp$ is the orthogonal projection matrix onto the orthogonal complement of column space spanned by $\bX_{\setminus j}$. Then, we immediately obtain an unbiased estimator of $\tau_j^2$:
\begin{align*}
    \hat{\tau}_j^2=\frac{1}{n-p+1}\mathrm{RSS}_j.
\end{align*}
{\cc
By Lemma A.3 in \citep{dai2023scale}, this estimator is uniformly consistent over $j \in [p]$ as $n,p \to \infty$ with $p/n \to \kappa \in (0,1)$. Namely,
\begin{align*}
    \max_{j\in[p]}\,({\hat\tau_j}-{\tau_j})\pconv0.
\end{align*}
}

\subsection{Master Theorem of Generalized Approximate Message Passing (GAMP)}
\label{sec:master}

We present the GAMP algorithm by \cite{rangan2011generalized} and its associated theoretical result, which are key tools to inspect the limiting distributional behavior of the estimator {\ccc$\hat\bbeta^\lambda$ for $\lambda\ge0$.}

First, we provide the GAMP algorithm, which generates a sequence of parameters $\tilde{\bbeta}^k$ for an index $k \in \N \cup \{0\}$ based on the minimization problem in \eqref{eq:estimator} and a limit of the sequence corresponds to ${\ccc\sqrt{n}}\hat{\bbeta}^{\ccc\lambda}$.
Let 
$\bX=(\bX_1,\ldots,\bX_n)^\top\in\R^{n\times p}$, $\breve\bX=\bX/\sqrt{n}$, and ${\bY}=(Y_1,\ldots,Y_n)\in\R^n$. Given initial values $\bar{\eta}_0>0,\tilde\bbeta^0$, and $\tilde\bxi^0=\breve\bX\tilde\bbeta^0$, a GAMP recursion takes the form, for each $k\in\N\cup\{0\}$,
{\ccc\begin{align}
    \tilde\bbeta^{k+1}&=\frac{\bar\eta_{k+1}}{\kappa}\breve\bX^\top\cbr{\bY-g\rbr{\prox_{\bar\eta_kG}(\tilde\bxi^k+\bar\eta_k\bY)}}
    +\rbr{{{\ccc1-\kappa^{-1}\lambda\bar\eta_{k}}}}\tilde\bbeta^k,\\
    \tilde\bxi^{k+1}&=\breve\bX\tilde\bbeta^{k+1}-\bar\eta_{k+1}\cbr{\bY-g\rbr{\prox_{\bar\eta_kG}(\tilde\bxi^k+\bar\eta_k\bY)}}.
\end{align}
Here, $\bar\eta_k$ is updated with $\bar\mu_k$ and $\bar\sigma^2_k$ as following: 
\begin{align}
    \bar\eta_{k+1}&=\kappa\bar\eta_k\rbr{\lambda\bar\eta_k+1-\E_{(Q_1,Q_2,U)}\sbr{\frac{1}{1+\bar\eta_kg'\rbr{d_k}}}}^{-1},\\
    \bar\mu_{k+1}&=\frac{\bar\eta_{k+1}}{{\ccc\kappa}\gamma^2}\E_{(Q_1,Q_2,U)}\sbr{\gamma Q_1\cbr{h(\gamma Q_1,U)-g\rbr{d_k}}}+\bar\mu_k(1-\lambda\kappa^{-1}\bar{\eta}_{k+1}),\label{eq:SE-iter}\\
    \bar\sigma_{k+1}^2&=\frac{\bar\eta_{k+1}^2}{\kappa^2}\E_{(Q_1,Q_2,U)}\sbr{\cbr{h(\gamma Q_1,U)-g\rbr{d_k}}^2},
\end{align}
where 
\begin{align*}
    &d_k=\prox_{\bar\eta_k(\gamma)G}(\bar\mu_k\gamma Q_1+\sqrt{\kappa}\bar\sigma_kQ_2+\bar\eta_kh(\gamma Q_1,U)),\\
    &(Q_1,Q_2)\sim \mathcal{N}_2(\bzero,\bI_2),\qquad U\indep(Q_1,Q_2).
\end{align*}
Here, $U\in\R$ is a random variable independent of $(Q_1,Q_2)$. }
It is closely related to the linearized alternating direction method of multipliers \citep{rangan2016fixed}.

Second, we provide the following result for the estimator $\hat{\bbeta}^{\ccc\lambda}$, which is regarded as a limit of the parameter sequence generated by the GAMP algorithm.
The result follows the master theorem for the GAMP algorithm presented by \cite{feng2022unifying}.

\begin{lemma}
\label{lem:master}
Assume that $\{(\bX_i,Y_i)\}_{i=1}^n$ is an i.i.d. sample with $\bX_i\sim \mathcal{N}_p(\bzero,\bI_p)$ independent of $\bbeta$ and $\bvarepsilon=(\varepsilon_1,\ldots,\varepsilon_n)$. 
{\ccc
Let $(\mu_\lambda,\sigma_\lambda,\eta_\lambda)
\in\mathbb{R}\times(0,\infty)^2$ be a positive fixed point of the scalar state-evolution recursion \eqref{eq:SE-iter}.}
For $r\in[2,\infty)$, suppose that the empirical distributions $p^{-1}\sum_{j=1}^p\delta_{\sqrt{n}\beta_j}$ and $n^{-1}\sum_{i=1}^n\delta_{\varepsilon_i}$ converge in the $r$-Wasserstein sense to the distributions of $\Bar{\beta}$ and $\Bar{\varepsilon}$ with finite $r$-th order moment, respectively. Let $Z\sim \mathcal{N}(0,\gamma^2)$ be independent of $\Bar{\varepsilon}$, and $G,\tilde{G}\sim \mathcal{N}(0,1)$ be independent of $\Bar{\beta}$ and $Z$. 
{\ccc Let $p=p(n)$, and assume one of the following two regimes:
\begin{enumerate}
    \item[(i)] $\lambda>0\quad \text{and} \quad p/n \to \kappa\in(0,\infty),$
    \item[(ii)] $\lambda=0 \quad \text{and} \quad p/n \to \kappa\in(0,1).$
\end{enumerate}
Then, as $n\to\infty$,} for $\hat{\bbeta}^{\ccc\lambda}$ defined in \eqref{eq:estimator} such that $\|\hat{\bbeta}^{\ccc\lambda}\|=O_p(1)$ and any pseudo-Lipschitz function ${\psi}:\R^2\to\R$ and $\tilde{\psi}:\R^3\to\R$ of order $r$\footnote{A function $\psi:\R^m\to\R$ is said to be pseudo-Lipschitz of order $r$ if there exists a constant $L>0$ such that for any $t_0,t_1\in\R^m$, $\norm{\psi(t_0)-\psi(t_1)}_2\le L(1+\norm{t_0}_2^{r-1}+\norm{t_1}_2^{r-1})\norm{t_0-t_1}$.}, we have
{\ccc\begin{align*}
    \frac{1}{p}\sum_{j=1}^p{\psi}\rbr{\sqrt{n}\hat{\beta}_j^\lambda,\sqrt{n}\beta_j}&\overset{\rm a.s.}{\longrightarrow}\E\sbr{{\psi}\rbr{\mu_\lambda\Bar{\beta}+\sigma_\lambda G,\Bar{\beta}}},\\
    \frac{1}{n}\sum_{i=1}^n\tilde{\psi}\rbr{\bX_i^\top\hat\bbeta^\lambda,\bX_i^\top\bbeta,\varepsilon_i}&\overset{\rm a.s.}{\longrightarrow}\E\sbr{\tilde{\psi}\rbr{\prox_{\eta_\lambda\ell}\rbr{\mu_{\lambda} Z+\sqrt{\kappa}\sigma_{\lambda} \Tilde{G}},Z,\Bar{\varepsilon}}}.
\end{align*}}
\end{lemma}
{\ccc Remarkably, the $\lambda>0$ regime allows $\kappa>1$.
This is because, as shown in the proof below, the penalty term with $\lambda>0$ ensures that $\nabla_{\bb}^2\{\sum_{i=1}^n\ell(\bb;\bX_i,Y_i)+n\lambda\|\bb\|^2/2\}\succeq n\lambda\bI_p$ for any value of $\kappa$. 
In contrast, when
$\lambda=0$, such positive definiteness is no longer guaranteed, in
particular, in the overparameterized regime $\kappa>1$.}
% {\rcc [It is better to mention \textit{what} may fail in detail.]}
The convergence of pseudo-Lipschitz functions of order {\cc$r$} is equivalent to the convergence of the {\cc$r$}-Wasserstein distance between the empirical distributions of each coordinate, as established {\cc in Section 7.4 of \citep{feng2022unifying}}. 

{\cc
Note that the state evolution characterization at each GAMP iteration can be obtained by applying Theorem~4.1 in \citep{feng2022unifying}. The key step is then to transfer this iteration-wise characterization to the limit of the iterates, namely the estimator itself. The necessary argument for passing from finite iterations to the fixed-point limit is provided in the proof below.}

\begin{proof}[Proof of Lemma \ref{lem:master}] 
We prove this lemma by the general proof strategy discussed in Section 4.4 in \cite{feng2022unifying}.
It consists of three steps: (i) find a fixed point of the GAMP recursion; (ii) consider the stationary version of the GAMP recursion; (iii) show that the stationary version of the GAMP iterate converges to the penalized estimator $\hat\bbeta^\lambda$.
Following the strategy, the step (i) and the step (ii)
are simply achieved by the convexity and smoothness of the surrogate loss function by the assumption. 
Given the fixed point $\bar{\eta}_*$ of $\bar\eta_k$, and initial values $\bbeta^0, \bxi^0=\breve\bX\bbeta^0$, the stationary version of the GAMP algorithm takes the form, for each $k\in\N\cup\{0\}$,
\begin{align}
    \bbeta^{k+1}&=\frac{\bar\eta_*}{\kappa}\breve\bX^\top\cbr{\bY-g\rbr{\prox_{\bar\eta_*G}(\bxi^k+\bar\eta_*\bY)}}+\rbr{{{\ccc1-\kappa^{-1}\lambda\bar\eta_{*}}}}\bbeta^k\label{eq:GAMP-beta2}\\
    \bxi^{k+1}&=\breve\bX\bbeta^{k+1}-\bar\eta_*\cbr{\bY-g\rbr{\prox_{\bar\eta_*G}(\bxi^k+\bar\eta_*\bY)}}.\label{eq:GAMP-xi2}
\end{align}
In Step (iii), our goal is to ensure the algorithmic convergence of scaled GAMP iterates $\hat\bbeta^k=\bbeta^k/\sqrt{n}$ to the penalized estimator $\hat\bbeta^\lambda$ as $k\to\infty$. 
Denote $\cL(\bb)\equiv\sum_{i=1}^n\ell(\bb;\bX_i,Y_i)+{\ccc n\lambda\|\bb\|^2/2}$ for any $\bb=(b_1,\ldots,b_p)^\top\in\R^p$. By Taylor's theorem, we have
\begin{align*}
    \cL(\hat\bbeta^\lambda)
    &=\cL(\hat\bbeta^k)+\rbr{\hat\bbeta^\lambda-\hat\bbeta^k}^\top\nabla\cL(\hat\bbeta^k) \\
    & \quad +\frac{1}{2}\rbr{\hat\bbeta^\lambda-\hat\bbeta^k}^\top\nabla^2\cL\rbr{t\hat\bbeta^\lambda+(1-t)\hat\bbeta^k}\rbr{\hat\bbeta^\lambda-\hat\bbeta^k},
\end{align*}
for some $t\in(0,1)$. 

{\ccc\noindent$\bullet$ \textbf{$\lambda>0$ regime}:}
Since we have $\nabla^2\cL(\bb)=\sum_{i=1}^ng'(\bX_i^\top\bb)\bX_i\bX_i^\top+{\ccc n}\lambda\bI_p{\ccc\succeq n\lambda \bI_p}$ {\ccc for any $n$ and $p$ (i.e., $p/n$ can be larger than 1),  we have} that
\begin{align*}
    \cL(\hat\bbeta^\lambda)\ge\cL(\hat\bbeta^k)
    +\rbr{\hat\bbeta^\lambda-\hat\bbeta^k}^\top\nabla\cL(\hat\bbeta^k)
    +{\ccc n\lambda}\|\hat\bbeta^\lambda-\hat\bbeta^k\|^2.
\end{align*}
Using optimality of the estimator $\cL(\hat\bbeta^k)\ge\cL(\hat\bbeta^\lambda)$ and the Cauchy-Schwarz inequality as $(\hat\bbeta^\lambda-\hat\bbeta^k)^\top\nabla\cL(\hat\bbeta^k)\ge-\|{\nabla\cL(\hat\bbeta^k)\|}\|\hat\bbeta^\lambda-\hat\bbeta^k\|$ yields
\begin{align}
\label{eq:taylor-ineq}
    \|\hat\bbeta^\lambda-\hat\bbeta^k\|
    \le{\ccc\frac{1}{n\lambda}}\norm{\nabla\cL(\hat\bbeta^k)}.
\end{align}
Next, we consider controlling $\|\nabla\cL(\hat\bbeta^k)\|$. We have
\begin{align}
    \prox_{\bar\eta_*G}\rbr{\bxi^{k-1}+\bar\eta_*\bY}
    &=\bxi^{k-1}+\bar\eta_*\bY-\bar\eta_*g\rbr{\prox_{\bar\eta_*G}\rbr{\bxi^{k-1}+\bar\eta_*\bY}}\\
    &=\bxi^{k-1}-\bxi^k+\breve\bX\bbeta^k\\
    &=\bxi^{k-1}-\bxi^k+\bX\hat\bbeta^k,
\end{align}
by the definition of the proximal operator and \eqref{eq:GAMP-xi2}. Thus,
\begin{align*}
    \bbeta^k-\rbr{{{\ccc1-\kappa^{-1}\lambda\bar\eta_{*}}}}\bbeta^{k-1}
    &=\frac{\bar\eta_*}{\kappa}\breve\bX^\top\cbr{\bY-g\rbr{\prox_{\bar\eta_*G}(\bxi^{k-1}+\bar\eta_*\bY)}}\\
    &=\frac{\bar\eta_*}{\kappa}\breve\bX^\top\cbr{\bY-g\rbr{\bX\hat\bbeta^k+\bxi^{k-1}-\bxi^k}},
\end{align*}
by \eqref{eq:GAMP-beta2}. Using this and the triangle inequalities give
\begin{align*}
    \norm{\nabla\cL(\hat\bbeta^k)}
    &=\norm{\bX^\top \cbr{g(\bX^\top\hat\bbeta^k)-\bY}+{\ccc n}\lambda\hat\bbeta^k}\\
    &\le\norm{\bX^\top\cbr{\bY-g(\bX\hat\bbeta^k+\bxi^{k-1}-\bxi^k)}+{\ccc n}\lambda\hat\bbeta^{k-1}}+{\ccc n}\lambda\norm{\hat\bbeta^{k}-\hat\bbeta^{k-1}}\\
    &\quad+\norm{\bX^\top\cbr{g(\bX\hat\bbeta^k+\bxi^{k-1}-\bxi^k)-g(\bX^\top\hat\bbeta^k)}}\\
    &\le{\ccc n}\rbr{\frac{\kappa}{\bar\eta_*}+\lambda}\norm{\hat\bbeta^k-\hat\bbeta^{k-1}} \\
    & \quad +\norm{\bX}_{\rm op}\norm{g(\bX\hat\bbeta^k+\bxi^{k-1}-\bxi^k)-g(\bX^\top\hat\bbeta^k)}\\
    &\le{\ccc n}\rbr{\frac{\kappa}{\bar\eta_*}+\lambda}\norm{\hat\bbeta^k-\hat\bbeta^{k-1}}+L_{g}\norm{\bX}_{\rm op}\norm{\bxi^{k-1}-\bxi^k},
\end{align*}
where $L_g=\sup_zg'(z)$.
This establishes, combined with \eqref{eq:taylor-ineq},
\begin{align*}
    \|\hat\bbeta^\lambda-\hat\bbeta^k\|
    \le {\ccc \lambda}\cbr{\rbr{\frac{\kappa}{\bar\eta_*}+\lambda}\norm{\hat\bbeta^k-\hat\bbeta^{k-1}}
    +\frac{L_g}{n}\norm{\bX}_{\rm op}\norm{\bxi^{k-1}-\bxi^k}}.
\end{align*}
Finally, Lemma \ref{lem:gamp-cauchy-unified} implies
\begin{align*}
    \lim_{k\to\infty}\lim_{n\to\infty}\norm{\hat\bbeta^\lambda-\hat\bbeta^k}^{\ccc2}=_{\rm a.s.}0.
\end{align*}
{\ccc\noindent $\bullet$ \textbf{$\lambda=0$ regime}:} 
Lemma \ref{lem:str-convex} implies, for some non-increasing positive function $\omega:\R_+\to\R_+$,
\begin{align*}
    \ell(\hat\bbeta)\ge\ell(\hat\bbeta^k)
    +\rbr{\hat\bbeta-\hat\bbeta^k}^\top\nabla\ell(\hat\bbeta^k)
    +\frac{1}{2}n\omega\rbr{\max\cbr{\|\hat\bbeta\|,\|\hat\bbeta^k\|}}\|\hat\bbeta-\hat\bbeta^k\|^2,
\end{align*}
with probability at least $1-c_1\exp(-c_2n)$ where $c_1,c_2>0$ are some positive constants. 
Here, we use $\omega(\|t\hat\bbeta+(1-t)\hat\bbeta^k\|)\ge\omega(\max(\|\hat\bbeta\|,\|\hat\bbeta^k\|))$ since $\omega(\cdot)$ is non-increasing.
Using optimality of the estimator $\ell(\hat\bbeta^k)\ge\ell(\hat\bbeta)$ and the Cauchy-Schwarz inequality as $(\hat\bbeta-\hat\bbeta^k)^\top\nabla\ell(\hat\bbeta^k)\ge-\|{\nabla\ell(\hat\bbeta^k)\|}\|\hat\bbeta-\hat\bbeta^k\|$ yields, with probability at least $1-c_1\exp(-c_2n)$,
\begin{align*}
    \|\hat\bbeta-\hat\bbeta^k\|\le\frac{2}{\omega\rbr{\max\cbr{\|\hat\bbeta\|,\|\hat\bbeta^k\|}}}\norm{\frac{1}{n}\nabla\ell(\hat\bbeta^k)}
    \le\frac{2}{\omega(\|\hat\bbeta\|)\omega(\|\hat\bbeta^k\|)}\norm{\frac{1}{n}\nabla\ell(\hat\bbeta^k)},
\end{align*}
where the last inequality follows from the fact that $0<\omega(\cdot)<1$ and $\omega$ is non-increasing. Next, we consider controlling $\|\nabla\ell(\hat\bbeta^k)\|$. We have
\begin{align}
    \prox_{\bar\eta_*G}\rbr{\bxi^{k-1}+\bar\eta_*\bY}
    &=\bxi^{k-1}+\bar\eta_*\bY-\bar\eta_*g\rbr{\prox_{\bar\eta_*G}\rbr{\bxi^{k-1}+\bar\eta_*\bY}}\\
    &=\bxi^{k-1}-\bxi^k+\breve\bX\bbeta^k\\
    &=\bxi^{k-1}-\bxi^k+\bX\hat\bbeta^k,
\end{align}
by the definition of the proximal operator and \eqref{eq:GAMP-xi2}. Thus,
\begin{align*}
    \bbeta^k-\bbeta^{k-1}
    &=\frac{\bar\eta_*}{\kappa}\breve\bX^\top\cbr{\bY-g\rbr{\prox_{\bar\eta_*G}(\bxi^{k-1}+\bar\eta_*\bY)}}\\
    &=\frac{\bar\eta_*}{\kappa}\breve\bX^\top\cbr{\bY-g\rbr{\bX\hat\bbeta^k+\bxi^{k-1}-\bxi^k}},
\end{align*}
by \eqref{eq:GAMP-beta2}. Using this and triangle inequalities give
\begin{align*}
    \norm{\nabla\ell(\hat\bbeta^k)}
    &=\norm{\bX^\top \cbr{g(\bX^\top\hat\bbeta^k)-\bY}}\\
    &\le\norm{\bX^\top\cbr{\bY-g(\bX\hat\bbeta^k+\bxi^{k-1}-\bxi^k)}}\\
    &\quad+\norm{\bX^\top\cbr{g(\bX\hat\bbeta^k+\bxi^{k-1}-\bxi^k)-g(\bX^\top\hat\bbeta^k)}}\\
    &\le\frac{\kappa\sqrt{n}}{\bar\eta_*}\norm{\bbeta^k-\bbeta^{k-1}}+\norm{\bX}_{\rm op}\norm{g(\bX\hat\bbeta^k+\bxi^{k-1}-\bxi^k)-g(\bX^\top\hat\bbeta^k)}\\
    &\le\frac{\kappa\sqrt{n}}{\bar\eta_*}\norm{\bbeta^k-\bbeta^{k-1}}+L_{g}\norm{\bX}_{\rm op}\norm{\bxi^{k-1}-\bxi^k},
\end{align*}
where $L_g=\sup_zg'(z)$.
This establishes, with probability at least $1-c_1\exp(-c_2n)$,
\begin{align*}
    \|\hat\bbeta-\hat\bbeta^k\|
    \le c\cbr{\frac{\kappa}{\bar\eta_*}\norm{\hat\bbeta^k-\hat\bbeta^{k-1}}
    +\frac{C_g}{n}\norm{\bX}_{\rm op}\norm{\bxi^{k-1}-\bxi^k}},
\end{align*}
with $c={2}/\rbr{\omega(\|\hat\bbeta\|)\omega(\|\hat\bbeta^k\|)}$. Finally, Lemma \ref{lem:gamp-cauchy-unified} and the Borel-Cantelli lemma implies
\begin{align}
\label{eq:conv-beta^k}
    \lim_{k\to\infty}\lim_{n\to\infty}\norm{\hat\bbeta-\hat\bbeta^k}^{\ccc2}=_{\rm a.s.}0.
\end{align}
This completes the proof.
\end{proof}

Lemma \ref{lem:master} provides a statement on the convergence of the distance between 1-dimensional distributions, which is computed by averaging over the coordinates of $\hat{\bbeta}^{\ccc \lambda},\bbeta\in\R^p$. However, for statistical inference purposes, it is necessary to investigate the limiting behavior of the marginal distributions of each coordinate. To achieve this, we utilize a property of the surrogate loss function. 
\begin{lemma}
\label{lem:gen-cov}
For any invertible matrix $\bL\in\R^{p\times p}$ satisfying $\bSigma=\bL\bL^\top$, $\bL^\top\hat\bbeta$ minimizes the surrogate loss function in \eqref{eq:estimator} for the true coefficient $\bL^\top\bbeta$ and the covariate $\bL^{-1}\bX\sim \mathcal{N}_p(\bzero,\bI_p)$. 
\end{lemma}
\begin{proof}[Proof of Lemma \ref{lem:gen-cov}]
Since the surrogate loss depends on $\bX\in\R^p$ and $\bb\in\R^p$ only through their inner product $\bX^\top\bb$, we have $\ell(\bb;\bX,Y)=\ell(\bL^\top \bb;\bL^{-1}\bX,Y)$. 
\end{proof}

Using this in reverse, once we show Lemma~\ref{lem:master} for a design with identity covariance, we can rewrite the estimator corresponding to the unit covariance in $\bL^\top\hat{\bbeta}$ to obtain the general covariance result.

\subsubsection{Limit of Estimation and Classification Errors}
As a consequence of the master theorem of GAMP (Lemma \ref{lem:master}), we obtain the convergence limits of the mean squared error (MSE) and cosine similarity-like classification error. If we set $\psi(s,t)=(s-t)^2$, then Lemma \ref{lem:master} implies the MSE limit 
\begin{align*}
    \frac{1}{\kappa}\|\hat{\bbeta}^{\ccc\lambda}-\bbeta\|_2^2\cvas(\mu_{\ccc\lambda}-1)^2\mathbb{E}[\bar{\beta}^2 ]+\sigma_{\lambda}^2.
\end{align*}
By setting $\psi(s,t)=(s-\mu t)^2$, we have a corrected MSE limit 
\begin{align*}
    \frac{1}{\kappa}\|\hat{\bbeta}^{\ccc\lambda}-\mu_{\ccc\lambda}\bbeta\|_2^2\cvas\sigma_{\lambda}^2.
\end{align*}
For classification errors, by taking the ratio of $\psi(s,t)=st$ and $\psi(s,t)=t^2$, we can obtain the convergence limit of a cosine similarity-like measure 
\begin{align*}
    \frac{\bbeta^\top\hat{\bbeta}^{\ccc\lambda}}{\|\bbeta\|^2}\cvas\mu_{\ccc\lambda}.
\end{align*}

\subsection{GLM extension of the SLOE estimator}
\label{sec:SLOE-GLM}
This section provides an extension of the SLOE estimator \cite{yadlowsky2021sloe} from logistic regression to the GLM. First, it considers another representation of the state evolution parameters by using \textit{corrupted signal strength} $\gamma_c^2=\lim_{n\to\infty}\mathrm{Var}(\bX_1(n)^\top\hat\bbeta(n))$ instead of $\gamma^2$:
\begin{align}
\label{eq:cSE-system}
\begin{cases}
\kappa^2\sigma^2&=\eta^2\E_{(Q_1',Q_2',\Bar{Y})}\sbr{\rbr{\Bar{Y}-g\rbr{\prox_{\eta G}(Q_2')}}^2}, \\
0&=\E_{(Q_1',Q_2',\Bar{Y})}\sbr{Q_1'\rbr{\Bar{Y}-g\rbr{\prox_{\eta G}(Q_2')}}},\\
    1-\kappa&=\E_{(Q_1',Q_2',\Bar{Y})}\sbr{( 1+\eta g'\rbr{\prox_{\eta G}(Q_2')})^{-1}},
\end{cases}
\end{align}
where 
\begin{align}
\label{eq:nonlin-dist}
   \left(\begin{array}{c}Q_1' \\ Q_2' \end{array}\right) \sim \mathcal{N}_2\left(0, \left[ \begin{array}{cc}
    \mu^{-2}(\gamma_c^2-\kappa\sigma^2) & -\mu^{-1}(\gamma_c^2-\kappa\sigma^2) \\
    -\mu^{-1}(\gamma_c^2-\kappa\sigma^2) & \gamma_c^2
\end{array} \right] \right). 
\end{align}
Then, since $\bX_1^\top\hat\bbeta$ is computable unlike $\bX_1^\top\bbeta$, we can construct the SLOE-like estimator for $\gamma_c^2$ as follows.
\begin{align*}
    \hat\gamma_c^2&=\frac{1}{n}\sum_{i=1}^nS_i^2 - \rbr{\frac{1}{n}\sum_{i=1}^nS_i}^2,\\
    S_i&=\bX_i^\top\hat\bbeta+\frac{U_i}{1+g'(\bX_i^\top\hat\bbeta)}(Y_i-g(\bX_i^\top\hat\bbeta)),\\
    U_i&=-(\bX(\bX^\top \bD\bX)^{-1}\bX^\top)_{ii},
\end{align*}
where $\bD=\mathrm{diag}(g'(\bX_1^\top\hat\bbeta),\ldots,g'(\bX_n^\top\hat\bbeta))\in\R^{n\times n}$. 
This is a consistent estimator of $\gamma_c^2$.

\begin{proposition}
    Suppose that $\hat\bbeta$ exists and $\kappa\in(0,1)$ is fixed.
    Under Assumption \ref{asmp:basic} (A1)--(A2), $\hat\gamma_c^2\overset{\rm a.s.}{\to}\gamma_c^2$ as $n\to\infty$.
\end{proposition}
This follows from the direct application of the original proof of Proposition 2 in \cite{yadlowsky2021sloe} because of the form of the surrogate loss and the monotonicity of the link function $g(\cdot)$.

Unfortunately, the estimator $\hat\gamma_c^2$ is unstable when $p/n\to1$. This property inherits from two components: $\hat\bbeta$ and $(\bX^\top \bD\bX)^{-1}$ in $U_i$. 
First, for large $p/n$, $\hat\bbeta$ is likely to diverge as evident from \citep{candes2020phase,koriyama2025phase}. This phenomenon is observed in many estimators.
Second, the upper bound of $|\hat\gamma_c^2-\gamma_c^2|$ depends on a constant $1/\lambda_{\rm min}(\bX^\top \bD\bX)$. This constant diverges as $\kappa\uparrow1$.

{\cc
One can compute the solution to the fixed-point equations \eqref{eq:cSE-system} by iterating the following: 
\begin{align}
    \bar\eta_{k+1}&=\kappa\bar\eta_k\rbr{1-\E_{(Q_1,Q_2,U)}\sbr{\frac{1}{1+\bar\eta_kg'\rbr{d_k}}}}^{-1},\\
    \bar\mu_{k+1}&=\frac{\bar\eta_{k+1}}{{\ccc\kappa}\gamma^2}\E_{(Q_1,Q_2,U)}\sbr{\gamma Q_1\cbr{h(\gamma Q_1,U)-g\rbr{d_k}}}+\bar\mu_{k},\\
    \bar\sigma_{k+1}^2&=\frac{\bar\eta_{k+1}^2}{\kappa^2}\E_{(Q_1,Q_2,U)}\sbr{\cbr{h(\gamma Q_1,U)-g\rbr{d_k}}^2},
\end{align}
where 
\begin{align*}
    &d_k=\prox_{\bar\eta_k(\gamma)G}(\bar\mu_k\gamma Q_1+\sqrt{\kappa}\bar\sigma_kQ_2+\eta_kh(\gamma Q_1,U)),\\
    &\gamma^2=\frac{{\gamma_c^2-\kappa\bar\sigma_k^2}}{\bar\mu_k^2},\quad(Q_1,Q_2)\sim \mathcal{N}_2(\bzero,\bI_2),\qquad U\indep(Q_1,Q_2).
\end{align*}
Here, $U\in\R$ is a random variable independent of $(Q_1,Q_2)$ determined by the model.}

\subsection{GLM Form with \texorpdfstring{$h(Z,\bar{\varepsilon})$}{hzeps}}
\label{sec:h-example}
In this section, we discuss the other form of GLM used in (A3) in Assumption \ref{asmp:basic} as follows:
\begin{align}
    Y_i&=h\rbr{\bX_i^\top\bbeta,\varepsilon_i},\label{eq:h-SE}
\end{align}
where $h:\R\times\R\to\R$ is a deterministic function determined by the distribution of $Y\mid\bX$, and $\varepsilon_i\in\R$ is an error variable \textit{independent} of $\bX$.
While the GLM \eqref{eq:model} only specifies the conditional mean, \eqref{eq:h-SE} completely determines the distributional behavior of $Y$ depending on $\bX$. 
The major difference from the additive model \eqref{eq:glm_model2} used in our estimation is that the random variable $\varepsilon_i$ is independent of $\bX_i^\top$.

This model \eqref{eq:h-SE} is flexible in design and allows for an intuitive representation of GLMs.
We see several examples of $h(Z,\bar{\varepsilon})$ introduced in \eqref{eq:SE-system}.

\noindent $\bullet$ \textit{Bernoulli case (binary choice model)}. When we use a model $\Bar{Y}\mid Z\sim\mathrm{Ber}(g(Z))$ with some $g:\R\to[0,1]$, the inverse transformation method yields
\begin{align*}
    \bar{Y}=1\{g(Z)\ge\Bar{\varepsilon}\},\quad\Bar{\varepsilon}\sim\mathrm{Unif}[0,1].
\end{align*}

\noindent $\bullet$ \textit{Exponential case}. If $\Bar{Y}\mid Z\sim\mathrm{Exp}(g(Z))$ with some $g:\R\to(0,\infty)$, the inverse transformation method yields
\begin{align*}
    \Bar{Y}=-\frac{1}{g(Z)}\log(\Bar{\varepsilon}),\quad \Bar{\varepsilon}\sim\mathrm{Unif}[0,1].
\end{align*}

\noindent $\bullet$ \textit{Poisson case}. If $\Bar{Y}\mid Z\sim\mathrm{Pois}(g(Z))$ with some $g:\R\to(0,\infty)$, we have
\begin{align*}
    \Bar{Y}=\min\cbr{k\in\N\cup\{0\}\left|\sum_{l=1}^{k+1}\Bar{\varepsilon}_l> g(Z)\right.},\quad \bar{\varepsilon}_l\overset{\rm iid}{\sim}\mathrm{Exp}(1).
\end{align*}

\noindent $\bullet$ \textit{Gaussian case}. If $\Bar{Y}\mid Z\sim \mathcal{N}(g(Z),\sigma_{\Bar{\varepsilon}}^2)$ with some $g:\R\to\R$, we have
\begin{align*}
    \Bar{Y}=g(Z)+\Bar{\varepsilon},\quad \Bar{\varepsilon}\sim \mathcal{N}(0,\sigma_{\Bar{\varepsilon}}^2).
\end{align*}

As evident from the aforementioned examples, in the case of modeling using a one-parameter distribution, the distribution of $\Bar{\varepsilon}$ can be fully determined without any additional parameters.

\subsection{Sufficient Condition for Assumption \ref{asmp:g_estimation}} \label{sec:identification}

We discuss a sufficient condition for Assumption \ref{asmp:g_estimation} (A5), which assumes the uniqueness of the SS estimators provided in Section \ref{sec:moment_other}.
{\cc
\subsubsection{On the condition on $\tilde\Psi$.}

At first, we can see that the monotonicity of $h_1:\varsigma\mapsto\E[g(\varsigma Z)],\ Z\sim\mathcal{N}(0,1)$ is sufficient for the uniqueness of the solution to $\tilde\Psi=0$.

\begin{lemma}
\label{lem:unique-by-monotone}
    Suppose that $h_1:(0,\infty)\to\R$ is continuous and strictly monotonic, and that there exists $\gamma>0$ such that $\E\bar Y=h_1(\gamma)$.
    Then, for any $\epsilon>0$, we have

    \begin{align}
        \inf_{\varsigma:|\varsigma-\gamma|>\epsilon}|\tilde\Psi(\varsigma)|>0.
    \end{align}
\end{lemma}
\begin{proof}
    We assume that $h_1$ is strictly increasing without loss of generality.
    If $\varsigma\ge\gamma+\epsilon$, then $h_1(\varsigma)\ge h_1(\gamma+\epsilon)>h_1(\gamma)=\E\bar Y$. 
    This implies 
    \begin{align}
        |\E\bar Y-h_1(\varsigma)|=h_1(\varsigma)-\E\bar Y\ge h_1(\gamma+\epsilon)-h_1(\gamma)>0.
    \end{align}
    Also, if $\varsigma\le\gamma-\epsilon$, then $h_1(\varsigma)\le h_1(\gamma-\epsilon)<\E\bar Y$. This implies
    \begin{align}
        |\E\bar Y-h_1(\varsigma)|=\E\bar Y-h_1(\varsigma)\le h_1(\gamma)-h_1(\gamma-\epsilon)>0.
    \end{align}
    Thus, combining the above, we have
    \begin{align}
        \inf_{\varsigma:|\varsigma-\gamma|>\epsilon}|\tilde\Psi(\varsigma)|
        \ge\min\rbr{h_1(\gamma+\epsilon)-h_1(\gamma),h_1(\gamma)-h_1(\gamma-\epsilon)}>0.
    \end{align}
\end{proof}

Using this, we can construct two different sufficient conditions for (A5).

\begin{proposition}
\label{asmp:asymmetric}
If one of the following assumptions holds, then we have 
\begin{align}
    \inf_{\varsigma:|\varsigma-\gamma|>\epsilon}|\tilde\Psi(\varsigma)|>0.
\end{align}
for any $\epsilon>0$.
\begin{enumerate}
    \item[(i)] $g(\cdot)$ is twice-continuously differentiable and strictly convex. 
    {\bc\item[(ii)] The even part $g_0(x):=(g(x)+g(-x))/2$ of $g(\cdot)$ is strictly monotonic for $\R_+$.}
\end{enumerate} 
\end{proposition}

\begin{proof}
(i) By Assumption \ref{asmp:g_estimation} (A6), we have $h_1'(\varsigma)=\E[Zg'(\varsigma Z)]$. Then, Stein's lemma implies that $h_1'(\varsigma)=\varsigma\E[g''(\varsigma Z)]$.
From this, the strict convexity of $g(\cdot)$ implies $g''(\cdot)>0$, which means $h_1'(\varsigma)>0$ by $\varsigma>0$.
Thus, Lemma \ref{lem:unique-by-monotone} concludes the proof.

{\bc(ii) We write the cumulative distribution and density functions of the standard normal distribution by $\Phi:\R\to\R$ and $\phi:\R\to\R$, respectively.
Note that the density function of Gaussian distribution with mean zero is an even function. Since the product of even functions is even and the product of an even function and an odd function is odd, we have
\begin{align*}
    \E[g(\varsigma Q)]
    =\frac{1}{\varsigma}\int_{-\infty}^\infty g(x)\phi\rbr{\frac{x}{\varsigma}}dx
    =\int_{-\infty}^\infty g(\varsigma y)\phi\rbr{y}dy
    =2\int_{0}^\infty g_0(\varsigma y)\phi\rbr{y}dy,
\end{align*}
where the second identity is from a change of variables $y=x/\varsigma$. Thus, {\cc for} any $\varsigma_1>\varsigma_0>0$, we can say that
\begin{align*}
    \E[g(Z_{\varsigma_1})]-\E[g(Z_{\varsigma_0})]
    =2\int_{0}^\infty \cbr{g_0(\varsigma_1 y)-g_0(\varsigma_0 y)}\phi\rbr{y}dy.
\end{align*}
is strictly positive or negative by $\phi(\cdot)>0$. Therefore, if the even part $g_0(\cdot)$ of $g(\cdot)$ is monotone on $\R_+$, then $\E[g(\varsigma Q)]$ is monotone in $\varsigma>0$. 
Then, Lemma \ref{lem:unique-by-monotone} concludes the proof.}
\end{proof}

{\cc
\begin{example}[Poisson regression]
    The MLE for Poisson regression model has $g(t)=e^t$, and thus satisfies the differentiability and convexity assumption.
    As a result, it allows for the unique solution to the moment equation.
\end{example}
}

\subsubsection{On the condition on $\Psi$}

We first define the following variables:
\[
m_2 := \E[\bar Y^2], \quad
m_4 := \E[\bar Y^4], \quad
h_2(\varsigma) := \E[g(\varsigma Z)^2], \quad
h_4(\varsigma) := \E[g(\varsigma Z)^4],
\]
where $Z\sim\mathcal N(0,1)$.
Recall that the population moment equations are given by
\[
\Psi(\varsigma,\varsigma_e)
=
\begin{pmatrix}
\Psi_1(\varsigma,\varsigma_e)\\
\Psi_2(\varsigma,\varsigma_e)
\end{pmatrix}
=
\begin{pmatrix}
m_2-h_2(\varsigma)-\varsigma_e^2\\
m_4-h_4(\varsigma)-6h_2(\varsigma)\varsigma_e^2-3\varsigma_e^4
\end{pmatrix}.
\]

\begin{lemma}
\label{lem:psi_separation}
Assume that the map
\[
\varsigma \mapsto K(\varsigma):=h_4(\varsigma)-3h_2(\varsigma)^2
\]
is strictly monotone on an interval containing $\gamma$.
Then, for any $\epsilon>0$, there exists a constant $c_\epsilon>0$ such that
\[
\inf_{(\varsigma,\varsigma_e): |\varsigma-\gamma|+|\varsigma_e-\sigma_e|>\epsilon}
\|\Psi(\varsigma,\varsigma_e)\| \ \ge\ c_\epsilon.
\]
\end{lemma}

\begin{proof}
We first characterize the solutions to the population equations $\Psi(\varsigma,\varsigma_e)=0$.
From $\Psi_1(\varsigma,\varsigma_e)=0$, we obtain $\varsigma_e^2 = m_2 - h_2(\varsigma).$
Substituting this identity into $\Psi_2(\varsigma,\varsigma_e)=0$ yields
\[
m_4 - 3m_2^2 = h_4(\varsigma)-3h_2(\varsigma)^2 = K(\varsigma).
\]
Since $K$ is strictly monotone by assumption, the above equation admits the unique solution
$\varsigma=\gamma$, and consequently $\varsigma_e^2 = m_2-h_2(\gamma)=\sigma_e^2$.
Hence $(\gamma,\sigma_e)$ is the unique solution to $\Psi(\varsigma,\varsigma_e)=0$.

Next, note that $h_2$ and $h_4$ are continuous functions of $\varsigma$ under the assumed moment
conditions on $g$, and therefore $\Psi$ is continuous on $(0,\infty)^2$.
Fix any compact set $\Theta\subset(0,\infty)^2$ containing $(\gamma,\sigma_e)$, and define
\[
\Theta_\epsilon
:=
\{(\varsigma,\varsigma_e)\in\Theta:\ |\varsigma-\gamma|+|\varsigma_e-\sigma_e|\ge\epsilon\}.
\]
Then $\Theta_\epsilon$ is compact and does not contain $(\gamma,\sigma_e)$.
By continuity, the function $(\varsigma,\varsigma_e)\mapsto\|\Psi(\varsigma,\varsigma_e)\|$
attains its minimum on $\Theta_\epsilon$.
Since $(\gamma,\sigma_e)$ is the unique zero of $\Psi$, this minimum must be strictly positive,
which proves the claim.
\end{proof}

We next provide concrete classes of link functions $g$ for which the strict monotonicity of $K(\varsigma)=h_4(\varsigma)-3h_2(\varsigma)^2$ can be verified analytically.
These examples demonstrate that the condition in Lemma~\ref{lem:psi_separation} is satisfied by a broad and practically relevant family of GLMs.

\begin{proposition}[Odd-power polynomial links]
\label{prop:odd_polynomial}
Let
\[
g(x)=a+b x^{q},
\]
where $a\in\mathbb R$, $b\neq 0$, and $q\ge 3$ is an odd integer.
Then the map $\varsigma\mapsto K(\varsigma)$ is strictly increasing on $(0,\infty)$.
Consequently, the population moment equations $\Psi(\varsigma,\varsigma_e)=0$ admit a unique solution $(\gamma,\sigma_e)$, and the uniqueness assertion in Lemma~\ref{lem:psi_separation} holds.
\end{proposition}

\begin{proof}
Let $\varsigma Q=\varsigma Z$ with $Z\sim\mathcal N(0,1)$.
Since $q$ is odd, the symmetry of the Gaussian distribution implies
\[
\E[\varsigma Q^{q}]=\E[\varsigma Q^{3q}]=0.
\]
A direct expansion yields
\begin{align*}
h_2(\varsigma)
&=\E[(a+b\varsigma Q^{q})^2]
=a^2+b^2\E[Z^{2q}] \varsigma^{2q},\\
h_4(\varsigma)
&=\E[(a+b\varsigma Q^{q})^4]
=a^4+6a^2b^2\E[Z^{2q}] \varsigma^{2q}
+b^4\E[Z^{4q}] \varsigma^{4q}.
\end{align*}
Therefore, we have
\[
K(\varsigma)
=h_4(\varsigma)-3h_2(\varsigma)^2
=-2a^4
+b^4\Big(\E[Z^{4q}]-3(\E[Z^{2q}])^2\Big)\varsigma^{4q}.
\]
For $q\ge 3$, the coefficient
\[
\E[Z^{4q}]-3(\E[Z^{2q}])^2=(4q-1)!!-3\big((2q-1)!!\big)^2
\]
is strictly positive.
Hence $K(\varsigma)$ is a strictly increasing function of $\varsigma$ on $(0,\infty)$.
\end{proof}

\begin{proposition}[Hyperbolic-sine links]
\label{prop:sinh}
Let
\[
g(x)=a+b\sinh(cx),
\]
where $a\in\mathbb R$ and $bc\neq 0$.
Then the map $\varsigma\mapsto K(\varsigma)$ is strictly increasing on $(0,\infty)$, and thus the conclusion of Lemma~\ref{lem:psi_separation} holds.
\end{proposition}

\begin{proof}
Let $\varsigma Q=\varsigma Z$ with $Z\sim\mathcal N(0,1)$ and set $t=c\varsigma$.
Since $\sinh$ is an odd function, we have
\[
\E[\sinh(tZ)]=\E[\sinh(tZ)^3]=0.
\]
Using the identities
\[
\sinh^2 x=\frac{\cosh(2x)-1}{2},
\qquad
\sinh^4 x=\frac{\cosh(4x)-4\cosh(2x)+3}{8},
\]
together with $\E[\cosh(\alpha Z)]=\exp(\alpha^2/2)$, we obtain
\begin{align*}
\E[\sinh(tZ)^2]
&=\frac{e^{2t^2}-1}{2},\\
\E[\sinh(tZ)^4]
&=\frac{e^{8t^2}-4e^{2t^2}+3}{8}.
\end{align*}
Consequently, we have
\[
\E[\sinh(tZ)^4]-3\E[\sinh(tZ)^2]^2
=\frac{e^{8t^2}-6e^{4t^2}+8e^{2t^2}-3}{8}
=:K_{\sinh}(t).
\]
Writing $s=t^2>0$, we compute
\[
\frac{\d}{\d s}K_{\sinh}(t)
=e^{8s}-3e^{4s}+2e^{2s}
=e^{2s}(e^{2s}-1)^2(e^{2s}+2)>0,
\]
which shows that $K_{\sinh}(t)$ is strictly increasing in $t>0$.
Since $t=c\varsigma$ with $|c|>0$, it follows that $K(\varsigma)$ is strictly increasing in $\varsigma$.
\end{proof}

\subsubsection{Necessary condition}
\label{sec:necessary_A5_both}

In this subsection, we provide a necessary-condition perspective that links the behavior of
the first moment equation $\tilde\Psi(\cdot)=0$ and the higher-order moment system $\Psi(\cdot,\cdot)=\bzero$.
The key observation is that a certain symmetry of the link function prevents identification
of the signal-strength parameter $\gamma$ for \emph{both} estimators.

{\ccc
\begin{lemma}
Assume the Gaussian-input limit model \eqref{eq:limit-additive}.

\begin{enumerate}
    \item[(i)] Suppose that
    \[
    g(x)=c+f(x),
    \]
    where $c\in\mathbb{R}$ and $f$ is odd. Then the first-moment map
    $\tilde{\Psi}(\varsigma)$ is constant in $\varsigma$ and does not identify
    $\gamma$.

    \item[(ii)] Independently of part~(i), suppose that $\sigma_e>0$, $h_2$ is
    continuous in a neighborhood of $\gamma$, and
    \[
    K(\varsigma)=h_4(\varsigma)-3h_2(\varsigma)^2
    \]
    is constant on an open interval containing $\gamma$. Then the
    second- and fourth-moment system does not locally identify
    $(\gamma,\sigma_e)$.
\end{enumerate}
\end{lemma}
}

\begin{proof}
We first prove (i).
Since $f$ is odd and $\varsigma Q\sim\mathcal N(0,\varsigma^2)$ is symmetric about zero, we have
\[
\E[f(\varsigma Q)]=0 \quad \text{for all }\varsigma>0.
\]
Therefore,
\[
\E[g(\varsigma Q)]
=\E[c+f(\varsigma Q)]
=c,
\]
which is independent of $\varsigma$.
Since $\E[\bar Y]=\E[g(\gamma Q)]$ under the model, it follows that
$\tilde\Psi(\varsigma)\equiv 0$ for all $\varsigma>0$, and hence
$\tilde\Psi$ cannot distinguish $\gamma$ from any other value.

We next prove (ii).
The population system
$\Psi(\varsigma,\varsigma_e)=0$ is equivalent to
\[
K(\varsigma)=m_4-3m_2^2,
\qquad
\varsigma_e^2=m_2-h_2(\varsigma),
\]
where $m_2=\E[\bar Y^2]$ and $m_4=\E[\bar Y^4]$.
If $K(\varsigma)$ is constant on an interval $I$ containing $\gamma$, then
$K(\varsigma)=K(\gamma)=m_4-3m_2^2$ for all $\varsigma\in I$.
Hence, for each such $\varsigma$, choosing
\[
\varsigma_e(\varsigma):=\sqrt{m_2-h_2(\varsigma)}
\]
(formally, whenever $m_2-h_2(\varsigma)\ge 0$) yields
$\Psi(\varsigma,\varsigma_e(\varsigma))=0$.
Thus the system admits infinitely many solutions arbitrarily close to $(\gamma,\sigma_e)$,
which rules out the separation property required by Assumption~(A5).
\end{proof}
Intuitively, suppose that $Y \mid \bm{X}$ follows a distribution that is fully determined by its mean $g(\bm{X}^\top \bm{\beta})$. 
If $g$ satisfies the symmetry condition $g(-t)=1-g(t),$
as in the logistic link, then under the zero-symmetry of the distribution of $\bm{X}^\top \bm{\beta}$, we have $\mathbb{E}[g(\bm{X}^\top \bm{\beta})]=1/2.$
Hence, the marginal distribution of $Y$ is fixed (at $\mathrm{Bernoulli}(1/2)$ in the logistic regression case) and does not depend on the signal strength $\gamma^2$. 
Consequently, the strength of the signal has no effect on the distribution of $Y$, making it impossible to estimate the signal strength based on $Y$ alone.
}

\section{Proofs of Main Results}
{\ccc
In what follows, for notational convenience, we suppress the dependence on $\lambda$ in $(\mu_\lambda,\sigma_\lambda,\eta_\lambda)$ and $\hat{\bbeta}^{\lambda}$ whenever it is clear from the context.

}

{\ccc
\begin{proof}[Proof of Proposition~\ref{prop:SE-GLM}]
We first derive a unified fixed-point system for $\lambda\ge0$ under the
isotropic design.

For each $y\in\mathbb{R}$, define the scalar surrogate loss
\[
\ell_y(t):=G(t)-yt.
\]
For any $u\in\mathbb{R}$ and $\eta>0$, completing the square gives
\begin{align}
\operatorname{prox}_{\eta\ell_y}(u)
&=
\argmin_{d\in\mathbb{R}}
\left\{
\eta G(d)-\eta yd+\frac12(d-u)^2
\right\}
\notag\\
&=
\argmin_{d\in\mathbb{R}}
\left\{
\eta G(d)
+
\frac12\{d-(u+\eta y)\}^2
\right\}
\notag\\
&=
\operatorname{prox}_{\eta G}(u+\eta y).
\label{eq:prox-loss-shift}
\end{align}
Thus, the output-side scalar variable appearing in
Lemma~\ref{lem:master} can equivalently be written using the proximal
operator of $G$.

Fix $\lambda\ge0$. Let
\[
(\mu_\lambda,\sigma_\lambda,\eta_\lambda)
\in\mathbb{R}\times(0,\infty)^2
\]
be the state-evolution fixed point appearing in
Lemma~\ref{lem:master}. Let
\[
Z:=\gamma Q_1,
\qquad
\bar Y:=h(Z,\bar\varepsilon),
\]
where $(Q_1,Q_2)\sim\mathcal{N}_2(\bm0,\bm I_2)$ and
$\bar\varepsilon$ is independent of $(Q_1,Q_2)$. Define
\begin{align}
D_\lambda
&:=
\operatorname{prox}_{\eta_\lambda G}
\left(
\mu_\lambda Z
+
\sqrt{\kappa}\sigma_\lambda Q_2
+
\eta_\lambda\bar Y
\right),
\label{eq:D-lambda-definition}\\
R_\lambda
&:=
\bar Y-g(D_\lambda).
\label{eq:R-lambda-definition}
\end{align}
All expectations below are taken with respect to the joint distribution of
$(Q_1,Q_2,\bar\varepsilon)$.

By Lemma~\ref{lem:master}, the limiting parameters of the estimator
are a fixed point of the scalar state-evolution recursion displayed in
Appendix~A.3. In the notation above, that recursion is
\begin{align}
\eta_{k+1}
&=
\kappa\eta_k
\left[
\lambda\eta_k+1
-
\mathbb{E}
\left\{
\frac{1}{1+\eta_k g'(D_k)}
\right\}
\right]^{-1},
\label{eq:eta-recursion-proof}\\
\mu_{k+1}
&=
\frac{\eta_{k+1}}{\kappa\gamma^2}
\mathbb{E}
\left[
Z\{\bar Y-g(D_k)\}
\right]
+
\left(
1-\frac{\lambda\eta_{k+1}}{\kappa}
\right)\mu_k,
\label{eq:mu-recursion-proof}\\
\sigma_{k+1}^2
&=
\frac{\eta_{k+1}^2}{\kappa^2}
\mathbb{E}
\left[
\{\bar Y-g(D_k)\}^2
\right],
\label{eq:sigma-recursion-proof}
\end{align}
where
$
D_k
=
\operatorname{prox}_{\eta_kG}
\left(
\mu_kZ+\sqrt{\kappa}\sigma_kQ_2+\eta_k\bar Y
\right).
$

At the fixed point, it holds that
\[
(\mu_{k+1},\sigma_{k+1},\eta_{k+1})
=
(\mu_k,\sigma_k,\eta_k)
=
(\mu_\lambda,\sigma_\lambda,\eta_\lambda),
\]
and hence $D_k=D_\lambda$ and
$\bar Y-g(D_k)=R_\lambda$.

First, substituting the fixed point into
\eqref{eq:sigma-recursion-proof} gives
$
\sigma_\lambda^2
=
\frac{\eta_\lambda^2}{\kappa^2}
\mathbb{E}[R_\lambda^2].
$
Equivalently, we obtain
\begin{align}
\kappa^2\sigma_\lambda^2
=
\eta_\lambda^2
\mathbb{E}
\left[
\{\bar Y-g(D_\lambda)\}^2
\right].
\label{eq:fixed-sigma-proof}
\end{align}

Second, substituting the fixed point into
\eqref{eq:mu-recursion-proof} yields
\begin{align*}
\mu_\lambda
&=
\frac{\eta_\lambda}{\kappa\gamma^2}
\mathbb{E}[ZR_\lambda]
+
\left(
1-\frac{\lambda\eta_\lambda}{\kappa}
\right)\mu_\lambda.
\end{align*}
Rearranging it yields
$
\frac{\lambda\eta_\lambda}{\kappa}\mu_\lambda
=
\frac{\eta_\lambda}{\kappa\gamma^2}
\mathbb{E}[ZR_\lambda].
$
Since $\eta_\lambda>0$ and $\gamma>0$, this is equivalent to
\begin{align}
\gamma^2\lambda\mu_\lambda
=
\mathbb{E}
\left[
Z\{\bar Y-g(D_\lambda)\}
\right].
\label{eq:fixed-mu-proof}
\end{align}

Finally, substituting the fixed point into
\eqref{eq:eta-recursion-proof} gives
\[
\eta_\lambda
=
\kappa\eta_\lambda
\left[
\lambda\eta_\lambda+1
-
\mathbb{E}
\left\{
\frac{1}{1+\eta_\lambda g'(D_\lambda)}
\right\}
\right]^{-1}.
\]
Since $\eta_\lambda>0$, division by $\eta_\lambda$ and rearrangement yield
\begin{align}
1-\kappa+\lambda\eta_\lambda
=
\mathbb{E}
\left[
\{1+\eta_\lambda g'(D_\lambda)\}^{-1}
\right].
\label{eq:fixed-eta-proof}
\end{align}
Equations
\eqref{eq:fixed-sigma-proof}--\eqref{eq:fixed-eta-proof}
constitute the unified fixed-point system for $\lambda\ge0$.

For the regime (ii), set $\lambda=0$ and write
\[
(\mu,\sigma,\eta)
:=
(\mu_0,\sigma_0,\eta_0),
\qquad
D:=D_0.
\]
Then \eqref{eq:fixed-sigma-proof}--\eqref{eq:fixed-eta-proof} reduce to
\[
\begin{cases}
\kappa^2\sigma^2
=
\eta^2
\mathbb{E}
\left[
\{\bar Y-g(D)\}^2
\right],
\\[1mm]
0
=
\mathbb{E}
\left[
Z\{\bar Y-g(D)\}
\right],
\\[1mm]
1-\kappa
=
\mathbb{E}
\left[
\{1+\eta g'(D)\}^{-1}
\right],
\end{cases}
\]
where
$
D
=
\operatorname{prox}_{\eta G}
\left(
\mu Z+\sqrt{\kappa}\sigma Q_2+\eta\bar Y
\right).
$
This is precisely the state-evolution system
\eqref{eq:SE-system} in the isotropic case.

To extend the unpenalized result to a general positive-definite covariance
matrix $\bm\Sigma$, write
$
\bm\Sigma=\bm L\bm L^\top
$
and define
$
\widetilde{\bm X}_i:=\bm L^{-1}\bm X_i,
\bm\theta:=\bm L^\top\bm\beta$, and 
$
\hat{\bm\theta}:=\bm L^\top\hat{\bm\beta}.
$
Then, we have
$
\widetilde{\bm X}_i\sim\mathcal N_p(\bm0,\bm I_p)
$
and
$
\widetilde{\bm X}_i^\top\bm\theta
=
\bm X_i^\top\bm\beta.
$
Moreover, for every candidate vector $\bm b$,
\[
\ell(\bm b;\bm X_i,Y_i)
=
\ell(\bm L^\top\bm b;\widetilde{\bm X}_i,Y_i).
\]
Therefore, by Lemma~\ref{lem:gen-cov},
$\hat{\bm\theta}$ is the isotropic-design surrogate estimator with true
coefficient $\bm\theta$. The signal strength is preserved:
\[
\|\bm\theta\|^2
=
\bm\beta^\top\bm L\bm L^\top\bm\beta
=
\bm\beta^\top\bm\Sigma\bm\beta
\to
\gamma^2.
\]
Hence the same scalar fixed-point system applies, proving
Proposition~\ref{prop:SE-GLM} for general $\bm\Sigma$.
\end{proof}
}

\begin{proof}[Proof of Proposition \ref{prop:asym_norm_betahat}]
Lemma \ref{lem:gen-cov} implies that, for any $j=1,\ldots,p$,
\begin{align*}
    \frac{\hat{\beta}_j-\mu\beta_j}{{\sigma}/\tau_j}=\frac{\hat{\theta}_j-\mu\theta_j}{{\sigma}},
\end{align*}
where we define
\begin{align}
\label{eq:theta}
    \btheta=\bL^\top\bbeta,\quad\hat\btheta=\bL^\top\hat\bbeta,
\end{align}
by a Cholesky factorization $\bSigma=\bL\bL^\top$. We have $\hat\theta_j=\tau_j\hat\beta_j$ and $\theta_j=\tau_j\beta_j$ by the rearrangement of indices. Here, define
\begin{align}
\label{eq:SE-n}
    \mu_n=\frac{\hat\btheta^\top\btheta}{\norm{\btheta}^2},\qquad \mbox{and} \qquad
    \sigma_n^2=\frac{1}{\kappa}\|\hat\btheta-\mu_n\btheta\|^2.
\end{align}
Then, we have
\begin{align}
\label{eq:decom:mu_n-mu}
    \sqrt{n}\frac{\hat\theta_j-\mu\theta_j}{\sigma}
    =\sqrt{n}\frac{\hat\theta_j-\mu_n\theta_j}{\sigma_n}\frac{\sigma_n}{\sigma}
    +\sqrt{n}\frac{(\mu_n-\mu)\theta_j}{\sigma}.
\end{align}
Lemma \ref{lem:theta-rot} gives us
\begin{align}
\label{eq:deq-PZ}
    \frac{\hat\btheta-\mu_n\btheta}{\sigma_n}
    \overset{\rm d}{=}\frac{\bP_{\btheta}^\perp\bZ}{\|\bP_{\btheta}^\perp\bZ\|},
\end{align}
where $\bZ=(Z_1,\ldots,Z_p)\sim\mathcal{N}_p(\bzero,\bI_p)$ {\cc and $\bP_{\btheta}^\perp=\bI_p-\btheta\btheta^\top/\|\btheta\|^2$}. Triangle inequalities yield that
\begin{align*}
    \frac{\norm{\bZ}}{\sqrt{p}}-\frac{|\btheta^\top\bZ|}{\sqrt{p}\norm{\btheta}}
    \le\frac{\|\bP_{\btheta}^\perp\bZ\|}{\sqrt{p}}
    \le\frac{\norm{\bZ}}{\sqrt{p}}+\frac{|\btheta^\top\bZ|}{\sqrt{p}\norm{\btheta}}.
\end{align*}
Since $|\btheta^\top\bZ|/(\sqrt{p}\norm{\btheta})\as0$ and $\norm{\bZ}/\sqrt{p}\as1$, we have $\|\bP_{\btheta}^\perp\bZ\|/\sqrt{p}\as1$. Then, this fact and \eqref{eq:deq-PZ} imply that, for each $j\in[p],$
\begin{align}
\label{eq:deq-tstat}
    \sqrt{n}\frac{\hat\theta_j-\mu_n\theta_j}{\sigma_n}
    \overset{\rm d}{=}\sigma_j Q+o_p(1),\qquad\sigma_j^2=1-\frac{\theta_j^2}{\norm{\btheta}^2},
\end{align}
where $Q\sim\mathcal{N}(0,1)$. 
Here we use the fact that the covariance matrix of $\bP_{\btheta}^\perp \bZ$ is $\bP_{\btheta}^\perp \bP_{\btheta}^\perp=\bI_p-\btheta\btheta^\top/\|\btheta\|^2$.
Thus, the facts that $\mu_n\as\mu$ and $\sigma_n^2\as\sigma^2$ by Lemma \ref{lem:SE-approx} conclude the proof of the first result.
{\ccc
Since the ridge penalty $\bb\mapsto\lambda/2\|\bb\|^2$ is orthogonal invariant in that $\|\bU\bb\|^2=\bb^\top\bU^\top\bU\bb=\|\bb\|^2$ for any orthogonal matrix $\bU\in\R^{p\times p}$, the same proof holds for $\lambda>0$ and $\bSigma=\bI_p$.}

{\ccc
We next establish the uniform convergence over $j\in[p]$.  
Define
\[
S_{j,n}:=
\frac{\hat\theta_j-\mu\theta_j}{\sigma/\sqrt n},
\quad
T_{j,n}:=
\frac{\hat\theta_j-\mu_n\theta_j}{\sigma_n/\sqrt n},
\quad
a_n:=\frac{\sigma_n}{\sigma},
\quad
b_{j,n}
:=
\frac{(\mu_n-\mu)\theta_j}{\sigma/\sqrt n}.
\]
Then we have $S_{j,n}=a_nT_{j,n}+b_{j,n}.$
By Lemma~\ref{lem:SE-approx}, $a_n\pconv1$ and $\mu_n\pconv\mu.$
Moreover, the uniform coordinate condition implies
\[
\max_{j\in[p]}|b_{j,n}|
\le
\frac{|\mu_n-\mu|}{\sigma}
\max_{j\in[p]}\sqrt n\,|\theta_j|
\pconv0.
\]
For any fixed $\epsilon\in(0,1/2)$, define an event
\[
E_{n,\epsilon}
:=\left\{
|a_n-1|\le\epsilon
\right\}
\cap
\left\{
\max_{j\in[p]}|b_{j,n}|\le\epsilon
\right\}.
\]
Then it follows that $\pr(E_{n,\epsilon}^c)\to0$.

For $t\in\mathbb R$, define the deterministic lower and upper thresholds
\[
\underline r_\epsilon(t)
:=
\inf_{\substack{|a-1|\le\epsilon\\|b|\le\epsilon}}
\frac{t-b}{a},
\qquad
\overline r_\epsilon(t)
:=
\sup_{\substack{|a-1|\le\epsilon\\|b|\le\epsilon}}
\frac{t-b}{a}.
\]
Since $a_n>0$ on $E_{n,\epsilon}$, we have, simultaneously for every
$j\in[p]$,
\[
\{T_{j,n}<\underline r_\epsilon(t)\}
\cap E_{n,\epsilon}
\subseteq
\{S_{j,n}<t\}
\cap E_{n,\epsilon}
\subseteq
\{T_{j,n}<\overline r_\epsilon(t)\}
\cap E_{n,\epsilon}.
\]
Therefore,
\begin{align}
\pr(T_{j,n}<\underline r_\epsilon(t)))
-\pr(E_{n,\epsilon}^c)
&\le
\pr(S_{j,n}<t)
\le
\pr(T_{j,n}<\overline r_\epsilon(t)))
+\pr(E_{n,\epsilon}^c).
\label{eq:uniform-sandwich}
\end{align}

We next control the perturbation of the Gaussian distribution function.
For any $a\in[1-\epsilon,1+\epsilon]$, $|b|\le\epsilon$, and
$t\in\mathbb R$,
\begin{align*}
\left|
\Phi\left(\frac{t-b}{a}\right)-\Phi(t)
\right|
&\le
\left|
\Phi\left(\frac{t-b}{a}\right)
-
\Phi\left(\frac{t}{a}\right)
\right|
+
\left|
\Phi\left(\frac{t}{a}\right)-\Phi(t)
\right|\\
&\le
\frac{|b|}{a\sqrt{2\pi}}
+
C|a-1|\\
&\le
C_0\epsilon,
\end{align*}
where $C,C_0>0$ are universal constants. Hence,
\[
\sup_{t\in\mathbb R}
\max\left\{
\left|\Phi(\underline r_\epsilon(t))-\Phi(t)\right|,
\left|\Phi(\overline r_\epsilon(t))-\Phi(t)\right|
\right\}
\le C_0\epsilon.
\]

Combining this bound with \eqref{eq:uniform-sandwich} yields
\begin{align*}
\max_{j\in[p]}
\sup_{t\in\mathbb R}
\left|
\pr(S_{j,n}<t)-\Phi(t)
\right|
\le
\Delta_n
+
\pr(E_{n,\epsilon}^c)
+
C_0\epsilon,
\end{align*}
where $\Delta_n:=\max_{j\in[p]}
\sup_{t\in\mathbb R}
\left|\pr(T_{j,n}<t)-\Phi(t)
\right|.$
Indeed, the upper bound in \eqref{eq:uniform-sandwich} gives
\begin{align*}
\pr(S_{j,n}<t)-\Phi(t)
&\le
\pr\left\{
T_{j,n}<\overline r_\epsilon(t)
\right\}
-\Phi\left(\overline r_\epsilon(t)\right)\\
&\quad+
\Phi\left(\overline r_\epsilon(t)\right)-\Phi(t)
+
\pr(E_{n,\epsilon}^c),
\end{align*}
and the lower bound gives the analogous inequality for
$\Phi(t)-\pr(S_{j,n}<t)$.

For every fixed $\epsilon>0$, letting $n,p\to\infty$ with
$p/n\to\kappa$ gives
\[
\limsup_{n,p\to\infty}
\max_{j\in[p]}
\sup_{t\in\mathbb R}
\left|
\pr(S_{j,n}<t)-\Phi(t)
\right|
\le C_0\epsilon,
\]
because $\Delta_n\to0$ and
$\pr(E_{n,\epsilon}^c)\to0$, where 
$\Delta_n\to0$ follows from \eqref{eq:deq-tstat} and $|1-\max_{j\in[p]}\theta_j^2/\|\btheta\|^2|\to0$.

Finally, letting $\epsilon\downarrow0$ completes the proof of the uniform convergence.
}
\end{proof}

\begin{proof}[Proof of Theorem \ref{thm:consistent_moment}] \textbf{Case of \eqref{eq:gamma-hat}}. The law of large numbers yields
\begin{align}
\label{eq:vdv1}
    &\sup_{\varsigma^2>0}\abs{\tilde\Psi_n(\varsigma^2)-\tilde\Psi(\varsigma^2)}=\abs{\frac{1}{n}\sum_{i=1}^nY_i-\E[\bar{Y}]}\overset{\rm a.s.}{\longrightarrow}0,
\end{align}
Thus, \eqref{eq:vdv1}, Assumption (A5), and $\tilde\Psi(\gamma)=0$ imply $\hat{\gamma}^2\overset{\rm a.s.}{\longrightarrow}\gamma^2$ by Theorem 5.9 in \cite{van2000asymptotic}.

{\ccc
Since $\hat\gamma$ is defined to be a exact root of $\tilde\Psi_n$, there exists $\bar\gamma_n$ in between $\gamma$ and $\hat\gamma$ such that
\begin{align*}
    0=\tilde\Psi_n(\hat\gamma)
    =\tilde\Psi_n(\gamma)-h_1'(\bar\gamma_n)(\hat\gamma-\gamma).
\end{align*}
Thus, we have
\begin{align*}
    \sqrt{n}(\hat\gamma-\gamma)
    =\frac{1}{h_1'(\bar\gamma_n)}\frac{1}{\sqrt{n}}\sum_{i=1}^n\{Y_1-h_1(\gamma)\}.
\end{align*}
Therefore, $\hat\gamma\pconv\gamma$, the continuity of $h_1'$, and $h_1'(\gamma)\neq0$ imply the result.
}

\textbf{Case of \eqref{eq:empirical_equation}} By the definition, we have
\begin{align}
    \sup_{\varsigma,\varsigma_e}\|\Psi_n(\varsigma,\varsigma_e)-\Psi(\varsigma,\varsigma_e)\|
    =\norm{\rbr{n^{-1}\sum_{i=1}^n Y_i^2-\E[\bar{Y}^2],n^{-1}\sum_{i=1}^n Y_i^4-\E[\bar{Y}^4]}}\overset{\rm a.s.}{\longrightarrow}0,
\end{align}
where the convergence follows from the law of large numbers. This uniform convergence, Assumption (A5), and $\Psi(\gamma)=\bzero$ imply $(\hat{\gamma}^2,\hat\sigma_e^2)\overset{\rm a.s.}{\longrightarrow}(\gamma^2,\sigma_e^2)$ by Theorem 5.9 in \cite{van2000asymptotic}.

{\ccc
Let $\vartheta=(\gamma,\sigma_e)^\top$ and $\hat\vartheta=(\hat\gamma,\hat\sigma_e)^\top$.
Then, by Taylor expansion,
\begin{align*}
    \bzero_2
    =\Psi_n(\hat\vartheta)
    =\Psi_n(\vartheta)+A_n(\hat\vartheta-\vartheta),
\end{align*}
where
\begin{align*}
    A_n=\int_0^1D\Psi\{\vartheta+t(\hat\vartheta-\vartheta)\}{\rm d}t.
\end{align*}
By the consistency of $\hat\vartheta$ and the continuity of $D\Psi$, we have $A_n\pconv A$. 
Thus, since $A$ is invertible, we have
\begin{align*}
    \sqrt{n}(\hat\vartheta-\vartheta)
    &=-A_n^{-1}\sqrt{n}\Psi_n(\vartheta)\\
    &=-A^{-1}\frac{1}{\sqrt{n}}\sum_{i=1}^n\begin{pmatrix}
            Y_i^2-\E[\bar Y^2]\\
            Y_i^4-\E[\bar Y^4]
        \end{pmatrix}
        +o_p(1).
\end{align*}
}
 \end{proof}

\begin{proof}[Proof of Proposition \ref{prop:t-stat}] 
We have
\begin{align*}
    \sqrt{n}\hat\tau_j\frac{\hat{\beta}_j-\hat{\mu}\beta_j}{\hat{\sigma}}
    &=\sqrt{n}\hat\tau_j\frac{\hat{\beta}_j-\mu\beta_j}{{\sigma}}\frac{\sigma}{\hat{\sigma}}-\sqrt{n}\hat\tau_j\frac{(\mu-\hat{\mu})\beta_j}{\hat{\sigma}}\\
    &=\sqrt{n}\hat\tau_j\frac{\hat{\beta}_j-\mu\beta_j}{{\sigma}}-\sqrt{n}\hat\tau_j\frac{\hat{\beta}_j-\mu\beta_j}{{\sigma}}\frac{\hat{\sigma}-\sigma}{\hat{\sigma}}-\sqrt{n}\hat\tau_j\frac{(\mu-\hat{\mu})\beta_j}{\hat{\sigma}}\\
    &=\sqrt{n}\tau_j\frac{\hat{\beta}_j-\mu\beta_j}{{\sigma}}-o_p(1),
\end{align*}
where the last identity follows from Lemma \ref{lem:SE-conv}, $\sqrt{n}\tau_j\beta_j=O(1)$, and the assumption (A1). Finally, Proposition \ref{prop:asym_norm_betahat} concludes the proof of the first result.

{\cc
We next prove the uniform convergence across $j\in[p]$.
Define the oracle and plug-in statistics
\[
T_{j,n}:=\frac{\sqrt{n}(\hat\beta_j-\mu\beta_j)}{\sigma/\tau_j},
\qquad
\hat T_{j,n}:=\frac{\sqrt{n}(\hat\beta_j-\hat\mu\beta_j)}{\hat\sigma/\hat\tau_j}.
\]
A direct rearrangement yields the affine decomposition
\[
\hat T_{j,n}=a_{j,n} T_{j,n}-b_{j,n},
\qquad
a_{j,n}:=\frac{\hat\tau_j}{\tau_j}\frac{\sigma}{\hat\sigma},
\qquad
b_{j,n}:=\sqrt{n}\beta_j(\hat\mu-\mu)\frac{\hat\tau_j}{\hat\sigma}.
\]
Let $Q\sim\mathcal{N}(0,1)$ and denote the Kolmogorov distance by
$d_K(U,V):=\sup_{t\in\R}|\pr(U<t)-\pr(V<t)|$.
Since the map $x\mapsto a_{j,n}x-b_{j,n}$ is monotone for $a_{j,n}>0$, we have
\begin{align}
    d_K(\hat T_{j,n},Q)
&\le d_K\big(a_{j,n}T_{j,n}-b_{j,n}, a_{j,n}Q-b_{j,n}\big)
     +d_K\big(a_{j,n}Q-b_{j,n}, Q\big)\\
&= d_K(T_{j,n},Q) + d_K\big(\mathcal{N}(-b_{j,n},a_{j,n}^2), \mathcal{N}(0,1)\big).
\end{align}
Taking the maximum over $j\in[p]$ gives
\[
\max_{j\in[p]} d_K(\hat T_{j,n},Q)
\le
\max_{j\in[p]} d_K(T_{j,n},Q)
+
\max_{j\in[p]} d_K\big(\mathcal{N}(-b_{j,n},a_{j,n}^2), \mathcal{N}(0,1)\big).
\]
The first term converges to $0$ by Proposition~\ref{prop:asym_norm_betahat} (uniform version).

For the second term, note that for $|a-1|\le 1/2$, the triangle inequality and the mean-value theorem imply
\begin{align}
    d_K\big(\mathcal{N}(-b,a^2), \mathcal{N}(0,1)\big)
&\le
\sup_{t\in\R}\big|\Phi((t+b)/a)-\Phi(t/a)\big|
+
\sup_{t\in\R}\big|\Phi(t/a)-\Phi(t)\big|\\
&\le
\frac{|b|}{a\sqrt{2\pi}} + C|a-1|
\end{align}
for a universal constant $C>0$.
Hence,
\[
\max_{j\in[p]} d_K\big(\mathcal{N}(-b_{j,n},a_{j,n}^2), \mathcal{N}(0,1)\big)
\lesssim
\max_{j\in[p]}|b_{j,n}|
+
\max_{j\in[p]}|a_{j,n}-1|.
\]
By Lemma~\ref{lem:SE-conv} and the assumed uniform consistency of $\hat\tau_j$,
we have $\max_{j\in[p]}|a_{j,n}-1|=o_p(1)$.
Moreover,
\[
\max_{j\in[p]}|b_{j,n}|
\le
\frac{|\hat\mu-\mu|}{\hat\sigma} 
\max_{j\in[p]}\sqrt{n} |\beta_j| \hat\tau_j
= o_p(1),
\]
under $\hat\mu-\mu=o_p(1)$ and $\max_{j\in[p]}\sqrt{n}\tau_j|\beta_j|=O(1)$ (together with $\max_{j\in[p]}|\hat\tau_j/\tau_j-1|=o_p(1)$).
Combining the above bounds yields
\[
\max_{j\in[p]}\sup_{t\in\R}\Big|\pr\Big(\hat T_{j,n}<t\Big)-\Phi(t)\Big|
=
\max_{j\in[p]} d_K(\hat T_{j,n},Q)
\to 0,
\]
which completes the proof of the uniform convergence.
}
\end{proof}

\begin{proof}[Proof of Theorem \ref{thm:conf-intrvl}]
By Proposition \ref{prop:t-stat}, $\mathcal{CI}_{\alpha,j}$ is clearly asymptotically level $(1-\alpha)$ confidence interval, since we have
\begin{align*}
    &\mathrm{Pr}\rbr{\frac{\hat{\beta}_j}{\hat{\mu}}+z_{\alpha/2}\frac{\hat{\sigma}}{\sqrt{n}\hat{\mu}\hat{\tau}_j}\le\beta_j\le\frac{\hat{\beta}_j}{\hat{\mu}}-z_{\alpha/2}\frac{\hat{\sigma}}{\sqrt{n}\hat{\mu}\hat{\tau}_j}}\\
    &=\mathrm{Pr}\rbr{z_{\alpha/2}\le\frac{\sqrt{n}(\hat{\beta}_j-\hat{\mu}\beta_j)}{\hat{\sigma}/\hat{\tau}_j}\le z_{(1-\alpha/2)}}\\
    &\to(1-\alpha).
\end{align*}
\end{proof}

\section{Technical Lemmas}
For a matrix $A\in\R^{m\times m}$, $A\succeq0$ is defined to mean that $A$ is positive semi-definite.
\begin{lemma}
\label{lem:str-convex}
Suppose that $\bX_i\sim\mathcal{N}_p(\bzero,\bI_p)$. Under the assumption (A2), for some constant $\epsilon_0$ such that $0\le\epsilon\le\epsilon_0$, 
\begin{align*}
    \frac{1}{n}\sum_{i=1}^n\nabla^2\ell(\bb;\bX_i,Y_i)
    \succeq\rbr{\inf_{z:|z|\le\frac{3\norm{\bb}}{\sqrt{\epsilon}}}g'(z)}\rbr{\sqrt{1-\epsilon}-\sqrt{\kappa}-2\sqrt{\frac{H(\epsilon)}{1-\epsilon}}}^2\bI_p
\end{align*}
with $H(\epsilon)=-\epsilon\log\epsilon-(1-\epsilon)\log(1-\epsilon)$ holds for any $\bb\in\R^p$ with probability at least $1-2\exp(-nH(\epsilon))-2\exp(-n/2)$.
\end{lemma}
\begin{proof}[Proof of Lemma \ref{lem:str-convex}]
Since the proof of Lemma 3 in \cite{sur2019likelihood} only uses the specific structure $\nabla^2\ell(\bb;\bX_i,Y_i)=\rho''(\bX_i^\top\bb)\bX_i\bX_i^\top$ with $\rho'(t)=1/(1+\exp(-t))$ {\ccc and the monotonicity of $\rho'$}, we can immediately extend it to our general structure $\nabla^2\ell(\bb;\bX_i,Y_i)=g'(\bX_i^\top\bb)\bX_i\bX_i^\top$ {\ccc by the monotonicity of $g$}.
\end{proof}

\begin{remark}
    This lemma implies that, for sufficiently small $\epsilon>0$,
    \begin{align*}
        \frac{1}{n}\sum_{i=1}^n\nabla^2\ell(\bb;\bX_i,Y_i)\succeq\omega(\norm{\bb})\bI_p,
    \end{align*}
    for some non-increasing positive function $\omega:\R_+\to\R_+$ independent of $n$.
\end{remark}

{\ccc
\begin{lemma}
\label{lem:gamp-cauchy-unified}
Let $\lambda\ge 0$, and suppose that
$(\mu_\lambda,\sigma_\lambda,\eta_\lambda)$, with
$\sigma_\lambda,\eta_\lambda>0$, is a solution to the SE system \eqref{eq:SE-system}. Define $\ell_y(t):=G(t)-yt$
and
\[
\Psi_\lambda(u,y):=
u-{\prox}_{\eta_\lambda\ell_y}(u)
=\eta_\lambda
\left[
g\left(
{\prox}_{\eta_\lambda G}(u+\eta_\lambda y)
\right)-y
\right].
\]
Set $a_\lambda:=1-{\lambda\kappa^{-1}\eta_\lambda}.$
Consider the stationary GAMP recursion
\begin{align}
\bm b^{k+1}
&=-\frac{1}{\kappa}
\breve{\bm X}^{\top}
\Psi_\lambda(\bm\xi^k,\bm Y)
+a_\lambda\bm b^k,
\label{eq:stationary-gamp-b}\\
\bm\xi^{k+1}
&=\breve{\bm X}\bm b^{k+1}
+\Psi_\lambda(\bm\xi^k,\bm Y),
\label{eq:stationary-gamp-xi}
\end{align}
where $\breve{\bm X}=\bm X/\sqrt n$, and $\Psi_\lambda$ is applied
coordinatewise.
Suppose that the recursion is initialized by
\[
\bm b^0
=\mu_\lambda\sqrt n\,\bm\beta
+\sigma_\lambda\bm z^0,
\qquad
\bm z^0\sim\mathcal N_p(\bm0,\bm I_p),
\]
independently of $(\breve{\bm X},\bm\varepsilon)$, and
$\bm\xi^0=\breve{\bm X}\bm b^0$.
Under the assumptions of Lemma \ref{lem:master},
\begin{align}
\lim_{k\to\infty}
\lim_{n,p\to\infty}
\frac1p
\|\bm b^{k+1}-\bm b^k\|_2^2
&=0,
\label{eq:cauchy-b}\\
\lim_{k\to\infty}
\lim_{n,p\to\infty}
\frac1n
\|\bm\xi^{k+1}-\bm\xi^k\|_2^2
&=0
\label{eq:cauchy-xi}
\end{align}
almost surely.
\end{lemma}
\begin{proof}
This proof is inspired by Lemma 6.9 in \citep{donoho2016high}.
For brevity, write $(\mu,\sigma,\eta):=(\mu_\lambda,\sigma_\lambda,\eta_\lambda)$, $\Psi:=\Psi_\lambda$.
Let
\[
D(u,y)
:=
\operatorname{prox}_{\eta\ell_y}(u)
=
\operatorname{prox}_{\eta G}(u+\eta y).
\]
The first-order condition for the proximal operator gives
\[
D(u,y)+\eta g(D(u,y))=u+\eta y,
\]
and hence
\[
\Psi(u,y)
=
u-D(u,y)
=
\eta\{g(D(u,y))-y\}.
\]
Differentiating the first-order condition with respect to $u$ yields
\[
\partial_uD(u,y)
=
\frac{1}{1+\eta g'(D(u,y))},
\]
so that
\begin{align}
\partial_u\Psi(u,y)
=1-\partial_uD(u,y)
=\frac{\eta g'(D(u,y))}
     {1+\eta g'(D(u,y))}
\in[0,1).
\label{eq:Psi-derivative}
\end{align}

Let us set
\[
U=\mu Z+\sqrt{\kappa}\sigma Q,
\qquad
\bar Y=h(Z,\bar\varepsilon),
\]
where $Z\sim\mathcal N(0,\gamma^2)$,
$Q\sim\mathcal N(0,1)$, and
$\bar\varepsilon$ is independent of $(Z,Q)$.
The first and third state-evolution equations imply
\begin{align}
\mathbb E[\Psi(U,\bar Y)^2]
&=
\kappa^2\sigma^2,
\label{eq:Psi-second-moment}\\
\mathbb E[\partial_u\Psi(U,\bar Y)]
&=
\kappa-\lambda\eta.
\label{eq:Psi-derivative-mean}
\end{align}
In particular, we obtain
$
a_\lambda
=
1-\frac{\lambda\eta}{\kappa}
=
\frac{1}{\kappa}
\mathbb E[\partial_u\Psi(U,\bar Y)].
$

Next, by Theorem 4.2 in \citep{feng2022unifying}\footnote{Although the theorem is stated only for a single time point $t$, the paragraph following it notes that the result can be extended to iterates at arbitrary time points; indeed, the proof is carried out in this more general form.
}, for every fixed $s,t\ge0$, there exists $q_{s,t}\in[-1,1]$ such that
\begin{align}
\frac1p\left\langle\bm b^s-\mu\sqrt n\,\bm\beta,\,\bm b^t-\mu\sqrt n\,\bm\beta\right\rangle
&\as
\sigma^2q_{s,t},
\label{eq:joint-se-b}\\
\frac1n
\left\langle
\bm\xi^s-\mu\bm X\bm\beta,\,
\bm\xi^t-\mu\bm X\bm\beta
\right\rangle
&\as
\kappa\sigma^2q_{s,t},
\label{eq:joint-se-xi}
\end{align}
and $q_{t,t}=1$.

More explicitly, if $(Q_s,Q_t)$ is a centered bivariate Gaussian
vector with unit marginal variances and
$\mathbb E[Q_sQ_t]=q_{s,t}$, then the joint state evolution for
$(\xi^s,\xi^t)$ is represented by
\[
U_s=\mu Z+\sqrt{\kappa}\sigma Q_s,
\qquad
U_t=\mu Z+\sqrt{\kappa}\sigma Q_t.
\]
By the definition of the state evolution recursion (e.g., we can see in equation (10) in \citep{feng2022unifying}), we have
\begin{align}
\sigma^2q_{s+1,t+1}
=
\frac{1}{\kappa^2}
\mathbb E[
\Psi(U_s,\bar Y)\Psi(U_t,\bar Y)
].
\label{eq:q-recursion-general}
\end{align}
For $s=t$, \eqref{eq:q-recursion-general} reduces to
\eqref{eq:Psi-second-moment}.

For $q\in[0,1]$, let $(Q_1,Q_2)$ be jointly standard Gaussian with
$\mathbb E[Q_1Q_2]=q$, and define
\[
F(q):=\mathbb E[\Psi(\mu Z+\sqrt{\kappa}\sigma Q_1,\bar Y)
\Psi(\mu Z+\sqrt{\kappa}\sigma Q_2,\bar Y)],\qquad
H(q):=\frac{F(q)}{F(1)}.
\]
Then $q_{s+1,t+1}=H(q_{s,t}).$

Taking the limit of $\frac1p\left\langle\bm b^0-\mu\sqrt n\,\bm\beta,\,\bm b^1-\mu\sqrt n\,\bm\beta\right\rangle$ by using the initialization $\bm b^0
=\mu_\lambda\sqrt n\,\bm\beta
+\sigma_\lambda\bm z^0$, Stein's identity
conditional on $(Z,\bar\varepsilon)$ gives
\[
q_{0,1}=a_\lambda-\frac1\kappa
\mathbb E[\partial_u\Psi(U,\bar Y)]=0,
\]
where the second identity follows from \eqref{eq:Psi-derivative-mean}.
Hence, writing $q_k:=q_{k,k+1}$, we have
$
q_0=0$ and $q_{k+1}=H(q_k).$
To analyze $H$, condition on $(Z,\bar\varepsilon)$ and define
\[
\phi_{Z,\bar\varepsilon}(x)
:=
\Psi(
\mu Z+\sqrt{\kappa}\sigma x,
h(Z,\bar\varepsilon)
).
\]
By \eqref{eq:Psi-second-moment} and Fubini's theorem,
\[
\mathbb E_{Z,\bar\varepsilon}
\left[
\mathbb E_Q
\left\{
\phi_{Z,\bar\varepsilon}(Q)^2
\right\}
\right]
=
F(1)
<
\infty.
\]
Hence,
$\phi_{Z,\bar\varepsilon}\in L_2(\mathcal N(0,1))$
for almost every $(Z,\bar\varepsilon)$.

Let $\{\mathsf H_m\}_{m\ge0}$ denote the orthonormal Hermite polynomials,
\[
\mathsf H_m(x)
:=
\frac{\operatorname{He}_m(x)}{\sqrt{m!}},
\]
and define $c_m(Z,\bar\varepsilon):=\mathbb E_Q
\left[\phi_{Z,\bar\varepsilon}(Q)\mathsf H_m(Q)\right].$
Then
\[
\phi_{Z,\bar\varepsilon}(x)
=
\sum_{m=0}^{\infty}
c_m(Z,\bar\varepsilon)\mathsf H_m(x)
\]
in $L_2(\mathcal N(0,1))$.
If $(Q_1,Q_2)$ is jointly standard Gaussian with
$\mathbb E[Q_1Q_2]=q$, Mehler's formula in this normalization states that
\[
\mathbb E
\left[
\mathsf H_m(Q_1)\mathsf H_\ell(Q_2)
\right]
=
\bm 1\{m=\ell\}q^m.
\]
Therefore, conditioning on $(Z,\bar\varepsilon)$ and then taking
expectations gives
\[
F(q)
=
\sum_{m=0}^{\infty}A_mq^m,
\qquad
A_m
:=
\mathbb E[c_m(Z,\bar\varepsilon)^2]
\ge0.
\]
Thus $H$ maps $[0,1]$ into itself. Moreover, since
$1-q^m\le m(1-q)$ for $q\in[0,1]$,
\begin{align}
1-H(q)
&\le
\rho_\lambda(1-q),
\label{eq:H-contraction}
\end{align}
where
$
\rho_\lambda
:=
\frac{\sum_{m=1}^{\infty}mA_m}{F(1)}.
$

It remains to show that $\rho_\lambda$ is strictly smaller than one.
By the chain rule, for almost every $x\in\mathbb R$,
$\phi_{Z,\bar\varepsilon}'(x)
=\sqrt{\kappa}\sigma\partial_u\Psi
\left(\mu Z+\sqrt{\kappa}\sigma x,h(Z,\bar\varepsilon)\right).$
Since $0\le\partial_u\Psi<1$, it follows that $\left|
\phi_{Z,\bar\varepsilon}'(x)
\right|
\le
\sqrt{\kappa}\sigma,$
and therefore
$\phi_{Z,\bar\varepsilon}'\in L_2(\mathcal N(0,1))$.
For every $m\ge1$, Gaussian integration by parts and the fact that $Q\mathsf H_{m-1}(Q)-\mathsf H_{m-1}'(Q)=\sqrt m\,\mathsf H_m(Q)$ give
\begin{align}
\mathbb E_Q
\left[
\phi_{Z,\bar\varepsilon}'(Q)
\mathsf H_{m-1}(Q)
\right]
=
\sqrt m\,
\mathbb E_Q
\left[
\phi_{Z,\bar\varepsilon}(Q)
\mathsf H_m(Q)
\right]
=
\sqrt m\,c_m(Z,\bar\varepsilon).
\label{eq:Hermite-derivative-coefficient}
\end{align}
Using this, applying Parseval's identity to
$\phi_{Z,\bar\varepsilon}'$ gives
\[
\mathbb E_Q
\left[
\phi_{Z,\bar\varepsilon}'(Q)^2
\right]
=
\sum_{m=1}^{\infty}
m\,c_m(Z,\bar\varepsilon)^2.
\]
Taking expectations with respect to $(Z,\bar\varepsilon)$ and applying
Tonelli's theorem, we obtain
\[
\sum_{m=1}^{\infty}mA_m
=
\kappa\sigma^2
\mathbb E
\left[
\partial_u\Psi(U,\bar Y)^2
\right].
\]
Using $F(1)=\kappa^2\sigma^2$, we obtain
$
\rho_\lambda
=
\frac1\kappa
\mathbb E[
\partial_u\Psi(U,\bar Y)^2
].
$
Because $0\le\partial_u\Psi(U,\bar Y)<1$ almost surely and
$
\mathbb E[\partial_u\Psi(U,\bar Y)]
=
\kappa-\lambda\eta,
$
we have $\rho_\lambda<1.$
Consequently, \eqref{eq:H-contraction} gives
$
0\le1-q_{k+1}
\le
\rho_\lambda(1-q_k),
$
and hence we have
\[
1-q_k\le\rho_\lambda^k\to0.
\]

Finally, Theorem 4.2 in \citep{feng2022unifying} with the test function
$(x,y)\mapsto(x-y)^2$ gives
\begin{align}
\lim_{n,p\to\infty}
\frac1p
\|\bm b^{k+1}-\bm b^k\|_2^2
&=
2\sigma^2(1-q_k),\\
\lim_{n,p\to\infty}
\frac1n
\|\bm\xi^{k+1}-\bm\xi^k\|_2^2
&=
2\kappa\sigma^2(1-q_k).
\end{align}
Letting $k\to\infty$ proves
\eqref{eq:cauchy-b}--\eqref{eq:cauchy-xi}.
\end{proof}
}

{\ccc
\begin{lemma}
\label{lem:SE-conv}
Let $Q_1,Q_2\sim\mathcal{N}(0,1)$ be independent, and let
$\bar\varepsilon$ be independent of $(Q_1,Q_2)$. Suppose that
$g:\mathbb{R}\to\mathbb{R}$ is continuously differentiable and
nondecreasing, and let $G:\mathbb{R}\to\mathbb{R}$ satisfy $G'=g$.
Suppose further that
$h:\mathbb{R}^2\to\mathbb{R}$ is Lipschitz continuous in its first
argument and that $\mathbb{E}\left[h(0,\bar\varepsilon)^2\right]<\infty.$

Let $\gamma>0$, and let
$\mathcal I_\gamma\subset(0,\infty)$ be a compact interval such that $\gamma\in\operatorname{int}(\mathcal I_\gamma).$
Let $\Theta\subset\mathbb{R}\times(0,\infty)^2$ be compact. 
Define the state-evolution residual map $\mathcal F:
\Theta\times\mathcal I_\gamma
\rightarrow
\mathbb{R}^3$.
Let $\bm\vartheta_0:=(\mu,\sigma^2,\eta)\in\Theta$ be the unique solution in $\Theta$ to the oracle state-evolution system
\[
\mathcal F(\bm\vartheta,\gamma)=\bm 0.
\]
Let $\hat\gamma\ge0$ satisfies
\[
\hat\gamma^2\pconv\gamma^2,
\]
and let $\hat{\bm\vartheta}:=(\hat\mu,\hat\sigma^2,\hat\eta)$ be a measurable random vector satisfying
\begin{align}
\pr\left(
\hat\gamma\in\mathcal I_\gamma,\,
\hat{\bm\vartheta}\in\Theta
\right)
&\rightarrow1,
\label{eq:plugin-SE-compactness}\\
\left\|
\mathcal F(\hat{\bm\vartheta},\hat\gamma)
\right\|
&\pconv0.
\label{eq:plugin-SE-approx-root}
\end{align}
Then, we obtain
\[
\hat{\bm\vartheta}
\pconv
\bm\vartheta_0.
\]
\end{lemma}

\begin{proof}
For $\bm\vartheta=(m,v,a)\in\Theta$ and $\varsigma\in\mathcal I_\gamma$, we can write
\begin{align}
\mathcal F(\bm\vartheta,\varsigma)
:=
\begin{pmatrix}
\kappa^2v
-
a^2\mathbb{E}\left[
R_{\bm\vartheta,\varsigma}^2
\right]
\\[1mm]
\mathbb{E}\left[
\varsigma Q_1R_{\bm\vartheta,\varsigma}
\right]
\\[1mm]
1-\kappa
-
\mathbb{E}\left[
\left\{
1+a g'(D_{\bm\vartheta,\varsigma})
\right\}^{-1}
\right]
\end{pmatrix},
\label{eq:SE-residual-map}
\end{align}
where
\begin{align}
\bar Y_{\varsigma}
&:=
h(\varsigma Q_1,\bar\varepsilon),
\qquad S_{\bm\vartheta,\varsigma}
:=
m\varsigma Q_1+\sqrt{\kappa v}\,Q_2,
\label{eq:SE-linear-varsigma}\\
D_{\bm\vartheta,\varsigma}
&:=
{\prox}_{aG}
\left(
S_{\bm\vartheta,\varsigma}
+a\bar Y_{\varsigma}
\right),
\qquad R_{\bm\vartheta,\varsigma}
:=
\bar Y_{\varsigma}
-
g(D_{\bm\vartheta,\varsigma}).
\label{eq:SE-residual-varsigma}
\end{align}

\textbf{Step 1: Continuity}.
We first show that $(\bm\vartheta,\varsigma)
\mapsto
\mathcal F(\bm\vartheta,\varsigma)$
is continuous on
$\Theta\times\mathcal I_\gamma$.

Consider an arbitrary convergent sequence $(\bm\vartheta_\ell,\varsigma_\ell)
\rightarrow
(\bm\vartheta,\varsigma)$ in $\Theta\times\mathcal I_\gamma,$
where $\bm\vartheta_\ell=(m_\ell,v_\ell,a_\ell)$ and $\bm\vartheta=(m,v,a).$
For brevity, write
\[
\bar Y_\ell
:=
\bar Y_{\varsigma_\ell},
\qquad
S_\ell
:=
S_{\bm\vartheta_\ell,\varsigma_\ell},
\qquad
D_\ell
:=
D_{\bm\vartheta_\ell,\varsigma_\ell},
\qquad
R_\ell
:=
R_{\bm\vartheta_\ell,\varsigma_\ell},
\]
and define
$\bar Y$, $S$, $D$, and $R$ analogously at
$(\bm\vartheta,\varsigma)$.

Since $h$ is Lipschitz continuous in its first argument,
\[
\bar Y_\ell
=h(\varsigma_\ell Q_1,\bar\varepsilon)
\as
h(\varsigma Q_1,\bar\varepsilon)
=\bar Y.
\]
Moreover,
\[
S_\ell
=m_\ell\varsigma_\ell Q_1+\sqrt{\kappa v_\ell}\,Q_2
\as
m\varsigma Q_1+\sqrt{\kappa v}\,Q_2=
S.
\]

By the first-order condition for the proximal operator,
\begin{align}
D_\ell+a_\ell g(D_\ell)
&=
S_\ell+a_\ell\bar Y_\ell,
\label{eq:prox-FOC-seq}\\
D+a g(D)
&=
S+a\bar Y.
\label{eq:prox-FOC-limit}
\end{align}
Subtracting \eqref{eq:prox-FOC-limit} from
\eqref{eq:prox-FOC-seq} gives
\begin{align}
&D_\ell-D
+
a_\ell\{g(D_\ell)-g(D)\}
\notag\\
&\qquad=
S_\ell+a_\ell\bar Y_\ell
-
S-a\bar Y
+
(a-a_\ell)g(D).
\label{eq:prox-continuity-decomp}
\end{align}
Since $g$ is nondecreasing,
\[
(D_\ell-D)\{g(D_\ell)-g(D)\}\ge0.
\]
Consequently, the two terms on the left-hand side of
\eqref{eq:prox-continuity-decomp} have the same sign, and hence
\begin{align}
|D_\ell-D|
&\le
\left|
S_\ell+a_\ell\bar Y_\ell-S-a\bar Y\right|
+|a_\ell-a|\,|g(D)|
\as0.
\label{eq:prox-joint-continuity}
\end{align}
Also, it follows from the continuity of $g$ and $g'$ that
\begin{align}
R_\ell
&\as R,
\label{eq:residual-pointwise-conv}\\
\left\{
1+a_\ell g'(D_\ell)
\right\}^{-1}
&\as
\left\{
1+a g'(D)
\right\}^{-1}.
\label{eq:derivative-pointwise-conv}
\end{align}

We next construct integrable envelopes in order to pass the limits through
the expectations in \eqref{eq:SE-residual-map}. Since
$\Theta\times\mathcal I_\gamma$ is compact, there exist finite constants
$M,\bar v,\bar a,\bar\gamma>0$ and $\underline a>0$ such that
\[
|m|\le M,\qquad
v\le\bar v,\qquad
\underline a\le a\le\bar a,\qquad
\varsigma\le\bar\gamma
\]
for every
$(m,v,a)\in\Theta$ and
$\varsigma\in\mathcal I_\gamma$. Therefore, for a constant
$C_1>0$ independent of
$(\bm\vartheta,\varsigma)$,
\begin{align}
|S_{\bm\vartheta,\varsigma}|
\le
C_1(|Q_1|+|Q_2|).
\label{eq:linear-envelope}
\end{align}

For $a>0$, define $H_a(t):=t+a g(t).$
Since $g$ is nondecreasing, $H_a$ is strictly increasing and $|H_a(t)-H_a(s)|\ge |t-s|$ for all $s,t\in\mathbb{R}.$
Thus, $H_a^{-1}$ is $1$-Lipschitz. Since
\[
D_{\bm\vartheta,\varsigma}
=
H_a^{-1}\left(S_{\bm\vartheta,\varsigma}+a\bar Y_{\varsigma}\right)
\]
and
\[
H_a^{-1}(a g(0))=0,
\]
we obtain
\begin{align}
|D_{\bm\vartheta,\varsigma}|
&=|H_a^{-1}\left(S_{\bm\vartheta,\varsigma}+a\bar Y_{\varsigma}\right)-H_a^{-1}\left(ag(0)\right)|+|H_a^{-1}\left(ag(0)\right)| \\
&\le
\left|
S_{\bm\vartheta,\varsigma}+a\bar Y_{\varsigma}-a g(0)
\right|.
\label{eq:prox-envelope}
\end{align}
Furthermore, the proximal first-order condition implies
\[
D_{\bm\vartheta,\varsigma}-S_{\bm\vartheta,\varsigma}
=aR_{\bm\vartheta,\varsigma}.
\]
Combining this identity with \eqref{eq:linear-envelope} and
\eqref{eq:prox-envelope} yields, for a constant $C_2>0$,
\begin{align}
|R_{\bm\vartheta,\varsigma}|
\le
C_2
\left\{
1+|Q_1|+|Q_2|+|\bar Y_{\varsigma}|
\right\}
\label{eq:residual-envelope}
\end{align}
uniformly over
$(\bm\vartheta,\varsigma)
\in\Theta\times\mathcal I_\gamma$.

Let $L_h$ be a Lipschitz constant of $h$ with respect to its first
argument. Then
\begin{align}
\sup_{\varsigma\in\mathcal I_\gamma}
|\bar Y_{\varsigma}|
&\le
|h(0,\bar\varepsilon)|
+
L_h\bar\gamma |Q_1|.
\label{eq:Y-envelope}
\end{align}
By assumption, $\mathbb{E}\left[
h(0,\bar\varepsilon)^2
\right]<\infty.$
Hence, \eqref{eq:residual-envelope} and \eqref{eq:Y-envelope} imply that
\[
\sup_{\substack{
\bm\vartheta\in\Theta\\
\varsigma\in\mathcal I_\gamma
}}
|R_{\bm\vartheta,\varsigma}|^2
\]
is bounded by an integrable random variable. Similarly,
\[
\sup_{\substack{
\bm\vartheta\in\Theta\\
\varsigma\in\mathcal I_\gamma
}}
\left|
\varsigma Q_1R_{\bm\vartheta,\varsigma}
\right|
\]
is bounded by an integrable random variable, by the Cauchy--Schwarz
inequality. Finally, since $g'\ge0$,
\[
0<
\left\{
1+a g'(D_{\bm\vartheta,\varsigma})
\right\}^{-1}
\le1.
\]
Therefore, \eqref{eq:residual-pointwise-conv},
\eqref{eq:derivative-pointwise-conv}, and the dominated convergence theorem
show that every component of
$\mathcal F(\bm\vartheta,\varsigma)$ is continuous on
$\Theta\times\mathcal I_\gamma$.

\textbf{Step 2: Consistency}. 
We now prove the consistency of the plug-in solution. Since $\hat\gamma^2\pconv\gamma^2$, $\hat\gamma\ge0$, and $\gamma>0,$
the continuous mapping theorem gives
$\hat\gamma\pconv\gamma.$

Fix any $\epsilon>0$, and define
\[
\Theta_\epsilon
:=
\left\{
\bm\vartheta\in\Theta:
\|\bm\vartheta-\bm\vartheta_0\|\ge\epsilon
\right\}.
\]
If $\Theta_\epsilon=\varnothing$, the conclusion is immediate. Otherwise,
$\Theta_\epsilon$ is compact.

By assumption, $\bm\vartheta_0$ is the unique zero in $\Theta$ of $\bm\vartheta
\mapsto
\mathcal F(\bm\vartheta,\gamma).$
Hence, $\|\mathcal F(\bm\vartheta,\gamma)\|>0$ for every $\bm\vartheta\in\Theta_\epsilon.$
Since
$\bm\vartheta\mapsto
\|\mathcal F(\bm\vartheta,\gamma)\|$
is continuous and $\Theta_\epsilon$ is compact, it attains a strictly
positive minimum:
\begin{align}
c_\epsilon
:=
\inf_{\bm\vartheta\in\Theta_\epsilon}
\|\mathcal F(\bm\vartheta,\gamma)\|
>0.
\label{eq:SE-root-separation}
\end{align}

Because $\mathcal F$ is continuous on the compact set
$\Theta\times\mathcal I_\gamma$, it is uniformly continuous there.
Therefore, there exists $\delta_\epsilon^{\rm uc}>0$ such that
\[
|\varsigma-\gamma|<\delta_\epsilon
\quad\implies\quad
\sup_{\bm\vartheta\in\Theta}
\left\|
\mathcal F(\bm\vartheta,\varsigma)
-
\mathcal F(\bm\vartheta,\gamma)
\right\|
<
\frac{c_\epsilon}{3}.
\]
We choose positive $\delta_\epsilon\le\delta_\epsilon^{\rm uc}$ sufficiently small that $(\gamma-\delta_\epsilon,\gamma+\delta_\epsilon)
\subset\mathcal I_\gamma.$

Define the event
\begin{align}
A_{n,\epsilon}
:=
&
\left\{
\hat\gamma\in\mathcal I_\gamma
\right\}
\cap
\left\{
\hat{\bm\vartheta}\in\Theta
\right\}
\cap
\left\{
|\hat\gamma-\gamma|<\delta_\epsilon
\right\}
\cap
\left\{
\left\|
\mathcal F(\hat{\bm\vartheta},\hat\gamma)
\right\|
<
\frac{c_\epsilon}{3}
\right\}.
\label{eq:SE-consistency-event}
\end{align}
By \eqref{eq:plugin-SE-compactness},
\eqref{eq:plugin-SE-approx-root}, and
$\hat\gamma\pconv\gamma$, we have $\pr(A_{n,\epsilon})\rightarrow1.$

Suppose that $A_{n,\epsilon}$ occurs and that $\|\hat{\bm\vartheta}-\bm\vartheta_0\|
\ge\epsilon.$
Then we have
$\hat{\bm\vartheta}\in\Theta_\epsilon$, so that
\eqref{eq:SE-root-separation} gives
\[
c_\epsilon
\le
\left\|
\mathcal F(\hat{\bm\vartheta},\gamma)
\right\|.
\]
On the other hand, by the triangle inequality,
\begin{align*}
\left\|
\mathcal F(\hat{\bm\vartheta},\gamma)
\right\|
&\le
\left\|
\mathcal F(\hat{\bm\vartheta},\gamma)
-
\mathcal F(\hat{\bm\vartheta},\hat\gamma)
\right\|
+
\left\|
\mathcal F(\hat{\bm\vartheta},\hat\gamma)
\right\|\\
&<
\frac{c_\epsilon}{3}
+
\frac{c_\epsilon}{3}
<
c_\epsilon,
\end{align*}
which is a contradiction. Therefore, on $A_{n,\epsilon}$,
\[
\|\hat{\bm\vartheta}-\bm\vartheta_0\|<\epsilon.
\]
Consequently,
\[
\pr\left(
\|\hat{\bm\vartheta}-\bm\vartheta_0\|
\ge\epsilon
\right)
\le
\pr(A_{n,\epsilon}^c)
\rightarrow0.
\]
This completes the proof.
\end{proof}
}

\begin{lemma}
\label{lem:theta-rot}
Let $\hat\btheta$ be the estimator, i.e., surrogate loss minimizer, in a GLM with a true coefficient vector $\btheta$ and features drawn i.i.d. from $\mathcal{N}_p(\bzero,\bI_p)$. Define $(\mu_n,\sigma_n)$ as in \eqref{eq:SE-n}. Then, 
\begin{align*}
    \frac{\hat\btheta-\mu_n\btheta}{\sigma_n}
\end{align*}
is uniformly distributed on the unit sphere lying in $\btheta^\perp$.
\end{lemma}
\begin{proof}[Proof of Lemma \ref{lem:theta-rot}]
Define an orthogonal projection matrix $\bP_{\btheta}=\btheta\btheta^\top/\norm{\btheta}^2$ onto the line including $\btheta$, and an orthogonal projection matrix $\bP_{\btheta}^\perp=\bI_p-\bP_{\btheta}$ onto the orthogonal complement of the line including $\btheta$. Let $\bU\in\R^{p\times p}$ be any orthogonal matrix obeying $\bU\btheta=\btheta$, i.e. any rotation operator about $\btheta$. Then, since $\hat\btheta=\bP_{\btheta}\hat\btheta+\bP_{\btheta}^\perp\hat\btheta$, we have
\begin{align*}
    \bU\hat\btheta
    =\bU\bP_{\btheta}\hat\btheta+\bU\bP_{\btheta}^\perp\hat\btheta
    =\bP_{\btheta}\hat\btheta+\bU\bP_{\btheta}^\perp\hat\btheta.
\end{align*}
Using this, we obtain
\begin{align}
\label{eq:rotation}
    \frac{\bU\bP_{\btheta}^\perp\hat\btheta}{\|\bP_{\btheta}^\perp\hat\btheta\|}
    \overset{\rm d}{=}\frac{\bP_{\btheta}^\perp\hat\btheta}{\|\bP_{\btheta}^\perp\hat\btheta\|}
    =\frac{\hat\btheta-\mu_n\btheta}{\sigma_n},
\end{align}
where the first identity follows from the fact that $\bU\hat\btheta\overset{\rm d}{=}\hat\btheta$ since $\bU\hat\btheta$ is the estimator with a true coefficient $\bU\btheta=\btheta$ and features drawn i.i.d. from $\mathcal{N}(\bzero,\bI_p)$.
\eqref{eq:rotation} reveals that $(\hat\btheta-\mu_n\btheta)/\sigma_n$ is rotationally invariant about $\btheta$, lies in $\btheta^\perp$, and has a unit norm. These conclude the proof.
\end{proof}

\begin{lemma}
\label{lem:SE-approx}
    Suppose that $0<\|\bL^\top\hat\bbeta\|,\|\bL^\top\bbeta\|\le C$ for some constant $C>0$. 
    Let $\btheta,\hat\btheta$ be defined in \eqref{eq:theta} {\ccc and $\kappa_n=p/n$}. 
    If we define
    \begin{align}
        \mu_n=\frac{\hat\btheta^\top\btheta}{\norm{\btheta}^2},\qquad
    \sigma_n^2=\frac{1}{\kappa_{\ccc n}}\|\hat\btheta-\mu_n\btheta\|^2,
    \end{align}
    then, under assumptions (A1) and (A2), we have
    \begin{align*}
        \mu_n\as\mu,\qquad\sigma_n^2\as\sigma^2,
    \end{align*}
    as $n,p(n)\to\infty$ with $p(n)/n\to\kappa$.
\end{lemma}
\begin{proof}[Proof of Lemma \ref{lem:SE-approx}]
Since $\|\btheta\|$ is bounded in probability by an assumption, $\bU\btheta=(\norm{\btheta},\ldots,\norm\btheta)/\sqrt{p}$ with some orthogonal matrix $\bU\in\R^{p\times p}$ satisfies an assumption of Lemma \ref{lem:master} with $\Bar{\beta}\sim\delta_{\gamma/\sqrt{p}}$ by the fact that $\norm{\bU\btheta}^2=\norm{\btheta}^2=\bbeta^\top\bSigma\bbeta=\gamma^2$. Then applying Lemma \ref{lem:master} to $(\bU\btheta,\bU\hat\btheta)$ with $\psi(s,t)=st, t^2$ and considering their ratio gives
\begin{align*}
    \frac{\langle\bU\btheta,\bU\hat\btheta\rangle}{\norm{\bU\btheta}^2}=\mu_n\as\mu.
\end{align*}
Also, setting $\psi(s,t)=(s-\mu t)^2$ yields
\begin{align*}
    \frac{1}{\kappa}\|\bU\hat\btheta-\mu\bU\btheta\|^2=\frac{1}{\kappa}\|\hat\btheta-\mu\btheta\|^2\as\sigma^2.
\end{align*}
Thus, we obtain
\begin{align*}
    \sigma_n^2=\frac{1}{\kappa}\|\hat\btheta-\mu\btheta\|^2+2(\mu-\mu_n)\btheta^\top\hat\btheta-(\mu^2-\mu_n^2)\norm{\btheta}^2\as\sigma^2.
\end{align*}
\end{proof}

\begin{lemma}
\label{lem:lip-prox}
Define $\norm{f}_\infty:=\sup_{t\in\R}|f(t)|$ for a function $f:\R\to\R$. A proximal operator $\mathrm{prox}(x,f,b)\equiv\mathrm{prox}_{bF}(x)$ with $x\in\R$, $b>0$, and a monotone continuous function $f=F'$, is Lipschitz continuous with respect to
\begin{enumerate}
    \item[(i)] $x\in\R$ with constant $1$,
    \item[(ii)] a monotone continuous function $f:\R\rightarrow\R$ in terms of $L^\infty$ norm $\norm{\cdot}_\infty$ with constant $b>0$,
    \item[(iii)] and $b>0$ with constant $\norm{f}_\infty$.
\end{enumerate}
\end{lemma}
\begin{proof}[Proof of Lemma \ref{lem:lip-prox}] 
For any $x,y\in\R$ and fixed monotone $f(\cdot)$, suppose that $u=\text{prox}(x,f,b)$ and $v=\text{prox}(y,f,b)$ and $u>v$ without loss of generality. Note that $F(\cdot)$ is convex since its derivative $f(\cdot)$ is a monotonically increasing function. By the first-order condition, we have $x=u+bf(u)$ and $y=v+bf(v)$. Since $f(\cdot)$ is monotonically increasing, $0\le(u-v)(f(u)-f(v))$. Using this,
\begin{align*}
0\le(u-v)(bf(u)-bf(v))
=(u-v)(x-u-y+v)
=(u-v)(x-y)-(u-v)^2,
\end{align*}
by $b>0.$ This immediately implies $|u-v|\le|x-y|$, i.e., $\text{prox}(x,f,b)$ is 1-Lipschitz continuous with respect to the first argument.

Next, set any two continuous monotone functions $f(\cdot)$ and $g(\cdot)$. Suppose that $\text{prox}(x,f,b)\ge\text{prox}(x,g,b)$ with fixed $x\in\R$ without loss of generality. Then,
\begin{align}
\label{eq:lip-prox2}
\abs{{\text{prox}(x,f,b)-\text{prox}(x,g,b)}}
&={\text{prox}(x,f,b)-\text{prox}(x,g,b)}\nonumber\\
&=b\rbr{g(\text{prox}(x,g,b))-f(\text{prox}(x,f,b))}\nonumber\\
&\le b\rbr{g(\text{prox}(x,f,b))-f(\text{prox}(x,f,b))}\nonumber\\
&\le b\norm{f-g}_\infty,
\end{align}
where the second identity follows from the fact that $\text{prox}(x,f)=x-bf(\text{prox}(x,f,b))$, and the first inequality is from {\cc monotonicity} of $g(\cdot)$. This means $\text{prox}(x,f,b)$ is Lipschitz continuous for the second argument with constant $b>0$ in terms of $\norm{\cdot}_\infty$.

At last, using the fact that $\prox(x,g,b)=\prox(x,f,\alpha b)$ for $g(x)=\alpha f(x), \alpha>0$ by the definition of the proximal operator, we have, for any $b,b'>0$,
\begin{align*}
    \abs{\prox(x,f,b)-\prox(x,f,b')}
    &= \abs{\prox(x,f,b)-\prox\rbr{x,\frac{b'}{b}f,b}}\\
    &\le b\norm{f-\frac{b'}{b}f}_\infty\\
    &=\norm{bf-b'f}_\infty\\
    &=\abs{b-b'}\norm{f}_\infty,
\end{align*}
where the inequality follows from \eqref{eq:lip-prox2}.
\end{proof}

\begin{lemma}
\label{lem:prox-bounded}
Let $F:\R\to\R$ be a strictly convex function and $f=F'$. Suppose that $f(\cdot)$ takes bounded values on bounded domains. For any bounded $x\in\R$ and $c,c'>0$, we have
\begin{align*}
    \abs{\prox_{cF}(x)-\prox_{c'F}(x)}\le C\abs{c-c'},
\end{align*}
for some positive constant C.
\end{lemma}
\begin{proof}[Proof of Lemma \ref{lem:prox-bounded}] By the proofs of Lemma \ref{lem:lip-prox} (ii) and (iii), we can improve Lemma \ref{lem:lip-prox} (iii) as
\begin{align*}
    \abs{\prox_{cF}(x)-\prox_{c'F}(x)}
    \le \abs{f(\prox_{cF}(x))}\abs{c-c'}.
\end{align*}
Thus, we can complete the proof by showing that $\prox_{cF}(x)$ is bounded. Remind that the definition of the proximal operator is 
\begin{align*}
    \prox_{cF}(x)=\argmin_{z\in\R}\cbr{cF(z)+\frac{1}{2}(z-x)^2}.
\end{align*}
We denote the objective function as $H(z):=cF(z)+\frac{1}{2}(z-x)^2=H_1(z)+H_2(z)$. Obviously, $H_2$ has the minimum value at $x$. 

(i) For the case that a minimizer $\tilde{z}$ of $H_1$ is bounded.
Note that the minimizer is unique by the strict convexity of $H_1$.
Without loss of generality, suppose that $\tilde{z} \leq z$.
In this case, we have $H_1'(z') < 0$ and $H_2'(z')< 0$ for $z' < \tilde{z}$, and also have $H_1'(z') > 0$ and $H_2'(z') > 0$ for $z < z'$.
Hence, $H'(z') < 0$ holds for  $z' < \tilde{z}$ and $H'(z') > 0$ holds for  $z < z'$, thus $\prox_{cF}(x)  \notin [-\infty, \tilde{z}) \cup (z, \infty] $.
In contrast, $H'(z')$ may be both positive or negative for $z' \in [\tilde{z},z]$.
Hence, $\prox_{cF}(x) \in [\tilde{z},z]$ holds and thus bounded itself.

(ii) For the case when $H_1$ has a minimum at an unbounded point, such as $H_1(z)=ce^z$, we can also show that $\prox_{cF}(x)$ is bounded. In this case, we can assume that $H_1$ is monotonically increasing without loss of generality. For $z>x$, we have $H'(z)>0$ by $H_1'(z)>0$ and $H_2(z)>0$. Since $H_1'(z)$ is positive and monotonically increasing by the monotonicity and the convexity of $H_1(z)$, there exists a constant $\Tilde{C}<x$ such that $H_1'(z)=cf(z)<-H_2'(z)=-(z-x)$ for any $z<\Tilde{C}$. Putting together the results, we have $\prox_{cF}(x)\in[\Tilde{C},x]$.
\end{proof}

\bibliographystyle{alpha}
\bibliography{main}

\end{document}